\documentclass[10pt]{article}

\usepackage[a4paper, margin=1.2in]{geometry}
\usepackage[title]{appendix}

\usepackage{amsmath} 
\usepackage{enumitem} 
\usepackage{amssymb, amsthm} 
\usepackage{graphicx}

\usepackage[usenames, dvipsnames]{color}

         \newtheorem{theorem}{Theorem}[section]
	 \newtheorem{proposition}[theorem]{Proposition}
	 \newtheorem{lemma}[theorem]{Lemma}
         \newtheorem{corollary}[theorem]{Corollary}

\theoremstyle{definition}        
         \newtheorem{remark}[theorem]{Remark}

\numberwithin{equation}{section}

\newcommand \Abb {\mathbb{A}}
\newcommand \Cbb {\mathbb{C}}
\newcommand \Dbb {\mathbb{D}}
\newcommand \Nbb {\mathbb{N}}
\newcommand \Rbb {\mathbb{R}}
\newcommand \Zbb {\mathbb{Z}}
\newcommand \Tbb {\mathbb{T}}
\newcommand \Hbb {\mathbb{H}}

\newcommand \be {\begin{equation}}
\newcommand \ee {\end{equation}} 

\newcommand{\Ccont}{ C }

\def \h {\hat}
\def \t {\tilde} 

\def \v {\vec}

\newcommand \pa {\partial}
\newcommand \wh {\widehat}
\newcommand \wt {\widetilde}

\newcommand \al {\alpha}
\newcommand \de {\delta}
\newcommand \ep {\epsilon}

\newcommand \ga {\gamma}
\newcommand \Ga {\Gamma}

\newcommand \na {\nabla}

\newcommand \si {\sigma}
\newcommand \ka {\kappa}
\newcommand \om {\omega}
\newcommand \Om {\Omega}

\newcommand \Acal {\mathcal{A}}
\newcommand \Bcal {\mathcal{B}}
\newcommand \Ccal {\mathcal{C}}
\newcommand \Dcal {\mathcal{D}} 
\newcommand \Fcal {\mathcal{F}} 
\newcommand \Gcal {\mathcal{G}} 
\newcommand \Lcal {\mathcal{L}}
\newcommand \Scal {\mathcal{S}}

\newcommand \Ocal {\mathcal{O}}

\newcommand \Rcal {\mathcal{R}}

\title{Local wellposedness for the free boundary incompressible Euler equations  with interfaces that exhibit cusps and corners of nonconstant angle}
\author{Diego C\'ordoba, Alberto Enciso and Nastasia Grubic}
\date{\normalsize{Instituto de Ciencias Matem\'aticas\\ Consejo Superior de
  Investigaciones Cient\'\i ficas\\ 28049 Madrid, Spain\\[1ex]
  E-mail: dcg@icmat.es, aenciso@icmat.es, nastasia.grubic@icmat.es}}

\begin{document}

\maketitle

\begin{abstract}
We prove that free boundary incompressible Euler equations are locally well posed in a class
of solutions in which the interfaces can exhibit corners and cusps. Contrary to what happens in all the previously known non-$C^1$
water waves, the angle of these crests can change in time. 	
\end{abstract}

\tableofcontents


\section{Introduction}

Consider the motion of an inviscid incompressible irrotational fluid in
the plane with a free boundary. A time-dependent interface
$$
\Ga(t) := \{z(\al,t) \,|\, \alpha\in \Rbb\}
$$ 
separates the plane into two open sets: the water region, which we denote by
$\Omega(t)$, and the vacuum region,
$\Rbb^2\backslash\overline{\Omega(t)}$. 
The
evolution of the fluid is described by the Euler equations,
\begin{subequations}\label{equationsI}
	\begin{align}
		\pa_t v + (v \cdot \nabla )v = - \nabla p - ge_2\quad &\text{in} \quad \Omega(t),\\
		\nabla\cdot v=0  \quad \mathrm{and} \quad \nabla^{\perp}\cdot v=0  \quad &\text{in} \quad \Omega(t), \\
		(\pa_t z- v) \cdot (\pa_\al z)^\perp = 0 \quad &\mathrm{on} \quad \Ga(t),\\
		p=0\quad &\mathrm{on} \quad \Ga(t).
	\end{align} 
\end{subequations}
Here $v$ and $p$ are the water velocity and pressure
on~$\Omega(t)$, $e_2$ is the second vector of a Cartesian
basis and $g$ is the acceleration due to gravity. We are disregarding
the capillarity effects. It is standard that this system of PDEs
on~$\Rbb^2$, which are often referred to as the water waves equation,
can be formulated solely in terms of the interface curve, $z(\al,t)$,
and the vorticity density on the boundary,
defined through the formula
$$
\nabla^\perp \cdot v  =:  \om(\al,t)\, \delta (x - z(\al,t))\,.
$$

It is well known that the water waves system is locally well posed on
Sobolev spaces when the initial fluid configuration is sufficiently
smooth and the interface does not self-intersect. The first local
existence results for the free boundary incompressible Euler equations
are due to Nalimov \cite{Nalimov}, Yosihara \cite{Yosihara} and Craig
\cite{Craig} for near equilibrium initial data, and to
Wu~\cite{Wu:well-posedness-water-waves-2d,Wu:well-posedness-water-waves-3d}
for general initial data in Sobolev spaces. Local
wellposedness for initial data in low regularity Sobolev spaces was proven by
Alazard, Burq and Zuily~\cite{ABZ,ABZ2} and subsequently refined by
Hunter, Ifrim and Tataru~\cite{HIT}. For other variations and
results on local well-posedness, see \cite{Alazard-Burq-Zuily1,
	AM,Ch-Lindblad,CCG,Coutand-Shkoller:well-posedness-free-surface-incompressible,Lannes1,Lannes3,Lindblad,SZ1,CL,Zhang}. In
all these works, the
lowest regularity for the interfaces they consider is $C^{3/2}$, on which
the Rayleigh--Taylor stability condition $\partial_n p<0$ is assumed to hold.

When the initial
configuration is a suitably small perturbation of the stationary flat
interface, the system is in fact globally well
posed~\cite{IP,AlazardDelort,DIPP,Wu:global-wellposedness-3d,Germain2}. If the initial datum is not small, the equations can
develop splash singularities in finite time~\cite{ADCPJ2}. The two
essential features of this scenario of singularity formation (which remains valid in the case of
rotational fluids~\cite{CS1} and in the presence of viscosity~\cite{ADCPJ3,CS3}
or surface tension~\cite{ADCPJ4}) are that the velocity
and the interface remain smooth up to the singular time, and that the
self-intersecting interface does not pinch the water
region~\cite{FIL,CS4}. Stationary splash singularities that do pinch
the water region have been constructed in~\cite{CEG,CEG2}.

In this paper we are concerned with non-smooth interfaces that may present corners (thus preventing the interface
from being~$C^1$) or cusps. The study of this kind
of solutions hearkens back at least to Stokes, who formally
constructed traveling wave solutions which featured sharp crests with
a $120^\circ$ corner. In the 1980s, Amick, Fraenkel and
Toland~\cite{Amick} managed to rigorously establish the existence of
these solutions, and some 30 years later, under suitable technical
conditions Kobayashi~\cite{Ko2}
showed that these are in fact the only non-smooth traveling
waves. 

In a major recent work, Wu~\cite{W2} builds upon a priori energy estimates for the water waves system
previously derived with Kinsey~\cite{W,KW} to establish a local
existence result for a class of non-smooth initial data which allows
for interfaces featuring sharp crests (with any acute angle) or
cusps. This class includes the self-similar
solutions with angled crests she had
previously obtained in~\cite{Wu12}. Further study of the class of
singular solutions constructed by Wu was carried out by
Agrawal~\cite{SA}, who showed that these singularities are
``rigid''. More precisely, for this class of solutions, an initial
interface with an angled crest
remains angled crested and the angle does not change or
tilt. There are related rigidity results for cusped interfaces as well~\cite{SA}, and the effect of surface tension (which, in particular, makes it impossible to construct interfaces with angled crests in the energy class) has been studied in detail in~\cite{SA21,SA22}.

Our objective in this paper is to prove a local wellposedness result
for a wide class of initial data which allows for corners and
cusps and where these rigidity effects do not appear. We will work in the context of the 2D free boundary Euler equations, disregarding the gravity and capillarity effects. Physically, the
motivation is that, while one does expect to have sharp crests for
which the angle does not change, the direct observation of angled
crested waves in the ocean strongly suggests that there should also be
other fluid configurations where the angle changes in
time. As we will see later on, our main result rigorously establishes
this fact. From a mathematical point of view, one should observe that the
aforementioned rigidity results lay bare that a substantially
different approach to non-smooth water waves is required in order to
prove this result.

To make this precise, it is convenient to start by explaining the
problems that one must overcome to prove a local wellposedness
result for interfaces with sharp crests or low-regularity cusps. A
first issue is to understand what scales of weighted Sobolev
spaces can provide a good functional framework for this problem (and,
actually, if such weighted Sobolev spaces exist at all, which is not
obvious a priori). Once a choice of weighted spaces has been made, one
must construct an energy adapted to these spaces and show that one can
close the energy estimates. This presents two major difficulties. On the one hand, 
the non-smooth weights appearing in the energy become more and more
singular as one integrates by parts in the various integrals that
appear in the estimates, so there is no way to close the estimates
without a number of highly nontrivial cancellations. On the other
hand, the Rayleigh--Taylor condition fails at an angle point, so one can
only impose a degenerate stability condition of the form $\partial_n
p\leq 0$. To circumvent these difficulties, one needs to
start from a genuinely new basic idea and carry out the rather demanding
technical work necessary to implement it.

The basic idea underlying the approach to the motion of angle crested
interfaces developed in~\cite{KW,W,W2} is to map the singular
interface conformally to the half-space and control the regularity of
the interface through weighted norms of the conformal map. While the
local behavior of a conformal map from a wedge to the half-space then
gives a hint about the kind of weights one might use to define the
energy in this case, it is far from obvious, a priori, that one can
close the resulting energy estimates, and doing so is in fact a
technical tour de force. 

In contrast, the basic idea in our approach is to identify and control
a class of singular solutions where the vorticity density~$\om$ and a
certain number of its derivatives vanish at the singular point, at all
small enough times. The observation behind this philosophy is that
sufficiently smooth solutions to the equations do feature all these
zeros under suitable symmetry assumptions. In order to show that these
zeros exist and are preserved by the evolution for a certain class of
symmetric initial data with non-$C^1$ interfaces, and to effectively
use them to control the singular weights that appear in the energy
estimates, we carry out our analysis directly in the water region,
which is not smooth. The reason for which the singular interfaces
appearing in this class of solutions are not rigid is therefore that
we are not making any assumptions about the existence of a conformal
map from the water region to the half-space for which certain weighted
norms remain bounded.

From a technical standpoint, a drawback of this approach is that we can only employ real-variable methods in
all our key estimates. An upside of this is that, as we are
not using conformal maps in an essential way, these ideas
should carry over to three-dimensional problems and to the two-fluid
case. We will explore these and other directions in forthcoming
contributions; the gravity water waves problem will be considered too. Also, a technical point reflecting the differences in both
approaches is that the cusped interfaces that appear in our class of
solutions can be of $C^{2,\al}\backslash C^3$ H\"older regularity,
while those considered in~\cite{KW,W,W2,SA} are of class~$C^{3,\al}$.

In order to construct a suitable functional framework which allows for interfaces with angled crests and where one can close the energy estimates for the free boundary Euler system, we have built upon the work of Maz'ya and Soloviev about boundary value problems for the Laplacian on
domains with cusps~\cite{MazSo}. To avoid getting bogged down in technicalities at this stage, let us just say that we consider scales of Sobolev spaces $H^k_\beta(m)$ which involve power weights that vanish at the tip of cusp or corner, and that the strength of the weight depends on the geometry of the interface at the singular point. Both the position of the zero of the weight and its strength remain constant during the evolution of the fluid.

Let us now pass to state our main result, which ensures that
the free boundary Euler equations are well posed
within a class of initial data including interfaces with angled
crests (whose angle changes in time) and with cusps. In
terms of the aforementioned weighted Sobolev spaces $H^k_\beta(m)$, whose definition we prefer until later, this local existence result
can be informally stated as follows. Precise statements of this result in the case of interfaces that exhibit cusps or corners are presented below as Theorems~\ref{thm:T.main-corner} and~\ref{thm:T.main}, respectively.


\begin{theorem}\label{T.main}
	The 2D free boundary Euler equations, given by the system~\eqref{equationsI} with $g=0$, are locally well-posed in a suitable scale of weighted Sobolev spaces $H^k_\beta(m)$ that allows for interfaces with corners and cusps, provided that a suitable analog of the Rayleigh--Taylor stability condition holds.
\end{theorem}

\begin{remark}\label{L.angleintro}
	In Remark~\ref{R.opening} we obtain an explicit formula for the rate of change for the angle of the corner which shows, in particular, that the angle does indeed change for typical initial data with an angled crest.
\end{remark}

The paper is organized as follows.
Firstly, in Section~\ref{s.preliminaries} we
write the water waves problem as a system of equations for the
interface curve and the vorticity density. In Section~\ref{s.estimatesSingInt} we
present some estimates for singular integral operators on weighted
Lebesgue spaces that will be of use throughout the
paper. We focus on cusped interfaces, as the estimates are more complicated in that case. Sections~\ref{s.estEnergy} and~\ref{s.Corner} are respectively devoted to deriving
the essential a priori estimates for the water waves system with
non-smooth interfaces in the Lagrangian parametrization and to proving
our local wellposedness theorem. The proofs of several key results for boundary value problems on
domains with outer cusps (and corners), in the style of Maz’ya and Soloviev's results on singular integral operators in domains with cusps~\cite{MazSo}, are presented in Section~\ref{s:inverseOp}. To
streamline the presentation, the proofs of several important technical lemmas
are relegated to an Appendix.

\section{Preliminaries}
\label{s.preliminaries}

We consider the incompressible irrotational fluid flow in a fully symmetric bounded planar domain $\Om = \Om(t)$ governed by Euler equations. More precisely, the fluid velocity $v$ and the pressure $P$ satisfy
\begin{subequations}
\begin{align}
\pa_t v + (v \cdot \nabla) v = - \frac{1}{\varrho}\nabla P,  \label{eq:dP} \\
\nabla\cdot v=0,  \quad \nabla^{\perp}\cdot v=0,  \label{eq:Inc-Irr}
\end{align} 
\end{subequations}
in $\Om$, where the fluid density $\varrho>0$ is assumed to be constant. For simplicity, we take $\varrho \equiv 1$. The interface $\Ga = \pa\Om$ is a closed curve characterized by the condition
\be\label{ass:P}
P\big|_{\Ga} = 0,
\ee
which corresponds to setting $\varrho\equiv 0$ in the exterior domain $\Rbb^2 \setminus \Om$. Moreover, the parametrization of the interface 
$$
\Ga = \{z(\al,t) \in \Rbb^2 \,|\, \alpha\in [-\pi, \pi]\}
$$
satisfies the kinematic boundary condition, i.e. 
\be
(\pa_t z- v) \cdot (\pa_\al z)^\perp = 0 \label{eq:evol_zt}.
\ee
For all times $t \geq 0$, the domain $\Om$ is assumed to be a union of two simply connected, disjoint, bounded domains each with a curvilinear corner of opening $2\nu > 0$ or an outward cusp (corresponding to $2\nu = 0$), connected through a common tip situated at the origin. We assume $\Om$ is symmetric with respect to both axes, which implies the intersection point remains at the origin for all times.
	\begin{figure}[h]\centering
		\vspace{-0.5cm}
	\includegraphics[scale=0.2]{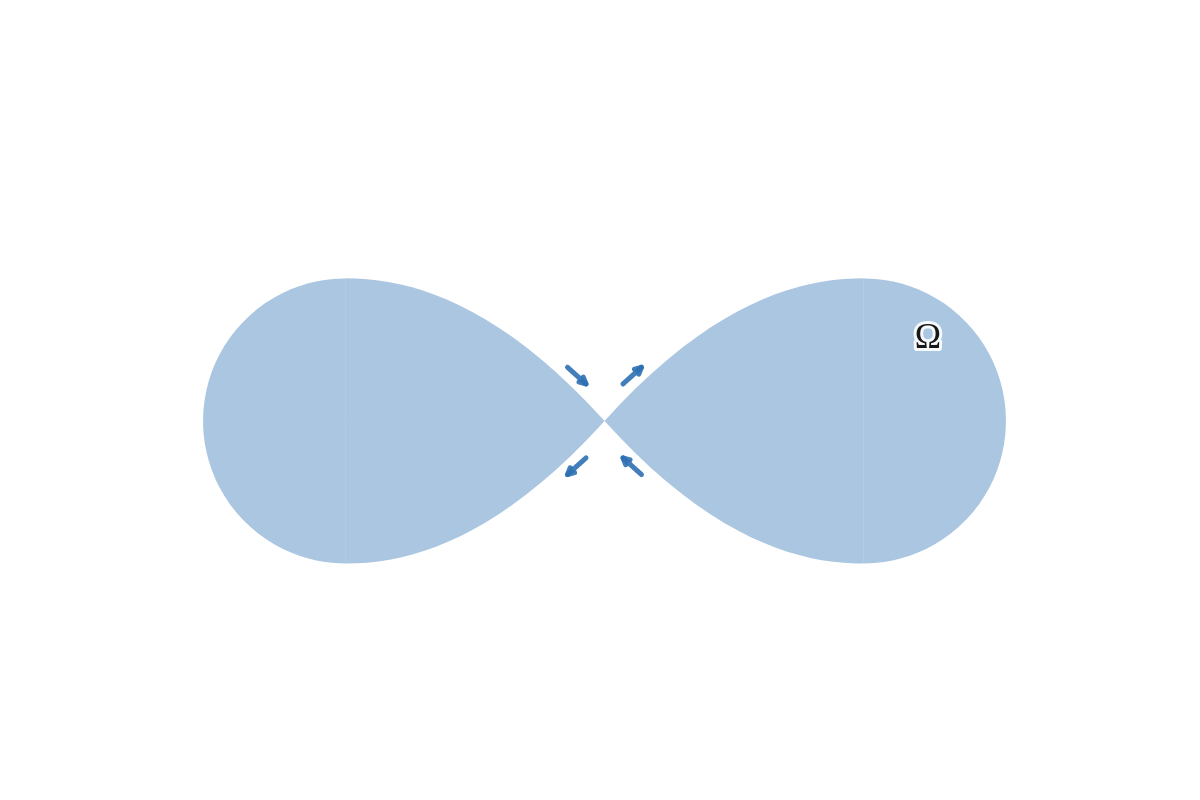}
		\vspace{-1cm}
	\caption{Fluid domain $\Om$ when $2\nu > 0$. The arrows specify the orientation.}
	\label{drop}
\end{figure}

Stated in terms of the parametrization (and letting $\theta(\al, t)$ denote the tangent angle at the corresponding point $z(\al, t)\in \Ga$), the intersection point in case of e.g. an outward cusp is characterized by
\be\label{eq:splash}
0 = z( \al_*, t)= z(-\al_*, t), \qquad 0 = \theta(\al_*, t) = \pi - \theta(-\al_*, t), \qquad \al_* = \pi/2, 
\ee   
but $\Ga$ is otherwise an arc-chord curve; i.e for any small $r>0$, we have 
\be\label{arcchord}
\Fcal_{r}(z) := \sup_{(\al, \beta)\in (B_r(\al_*)\times B_r(-\al_*))^c}\Fcal(z)(\al, \beta, t) < \infty,
\ee
where 
\be\label{eq:defFcal}
\Fcal(z)(\al, \beta, t) :=\begin{cases}
                        \frac{|e^{i\al} - e^{i\beta}|}{|z(\al, t)-z(\beta, t)|} \quad &\al\neq \beta\\
                        \frac{1}{|z_\al(\al, t)|}                       \quad &\al = \beta.
                       \end{cases}  
\ee

We consider the equations in the vorticity formulation. By assumption, the flow in $\Rbb^2\setminus \Ga$ is irrotational and we may assume the vorticity is a measure supported on $\Ga$, i.e.
$$
\om(z, t) = \om(\al, t) \de(z - z(\al, t)), \quad z\in \Rbb^2\setminus \Ga,
$$
where we slightly abuse the notation and denote both the vorticity and its amplitude by $\om$. Equations \eqref{eq:Inc-Irr} imply the complex conjugate of the velocity is  analytic in $\Om$, hence it can be written in terms of the boundary vorticity as  
\be\label{def:vom}
 v(z, t)^* = \frac{1}{2\pi i}\int_{-\pi}^\pi \frac{\om(\al', t)}{z-z(\al', t)} \,ds_{\al'}, \quad z\in \Rbb^2 \setminus \Ga,
\ee
where we have set $ds_{\al'} = |z_{\al}(\al', t)|d\al'$ and we use $^*$ to denote complex conjugation.  Approaching any regular point on $\Ga$ (in our setting, all points excluding the cusp/corner tip) from the inside of $\Om$ we obtain
$$
v(\al, t)^* = \frac{\om(\al, t)}{2 z_s(\al, t)} \,  + BR(z,\om)^*,
$$
where $z_s(\al, t)$ is the parametrization independent derivative of $z(\al, t)$ (i.e. $z_s(\al, t) \equiv \pa_s z(\al, t) = \frac{z_\al(\al, t)}{|z_\al(\al, t)|}$) and $BR(z,\om)$ is the Birkhoff-Rott integral whose complex-conjugate is defined via
$$
BR(z,\om)^* = \frac{1}{2\pi i} \, p.v. \int_{-\pi}^\pi \frac{\om(\al', t)}{z(\al, t)-z(\al', t)} \,ds_{\al'}
$$
where $p.v.$ stands for the principal value. 

We will frequently identify $z = (z_1, z_2)\in \Rbb^2$ with its complex representation $z  =z_1 + iz_2\in \Cbb$. In particular, we will use that the real scalar product $z\cdot w$ of two vectors $z, w \in\Rbb^2$, can be written as a product of two complex numbers, namely   
$$
z\cdot w = \Re \left(z^* w\right)
$$
Moreover, in complex notation we have $z^\perp = i z$, where $z^\perp$ is vector perpendicular to $z$.

Let us finish with some symmetry considerations. In terms of the parametrization of the interface $z(\al, t) = (z_1(\al, t),z_2(\al, t))$, the assumption of full symmetry implies
\be
\label{def:sym-x}
z_1(\al, t) = z_1(-\al, t), \qquad z_2(\al, t) = -z_2(-\al, t),
\ee
and
\be
\label{def:sym-y}
z_1(\al_* - \al, t) = - z_1(\al_*+ \al, t), \qquad z_2(\al_* - \al, t) = z_2(\al_* + \al, t).
\ee
The tangent angle $\theta = \arctan\big(\frac{z_{2\al}}{z_{1\al}}\big)$ then satisfies   
\be\label{def:sym-xy-theta}
\theta(\al, t) = \pi - \theta(-\al, t), \qquad \theta(\al_* - \al, t) = - \theta(\al_*+ \al, t).
\ee
On the other hand, the pressure must be invariant under both $z \leftrightarrow\pm  z^*$ and therefore 
$$
v(\pm z^*, t) = \pm v(z, t)^*
$$
by the Euler equations. In particular, the vorticity must be odd with respect to both axes, i.e.
\be\label{def:sym-om}
\om(-\al, t) = -\om(\al, t), \qquad \om(\al_* - \al,t) = -\om(\al_* + \al, t).
\ee

\subsection{Equations in vorticity formulation}

We now briefly state the relevant equations reformulated in terms of boundary vorticity $\om$ and parametrization of the interface $z$.  Equation \eqref{eq:evol_zt} gives 
\be\label{eq:z_t}
z_t^* = \left(\frac{\om}{2} - \t\varphi\right)\frac{1}{z_s} + BR(z,\om)^*, 
\ee
where $\t\varphi$ is a scalar function reflecting the freedom to choose the tangential component of $z_t$ and we use the notation $\pa_s := \frac{1}{|z_\al|}\pa_\al$ to denote the parametrization independent derivative. It will be convenient to work with the tangent angle $\theta$. Taking a time derivative of $z_\al = |z_\al|e^{i\theta}$, we obtain
$$
\theta_t = \frac{z_{t\al}\cdot z_\al^\perp}{|z_\al|^2} = z_{ts}\cdot z_s^\perp = \Re(iz_{ts}^* z_s) ,\quad \quad \frac{|z_{\al}|_t}{|z_\al|} = \frac{z_{t\al}\cdot z_\al}{|z_\al|^2} = z_{ts}\cdot z_s = \Re(z_{ts}^* z_s).
$$
Taking a $\pa_s$-derivative of equation \eqref{eq:z_t} then implies 
$$
\theta_t + \t\varphi\theta_s = \frac{\om\theta_s}{2} + BR(z,\om)_s\cdot z_s^\perp, \quad \quad \frac{|z_{\al}|_t}{|z_\al|} + \t\varphi_s = \frac{\om_s}{2}  + BR(z,\om)_s\cdot z_s.
$$
Since the pressure is constant on the interface, the tangential component of the pressure gradient \eqref{eq:dP} must be identically equal to zero, which combined with \eqref{eq:z_t} implies 
\be\label{eq:om_t}
\frac{\om_t}{2} + BR(z,\om)_t \cdot z_s = - \t\varphi \Big(\frac{\om_s}{2} + BR(z, \om)_s \cdot z_s \Big) = - \t\varphi \Big(\t\varphi_s + \frac{|z_\al|_t}{|z_\al|}\Big), 
\ee
while the equation for the normal component of the pressure gradient 
$$
\si:=-\na P \cdot z_s^\perp,
$$ 
reads  
\be\label{eq:si}
\si = \left(\frac{\om\theta_t}{2} + BR(z,\om)_t\cdot z_s^\perp\right) + \t\varphi \left(\frac{\om\theta_s}{2} + BR(z,\om)_s\cdot z_s^\perp\right)
\ee
The equation governing time evolution of $\t\varphi$ is now a matter of straightforward calculation. We have 
\be\label{eq:phi_tal}
\Big(\t\varphi_{s} + \frac{|z_\al|_t}{|z_\al|}\Big)_t =- \t\varphi \pa_s\Big(\t\varphi_s + \frac{|z_\al|_t}{|z_\al|}\Big) - \Big(\t\varphi_s + \frac{|z_\al|_t}{|z_\al|}\Big)^2 - \si \theta_s + \Big(BR_s\cdot z_s^\perp + \frac{\om\theta_s}{2}\Big)^2.   
\ee
At this point we fix the parametrization. We will consider the equations in the Lagrangian parametrization, i.e. we set $\t\varphi \equiv 0$ which means $z_t = v$. For later use, we introduce the notation
\be\label{eq:varphis}
\varphi_s := \frac{\om_s}{2} + BR(z,\om)_s\cdot z_s.
\ee
By the above, we have $\varphi_s = \frac{|z_\al|_t}{|z_\al|}$ which here takes the role $\t \varphi_s$ has in the arc-length parametrization (in which case $\varphi_s$ is a function depending only on time). To keep analogy with that case, we use the notation $\varphi_s$  although strictly speaking it is not a (parametrization-independent) derivative of any meaningful quantity. Except in this case subscript $s$ will always mean parametrization independent derivative.

\subsection{Weighted Sobolev spaces}
\label{ss.weightedSobolev}

Let us first introduce some notation. Given any two (non-negative) quantities $f$ and $g$, we say
$$
f \lesssim g \quad :\Leftrightarrow \quad \exists c>0: \ f \leq c g
$$
and similarly
$$
f \sim g\, \quad :\Leftrightarrow \quad \exists c>0: \ c^{-1} g \leq f \leq c g,
$$ 
where $c$ is either some absolute constant or can be controlled by some power of the energy. We will sometimes use the big-O notation, that is, we will write $f = O(g)$ if and only if $|f|\lesssim |g|$.

Let $\ga\in \Rbb$. We define the weighted Lebesgue space $\Lcal_{2, \ga}(\Ga, m)$ to be
$$
\Lcal_{2, \ga}(\Ga, m) := \{ \phi: \Ga \rightarrow \Rbb \, \, | \, \, m^\ga \phi \in L^2(\Ga) \}
$$ 
endowed with the norm 
$$
\|\phi \|_{2, \ga}^2 := \int_{\Ga} m(z)^{2\ga} |\phi(z)|^2 ds_z,
$$
where if not explicitly stated otherwise $m$ is always a power weight, i.e. we have  
$$
m(z) \sim |z|,  \quad z\in \Ga.
$$
If $0<\ga + 1/2 < 1$, we say that $m(\cdot)^\ga$ is a Muckenhaupt weight (it then satisfies the Muckenhaupt $(A_2)$-condition). Frequently we will consider weighted Lebesgue spaces on the one dimensional torus $\Tbb$ or on a particular interval $I\subseteq \Rbb$ instead of $\Ga$, however to simplify the notation we usually drop the space reference altogether and simply write $\Lcal_{2, \ga}(m)$. On the other hand, we may occasionally write $\Lcal_{2,\ga}(I)$ if we want to emphasize the particular interval of integration. 

We introduce two families of weighted Sobolev spaces:
$$
\aligned
\phi \in H^k_{\ga}(m) \quad &:\Leftrightarrow \quad \pa_s^j\phi\in \Lcal_{2, \ga}(m) , &\quad 0\leq j\leq k,\\
\phi \in \Lcal^k_{2, \ga}(m) \quad &:\Leftrightarrow \quad \pa_s^j\phi\in \Lcal_{2, \ga + (j - k)}(m), &\quad 0 \leq j\leq k,
\endaligned
$$
where $\pa_s = \frac{1}{|z_\al|}\pa_\al$ is the parametrization independent derivative. We will also need the subspace
$$
H^k_{\ga, 0}(m) := \{ \phi \in H^k_{\ga}(m) \, \, | \, \, m(z)^{-1}\phi(z) \in L^\infty \}.
$$
We clearly have 
\be\label{def:relsobolev}
\Lcal^k_{2, \ga}(m)  \subseteq H^k_{\ga}(m).
\ee
Hardy inequalities imply these can be identified, whenever $(\ga - k) + 1/2 > 0$ (for the convenience of the reader, we give these in the Appendix, Theorem \ref{thm:hardy} and Lemma \ref{lem:hardy}). The inclusion is proper otherwise. We give more details at the end of this section. 

We will frequently use the following: 
\begin{lemma}\label{lem:sobolev}
Let $m'= O(1)$. Then, we have  
$$
\phi \in \Lcal^1_{2,\ga + 1}(m)  \quad \Rightarrow \quad \phi = O\big(m^{-(\ga + \frac{1}{2})}\big). 
$$
\end{lemma}

\begin{proof}
The claim follows by integration since
$$
(m^{2(\ga + 1/2)} \phi^2)' = (2\ga + 1)m' m^{2\ga} \phi^2 + m^{2\ga + 1} \phi \phi' \in L^1.
$$
\end{proof}

To define fractional Sobolev spaces, let $\Lambda^{1/2}\phi := (-\Delta)^{1/4}\phi$ be the periodic fractional Laplacian defined (modulo a multiplicative constant factor) by  
$$
\Lambda^{1/2}\phi\,(\al) \, \sim \, p.v.\int_{-\pi}^\pi \frac{\phi(\al) - \phi(\al')}{|\sin \big(\frac{\al - \al'}{2}\big) |^{3/2}}\,d\al'.
$$
The weighted estimate for the Riesz integral suggests
$$
\pa_s\phi\in \Lcal_{2,\ga + \frac{1}{2}}(m)  \quad \Rightarrow \quad \Lambda^{1/2}\phi \in \Lcal_{2,\ga}(m)
$$
whenever $0<\ga + 1/2< 1/2$ (we prove this in the Appendix, Lemma \ref{lem:halfDer}). When $\ga$ does not satisfy this condition, we introduce a parameter $\lambda\in\Rbb$ such that   
$$
0 < (\ga - \lambda) + 1/2 < 1/2
$$
and define
$$
H^{k+1/2}_{\ga}(m) := \{\phi \in H^k_{\ga-1/2}(m) \, \, | \, \, \Lambda^{1/2}(m^\lambda\pa_s^k\phi) \in \Lcal_{2, \ga - \lambda}(m)\},
$$
endowed with the norm 
$$
\|\phi\|^2_{H^{k+1/2}_{\ga}(m)} := \|\phi\|^2_{H^{k}_{\ga-1/2}(m)} + \|\Lambda^{1/2}(m^\lambda\phi)\|^2_{2, \ga - \lambda}.
$$
Similarly, we define $\Lcal^{k+1/2}_{2,\ga}(m)$ to be the subspace of $\Lcal^{k}_{2,\ga-1/2}(m)$ such that $\Lambda^{1/2}(m^\lambda\pa_s^k\phi) \in \Lcal_{2, \ga - \lambda}(m)$. It is not difficult to see this definition is independent of the exact value of $\lambda$ (see the Appendix for similar results regarding commutators with $\Lambda^{1/2}$). 
Finally, we introduce the periodic Hilbert transform
$$
H\phi(\al) = \frac{1}{2\pi } p.v.\int_{-\pi}^\pi \phi(\al') \cot\Big(\frac{\al - \al'}{2}\Big)d\al'.
$$

When doing estimates on the singular integrals, it will be convenient to work in the graph parametrization. More precisely, we assume there exists a neighborhood $B$ of the origin such that $B\cap \Ga$ consists of exactly two connected components $\Ga^\pm$ which can be parametrized as a graph, that is
\be\label{def:Ga_pm}
B\cap \Ga = \Ga^+\cup \Ga^-, \qquad \Ga^\pm =\big\{x \pm i\ka(x, t) \,:\, |x|< 2\de \big\}, 
\ee 
where $\de = \de(t)>0$ (note the orientation is reversed on the lower branch). For the weight function, this corresponds to setting 
$$
m(z) \sim |z_1|,
$$
in which case, we use the notation 
$$
x = z_1(\al, t), \qquad m(x) = |x|, \qquad x\in I_{2\de} := (-2\de,2\de).
$$

To finish this section let us comment on the relation \eqref{def:relsobolev}, when $(\ga - k) + 1/2 < 0$. Let $I^+_\de=(0,\de)$ and let e.g. $k=1$. Then, using integration and Hardy inequalities (cf. Lemma \ref{lem:hardy}) we have for e.g. the right-hand side $\Ga^\pm \cap I^+_\de$  
$$
H^1_{\ga}(\Ga^\pm \cap I^+_\de) =\Lcal^1_{2,\ga}(\Ga^\pm \cap I^+_\de)\oplus \Rbb, \qquad \ga + 1/2 < 1
$$
and similarly for the left-hand side $\Ga^\pm \cap I^-_\de$, where we have set $I^-_\de := (-\de,0)$. When $\Ga^\pm$ are smooth enough (corresponding to two smoothly connected cusps, i.e.~$\nu = 0$), the constant must be the same, i.e. we have 
$$
H^1_{\ga}(\Ga^\pm) =\Lcal^1_{2,\ga}(\Ga^\pm)\oplus \Rbb.
$$
However, if $\Ga^\pm$ are only piecewise smooth (corresponding to the case $\nu > 0$), then $H^1_\ga(\Ga^\pm)$ is embedded in the space of piecewise continuous functions only, allowing for jumps when crossing the singular point. The generalization to higher $k$ is straightforward.



\section{Singular integrals on domains with cusps}
\label{s.estimatesSingInt}

Throughout this section we assume $\nu = 0$. At this point, let us discuss the regularity assumptions on the parametrization of the interface, which are assumed to hold throughout this section. Let $\mu\in (1/2, 1]$ be fixed and let $\beta\in\Rbb$ be such that 
\be\label{ass:beta-mu}
1 - \mu < \beta + 1/2 < \mu. 
\ee
In the neighborhood of the origin, cf.~\eqref{def:Ga_pm}, we assume 
\be\label{assump:rho''}
\ka(\cdot, t) \in \Lcal^4_{2, \beta + 2}(I_{2\de}), \qquad  |x|^{j - (1 + \mu)}|\pa^j_x\ka(x, t)| \, \lesssim \, 1, \qquad j = 0,1,2,
\ee
together with the lower bound 
\be\label{assump:lower-ka}
|x|^{- (1 + \mu)}\ka(x, t) \, \gtrsim \, 1.
\ee
Away from the origin the parametrization of the interface is $H^4$ and satisfies the arc-chord condition \eqref{arcchord}. Note that the lower bound in \eqref{ass:beta-mu} ensures the two assumptions on the derivatives of $\ka$ are consistent, since 
$$
1 - \mu < \beta + 1/2 \quad \Rightarrow \quad m^{\mu - 1}\in \Lcal_{2,\beta}(m),
$$
cf. Lemma \ref{lem:sobolev}. 

In view of the energy estimates in Section \ref{s.estEnergy}, we give conditions on the tangent angle $\theta$ and the length $|z_\al|$ of the tangent vector, which imply the above assumptions 
\be\label{regTheta}
\theta \in H^{3}_{\beta + 2}(\Tbb, m), \qquad \log|z_\al| \in H^2_{\beta + 1}(\Tbb, m),
\ee
In particular, we have $\log|z_\al|\in\Ccal^{0,\lambda'}(\Tbb)$ for some $\lambda'\in (0,1)$ and therefore $|z_\al| \sim 1$. Moreover, note that
$$
z_{s} = e^{i\theta} \quad \Rightarrow \quad  z_{1s} \in H^3_{\beta + 2}(m) , \quad z_{2s}  \in \Lcal^3_{2,\beta + 2}(m). 
$$
where we have taken into account $\theta(\al_*) = 0$ (resp. $\theta(-\al_*) = \pi$). When $\al$ is sufficiently close to $\al_*$, we further assume
\be\label{regTheta2}
m(\al)^{- \mu}|\theta(\al)| \, \sim \, 1, \qquad m(\al)^{1 - \mu}|\theta_\al(\al)| \, \lesssim \, 1.
\ee
Going over to the graph parametrization in the neighborhood of $\pm\al_*$, we have $z_{1s} = \frac{1}{|z_x|}$ resp. $\theta_s = \frac{\theta_x}{|z_x|} = \frac{\ka_{xx}}{|z_x|^3}$. By making $\de$ smaller we may assume $\ka_{xx}> 0$ on $I_{2\de}$ and therefore $\ka_x>0$ when $x>0$ (recall that $x=0$ is a local minimum for the upper branch).  

\begin{figure}[h]\centering	
	\includegraphics[scale=0.18]{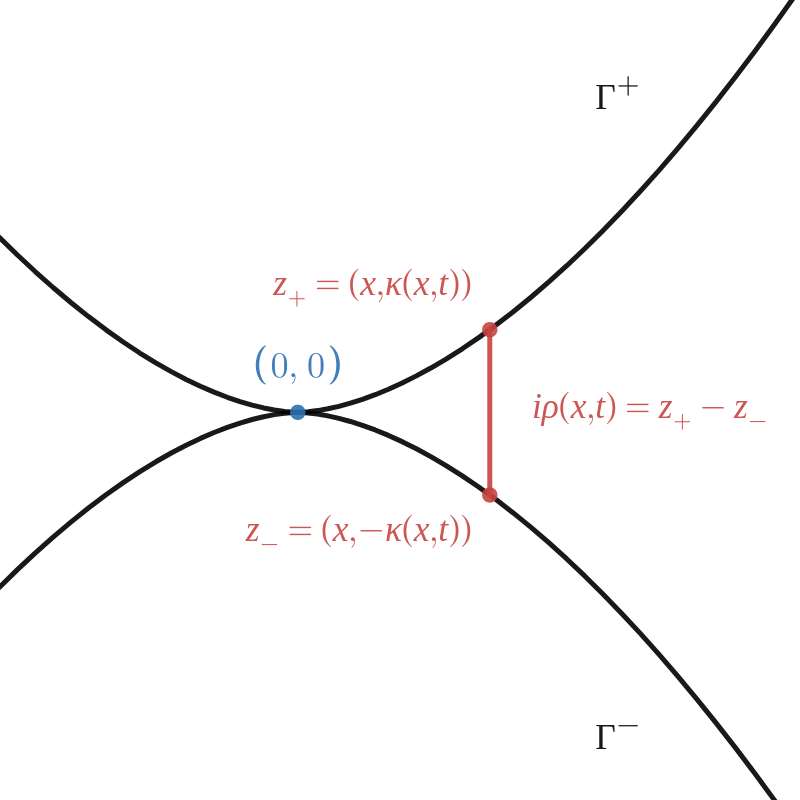}
	\caption{Interface near $z_* = 0$}
	\label{fig2}
\end{figure}

Finally, we introduce some more notation. We append subscript $\pm$ to a point $z\in\Ga$ whenever it is an element of $\Ga^\pm$ in the graph parametrization cf.~\eqref{def:Ga_pm}, i.e.
$$
z_\pm = x \pm i\ka(x, t) \in\Ga^\pm
$$ 
and we use the notation $\rho(x,t) := 2\ka(x,t)$ when we want to emphasize some quantity depends on the difference $z_+ - z_-$ and not on any symmetry assumptions.

For a fixed $x\in I_\de^+ := (0,\de)$, we split $I_{2\de} := (-2\de, 2\de)$ in three intervals 
$$
I_l(x):=[-x,(1-\varepsilon)x], \quad I_c(x):=((1-\varepsilon)x, (1+\varepsilon)x), \quad I_r(x):=(-2\de, -x]\cup [(1+\varepsilon)x, 2\de)
$$ 
for some small $\varepsilon = \varepsilon(t) >0$ and denote the corresponding parts of $\Ga^\pm$ by $\Ga^\pm_l(x)$, $\Ga^\pm_c(x)$ and $\Ga^\pm_r(x)$ respectively, i.e.
$$
\Ga^\pm_{l, c, r}(x) := \{ \, z \in \Ga^\pm \, : \,  x \in I_{l, c, r}(x)  \}.
$$
We will frequently make use of the notation 
$$
F(f)(x, u) := \frac{f(x) - f(u)}{x-u}.
$$

\subsection{Regularity properties of the Birkhoff-Rott integral}
\label{s.BR}

Most of the results of this section hold regardless of any symmetry assumptions. When they do depend on symmetry, it is stated explicitly. However, for simplicity, we work under the assumption that the interface is symmetric with respect to the $x$-axis, although it is mostly used for notational convenience only. In particular, the proofs work in the general case as is. We drop the time dependence throughout this section. Whenever we say a quantity belongs to some weighted Sobolev space depending on parameter $\beta$, we always mean $\beta$ from the regularity assumptions \eqref{regTheta} for the interface parametrization. 

\begin{lemma}\label{lem:BRbasic}
Let $0< \ga + 1/2 <1$ such that $\ga + 1/2 \neq \mu$. Then, 
$$
BR(z, \cdot )^*: \Lcal_{2,\ga}(m) \rightarrow \Lcal_{2,\ga}(m).
$$
The same is true for the operator $BR(z, \cdot )^*z_\al$.
\end{lemma}

\begin{proof} It is enough to prove the claim in a neighborhood of the origin where the interface can be parametrized as a graph (we will do the remaining regions in some detail when we consider derivatives in Lemma \ref{lem:BRcancel}). Let without loss of generality $z=z_+$. We show that   
$$
x \mapsto \frac{1}{2\pi i}\int_{\Ga^+ \cup \Ga^-} f(q)\, \bigg(\frac{z_+'}{z_+ - q}\bigg) \, ds_q \, \in \, \Lcal_{2, \ga}(I_\de)
$$
provided $f\in\Lcal_{2,\ga}(I_{2\de})$. Without loss of generality we may assume $0 < x < \de$. Then, we have 
$$
\bigg|\bigg(\frac{z'_+}{z_+-q_\pm}\bigg)\bigg|\,\lesssim\, \frac{1}{x}, \quad u\in I_l(x), \quad \quad
\bigg|\bigg(\frac{z'_+}{z_+-q_\pm}\bigg)\bigg| \,\lesssim\, \frac{1}{|u|}, \quad  u\in I_r(x) 
$$
and therefore 
$$
\bigg|\int_{\Ga^+_l(x)\cup \Ga^-_l(x)} f(q)\, \bigg(\frac{z'_+}{z_+-q}\bigg)\,ds_q\bigg| \, \lesssim \, \frac{1}{ x}\int_0^{x}|f(u)| + |f(-u)| \, du
$$
on $\Ga^\pm_l(x)$ respectively
$$
\bigg|\int_{\Ga^+_r(x)\cup \Ga^-_r(x)} f(q)\, \bigg(\frac{z'_+}{z_+-q}\bigg)\,ds_q\bigg| \, \lesssim \, \int_{x}^{2\de}\frac{|f(u)| + |f(-u)|}{u}\, du
$$
on $\Ga^\pm_r(x)$. Both are bounded in $\Lcal_{2,\ga}(I_\de)$ by the corresponding Hardy inequalities cf. Appendix, Lemma \ref{lem:hardy}.

It remains to consider $q = q_\pm\in \Ga^\pm_c(x)$. When $q = q_+$, a short calculation yields
\be\label{aux:gac+}
\frac{z_+'}{z_+ - q_+} = \frac{1}{x-u} + \frac{1}{x-u}\bigg(\ka'(x) - \frac{\ka(x) - \ka(u)}{x-u}\bigg)\frac{i}{1 + i \frac{\ka(x) - \ka(u)}{x-u}},
\ee
where, for fixed $u$ and $x$, there exists some $\xi$ such that
$$
\frac{1}{x-u}\bigg(\ka'(x) - \frac{\ka(x) - \ka(u)}{x - u}\bigg) = \pa_x F(\ka)(x,u) = \ka''(\xi) = O(\xi^{\mu - 1}).
$$
Since $x \sim u$, we conclude  
\be\label{aux:gac+o1}
\frac{z_+'}{z_+ - q_+} = \frac{1}{x-u} + O(x^{\mu - 1}).
\ee
In particular, the corresponding integral over $\Ga^+_c(x)$ is bounded in $\Lcal_{2,\ga}$ as required; the error term can be estimated by Hardy's inequality and the Hilbert transform is bounded on weighted Lebesgue spaces whose weight satisfies Muckenhaupt condition.

On the other hand, when $q = q_-$, we have  
$$
\frac{z_+'}{z_+ - q_-} = \frac{1}{(x-u) + i\rho(u)} + O(x^{- 1}),
$$
by Lemma \ref{lem:errorrk} (cf. Appendix), where recall that $\rho(u) = \Im(q_+ - q_-) = 2\ka(u)$. As before, the error term can be estimated by Hardy's inequality and it only remains to estimate  
\be\label{eq:aux-I}
I(x):=\int_{I_c(x)}\frac{1}{(x-u) + i\rho(u)} \, f(u)du.
\ee
At this point, we employ the variable change $h:\Rbb_+\rightarrow (0, 2\de)$ 
$$
h(\xi) = x, \quad h'(\xi) = -\rho(h(\xi))
$$ 
(cf. Appendix for details). The weight function transforms as 
\be\label{eq:VCweight}
x^{2\ga} dx \, \leadsto \, h(\xi)^{2\ga}|h'(\xi)|\, d\xi, \quad \quad h(\xi)^{2\ga}|h'(\xi)| \sim \t m(\xi)^{-2\t \ga}, \quad \quad \t m(\xi) := 1 + \xi,
\ee
where we have set $\t \ga := \mu^{-1}(\ga + 1/2)+1/2$. In particular,
$$
f\in \Lcal_{2,\ga}(m) \quad \Leftrightarrow \quad f\circ h \in \Lcal_{2,-\t\ga}(\t m).
$$ 
As for the kernel, Lemma \ref{lem:kernelDecomposition} implies
$$
\frac{1}{(x-u) + i\rho(u)}\, du \leadsto \frac{h'(\tau)}{(h(\xi) - h(\tau)) - i h'(\tau)}\, d\tau = \left[\frac{1}{(\xi - \tau) - i} + O(\t m(\tau)^{-1})\right] d\tau.
$$
In particular, we have
$$
\aligned
\int_{0}^\de x^{2\ga} |I(x)|^2\,dx \,\lesssim \,  \int_{0}^{\infty} \t m(\xi)^{-2\t \ga} \bigg(&\int_{\wt I_c(\xi)} f(\tau)\,\frac{1}{(\xi - \tau) - i}\, d\tau \bigg)^2 d\xi \\
&+ \int_{0}^{\infty} \t m(\xi)^{-2 \t \ga}\Big(\int_{0}^{(1+\t \varepsilon_+)\xi} \frac{|f(\tau)|}{\t m(\tau)} \, d\tau \Big)^2d\xi,
\endaligned
$$
where $\wt I_c(\xi) := h^{-1}(I_c(x))\subseteq \{\tau \, :\, |\tau - \xi |< \varepsilon_+ \xi \}$ (cf. \eqref{eq:veps12h} in the Appendix). The kernel of the main term is bounded on any weighted Lebesgue space whose weight satisfies the Muckenhaupt condition by the Fourier multiplier theorem. Since $-\t \ga + 1/2 <0$, this is not the case for the weight $\t m^{-\t\ga}$. However, we only integrate over a region where $\tau \sim \xi$, hence we can write  
$$
\int_{\wt I_c(\xi)} f(\tau)\,\frac{1}{(\xi - \tau) - i}\, d\tau = \t m(\xi)^k\int_{\wt I_c(\xi)} \frac{f(\tau)}{\t m(\tau)^k}\,\frac{1}{(\xi - \tau) - i}\, d\tau + O\bigg(\int_{0}^{(1 +  \varepsilon_+)\xi}\frac{|f(\tau)|}{\t m(\tau)}d\tau\bigg),
$$
where $k \in \Nbb$ is such that $0 < (k - \t\ga) + 1/2 < 1$. The error terms are bounded by Hardy inequalities (it is not difficult to see these hold for $\t m$ as well). Recall that we are excluding the case $\ga + 1/2 = \mu$ corresponding to the limiting case $- \t \ga + 1/2 = -1$ for which we don't have this type of Hardy inequalities.
\end{proof}


In general, when $f\in \Lcal_{2, \ga - k}(m)$ with $k\geq 1$ and $0<\ga + 1/2 <1$, we cannot expect $BR(z, f)^*z_\al$ to map $\Lcal_{2, \ga - k}(m)$ to itself. We therefore introduce the following correction:
\be\label{def:BR-k}
\aligned
BR_{-k}(z, f)^* &:= BR(z,f)^* +  \sum_{i = 0}^{k-1}\frac{z(\al)^{i}}{2\pi i}\int_\Ga \frac{1}{z(\al')^{i+1}} \,f(\al') ds_{\al'}\\
&=\frac{z(\al)^k}{2\pi i}\, p.v.\int_{\Ga} \frac{f(\al')}{z(\al')^k} \, \frac{1}{z(\al) - z(\al')} ds_{\al'}.
\endaligned
\ee
For later use we also introduce the notation
\be\label{def:bi}
b_i(f):= \frac{1}{2\pi i}\int_\Ga \frac{1}{z(\al')^{i+1}} \,f(\al') ds_{\al'}.
\ee
When $f\in \Lcal_{2, \ga + k}(m)$, we similarly set
\be\label{def:BR+k}
\aligned
BR_{+k}(z, f)^* &:= BR(z,f)^* -\frac{1}{2\pi i} \sum_{i = 0}^{k-1}\frac{1}{z(\al)^{i+1}}\int_\Ga z(\al')^{i} f(\al') \, ds_{\al'}\\
&=\frac{1}{2\pi i z(\al)^k}\, p.v.\int_{\Ga} z(\al')^k f(\al') \, \frac{1}{z(\al) - z(\al')} ds_{\al'}.
\endaligned
\ee
Lemma \ref{lem:BRbasic} implies 

\begin{corollary}\label{lem:BRcorrections}
Let $k\geq 1$. Under the assumptions of Lemma \ref{lem:BRbasic}, we have
$$
\aligned
f\in \Lcal_{2,\ga - k}(m) &\quad \Rightarrow \quad BR_{-k}(z, f)^*\in \Lcal_{2,\ga- k}(m), \\
f\in \Lcal_{2,\ga + k}(m) &\quad \Rightarrow \quad BR_{+k}(z, f)^*\in \Lcal_{2,\ga+k}(m).
\endaligned
$$
If the interface $\Ga$ is fully symmetric and $f$ is either even or odd with respect to both axes, then
$$
BR(z, f)^* = BR_{-1}(z, f)^*
$$
and, in particular,
$$
BR(z, \cdot )^*: \Lcal_{2, \ga-1}(m) \rightarrow \Lcal_{2, \ga-1}(m).
$$
\end{corollary}


We now give a few results on the derivatives of the Birkhoff-Rott integral $BR(z,f)^*$ when $f$ belongs to certain weighted Sobolev space. We have

\begin{lemma}\label{lem:BRdersI}
Let $\mu \in (1/2, 1]$ and let $1-\mu < \ga + 1/2 < 1$. When $f\in H^1_{\ga}(m)$, we have 
$$
BR(z, f)^*\in H^1_{\ga}(m), \qquad BR(z, f)^*z_s\in H^1_{\ga}(m),
$$ 
where  
\be\label{eq:BR-ipp}
BR(z,f)_s ^* = z_s  BR(z,D_s f)^*, \qquad D_s f := \pa_s\Big(\frac{f}{z_s}\Big).
\ee
Similarly, when $f \in \Lcal^1_{2, \ga}(m)$, we have
$$
BR_{-1}(z, f)^* \in \Lcal^1_{2, \ga}(m), \qquad BR_{-1}(z, f)^*z_s\in \Lcal^1_{2, \ga}(m).
$$
On the other hand, when $f\in \Lcal^1_{2,\ga + 1}(m)$, we have 
$$
BR(z, f)^*\in \Lcal^1_{2, \ga + 1}(m), \qquad BR(z, f)^*z_s\in \Lcal^1_{2, \ga + 1}(m)
$$ 
where now
$$
BR(z,f)_s ^* = z_s  BR_{+1}(z,D_s f)^*.
$$
\end{lemma}

\begin{proof}
Since $1 - \mu < \ga + 1/2$, we have $m^{\mu - 1}\in \Lcal_{2,\ga}(m)$ and therefore
$$
f\in H^1_{\ga}(m) \quad \Rightarrow \quad z_s D_s f= f_s - if\theta_s \, \in \,\Lcal_{2, \ga}(m).
$$
Since we have two 'smoothly' connected cusps and the interface is regular away from the origin, integration by parts and Lemma \ref{lem:BRbasic} yield  
$$
BR(z,f)^*_s = z_s BR(z, D_s f)^* \in \Lcal_{2, \ga}(m)
$$
(see also Lemma \ref{lem:corner-ipp} from Section \ref{s.Corner}). Similarly, we have $BR(z, f)^* \in \Lcal_{2, \ga}(m)$, hence $BR(z, f)^* \in L^\infty$ and we conclude
$$
(BR(z, f)^*z_s)_s =  BR(z, f)^*_s\, z_s + i\theta_s BR(z, f)^*z_s \in \Lcal_{2, \ga}(m).
$$
The remaining statements follow from Corollary \ref{lem:BRcorrections}. Note that 
\be\label{eq:correction-ipp}
\pa_s \left(\frac{1}{z(\al) - z(\al')}\right) = -\frac{z_s(\al)}{z_s(\al')}\, \pa_{s'} \left(\frac{1}{z(\al) - z(\al')} - \frac{1}{z(\al)}\right). 
\ee
\end{proof}

When higher order derivatives of $f$ belong to $\Lcal_{2, \ga}(m)$, we can proceed similarly. However, when $\mu <1$, the resulting space will depend on the values of $f$ at the singular points (assume e.g. that $f\in H^2_\ga(m)$, then, in general, we only have $D_s f = O(m^{\mu - 1})$). As we are only interested in the fully symmetric case, with $f$ typically odd w.r.t to both axes and therefore vanishing at $\pm\al_*$ for all times, we limit ourselves to the following:

\begin{lemma}\label{lem:BRdersII}
Let $\Ga$ be fully symmetric and let $f\in H^2_{\ga,0}(m)$ be odd with respect to both axes. Then, under the assumptions of Lemma \ref{lem:BRdersI}, we have
$$
BR(z,f)^* \in H^2_{\ga,0}(m), \qquad BR(z, f)^*z_s \in H^2_{\max\{\ga, \beta\}}(m). 
$$
\end{lemma}

\begin{proof}
By assumption $1-\mu < \ga + 1/2$ and therefore
$$
f\in H^2_{\ga, 0}(m) \quad \Rightarrow \quad z^2_s D^2_s f = (f_s - i \theta_s f)_s - 2i \theta_s (f_s - i \theta_s f)  \in \Lcal_{2, \ga}(m). 
$$
We can now apply Lemma \ref{lem:BRdersI} to $D_s f$ to conclude 
$$
BR(z,D_s f)^* \in H^1_\ga(m) \quad \Rightarrow \quad BR(z,D_s f)^* \in L^\infty.
$$
In particular, we have $BR(z, f)^*$ and $BR(z, f)^*_s \in L^\infty$ and since the interface is fully symmetric and $f$ is odd w.r.t to both axes we must have $BR(z, f)^* = O(m)$ by integration (the constant of integration vanishes by symmetry).  
\end{proof}



In the next few Lemmas, we consider what happens when we put a derivative on the kernel of the Birkhoff-Rott integral. When $\Ga$ is a sufficiently regular arc-chord curve and $f\in \Lcal_{2,\ga}(m)$, we expect $BR(z,f)^*iz_s - Hf\in\Lcal^1_{2,\ga + 1}(m)$. This is consistent with integrating by parts when the derivative of $f$ is available, i.e. when $f\in \Lcal^1_{2,\ga + 1}(m)$. However, when the interface has a cusp singularity this is no longer true, as putting a derivative on the kernel of $BR(z,f)^*iz_s - Hf$ incurs a loss of $O(m^{\mu +1})$, i.e. the derivative belongs to $\Lcal_{2,\ga + \mu + 1}(m)$ only. In order to cancel the extra $O(m^\mu)$ factor, an additional term is necessary. More precisely, we have

\begin{lemma}\label{lem:BRcancel} Let $1 - \mu < \ga + 1/2 < 1$ with $\ga + 1/2 \neq  \mu$ and let
$$
\Acal f := \left(\frac{f}{2z_s} + BR(z,f)^*\right)iz_\al - H\left(\left(\frac{f}{2z_s} + BR(z,f)^*\right)z_\al\right),
$$
where $H$ denotes the periodic Hilbert transform. Then, we have 
$$
f \in \Lcal_{2, \ga}(m) \quad  \Rightarrow \quad  \Acal f \in H^1_{\ga + 1}(m).
$$
\end{lemma}

\begin{proof} Let us define 
$$
\Rcal f :=2BR(z, f)^*iz_\al - H(|z_\al|f), \qquad \Lcal f :=  \frac{1}{i}H\Rcal f.
$$  
We then have $2\Acal f = (\Rcal -\Lcal)f$. It is enough to consider
$$
\Rcal f(\al) = \frac{1}{2\pi}\, p.v.\int_{-\pi}^\pi f(\al')\bigg(\frac{z_\al}{z(\al) - z(\al')} - \frac{1}{\al - \al'}\bigg)d\al' 
$$
(we absorb $|z_{\al'}|$ into the definition of $f$). By Lemma \ref{lem:BRbasic}, we know that $\Rcal f \in \Lcal_{2, \ga}(m)$. 

We first assume $\al\in [-\pi, \pi]$ is far away from the singular points, i.e. $|\al \pm \al_*|> \epsilon$ and $\epsilon>0$ is so small that $\Ga$ can be parametrized as a graph on $B_{2\epsilon}(\al_*)\cup B_{2\epsilon}(-\al_*)$. Without loss of generality we may assume $\al$ belongs to the ball $A_{\epsilon}:= B_{\al_*- \epsilon}(0)$, i.e. that $|\al| < \al_* - \epsilon$. 
When $\al'\in A_{\epsilon/2}$, we have 
\be\label{eq:k1}
\frac{z_\al}{z(\al) - z(\al')} - \frac{1}{\al - \al'} = \frac{\pa_\al F(z)(\al,\al')}{F(z)(\al,\al')},
\ee
where 
$$
F(z)(\al, \al'):= \frac{z(\al) - z(\al')}{\al -\al'}.
$$
Since $m \sim 1$ on $A_{\epsilon/2}$, we have $z \in H^3(A_{\epsilon/2})$ and therefore we may take one derivatives of \eqref{eq:k1}. In fact, the Taylor development of $z$ gives
$$
\pa_\al^k F(z)(\al,\al') = \frac{1}{(\al - \al')^{k+1}}\int_\al^{\al'}\pa_{\al}^{k+1} z(\tau)\,\frac{(\al' - \tau)^k}{k!}d\tau
$$
and
$$
\inf_{A_{\epsilon}\times A_{\epsilon/2}}|F(z)(\al, \al')| \geq \Fcal_{\epsilon/2}(z)^{-1}.
$$
On the other hand, if $\al'\in ({A_{\epsilon/2}})^c = [-\pi, \pi]\setminus A_{\epsilon/2}$, we have  
$$
|\al - \al'| > \epsilon/2, \qquad \inf_{A_{\epsilon}\times (A_{\epsilon/2})^c}|F(z)(\al, \al')| \geq \Fcal_{\epsilon/2}(z)^{-1}
$$
and we may take the required derivatives directly on the kernel (note that $f\in \Lcal_{2,\ga}(m)\subseteq L^1$). In particular, we have  
$$
\Rcal f \,  \in H^1(A_{\epsilon})
$$
for any small $\epsilon>0$ and therefore also
$$
\Lcal f(\al)  = \frac{1}{\pi i} \, p.v.\int_{A_{\epsilon/2}}\frac{1}{\al - \al'}\, \Rcal f(\al')\,d\al' + \frac{1}{\pi i} \int_{(A_{\epsilon/2})^c} \frac{1}{\al - \al'}\, \Rcal f(\al')\,d\al' \in H^1(A_\epsilon)
$$
where we replace $\epsilon$ by $\epsilon/4$ and use $\Rcal f \in H^1(A_{\epsilon/4})$ when considering derivatives of the principal value part of the integral and we use  $\Rcal f \in \Lcal_{2, \ga}(m)$ when estimating derivatives of the second.

Let now $\al\in \overline{B_{\epsilon}(\al_*)}$. We first claim  
$$
\Rcal f(\al) = \frac{1}{2\pi}\int_{B_{2\epsilon}(-\al_*)} f(\al)\bigg(\frac{z_{\al}}{z(\al) - z(\al')}\bigg) d\al' + \t \Rcal f(\al), \quad  \t \Rcal f\in H^1_{\ga + 1}(B_{\epsilon}(\al_*)).
$$
Indeed, far away from the singular points the arc-chord condition holds and $|\al - \al'|$ is bounded away from zero, hence the derivative of the corresponding part of the integral has the same regularity as $z_{\al\al}$. When $\al'\in B_{2\epsilon}(\al_*)$, Lemma \ref{lem:derFphiI} from the Appendix implies
$$
\al\mapsto \int_{B_{2\epsilon}(\al_*)} \frac{\pa_\al F(z)(\al,\al')}{F(z)(\al,\al')}\, f(\al')d\al'\in H^1_{\ga + 1}(B_{\epsilon}(\al_*)).
$$
The claim follows, since $|\al - \al'|$ is bounded away from zero when $\al'\in B_{2\epsilon}(-\al_*)$. In particular, it is not difficult to see that we have
$$
\Lcal f(\al) = \frac{1}{\pi i}\,p.v.\int_{B_{2\epsilon}(\al_*)}\frac{1}{\al - \al'} \bigg(\frac{1}{2\pi}\, p.v.\int_{B_{4\epsilon}(-\al_*)} f(\al'') \, \frac{z_{\al'}}{z(\al') - z(\al'')} d\al''\bigg)d\al' + H^1_{\ga + 1}(B_{\epsilon}(\al_*))
$$
where we have used that 
$$
p.v. \int_{B_{2\epsilon}(\al_*)} \frac{\t \Rcal f(\al')}{\al - \al'}\, d\al'\in H^1_{\ga + 1}(B_{\epsilon}(\al_*)),
$$
cf. Appendix, Lemma \ref{lem:trHilbert} (where we replace $\epsilon$ by $2\epsilon$ and use $\t \Rcal f \in H^1(B_{4\epsilon})$). 


At this point, we go over to the graph parametrization, where we set $x=z_1(\al)$ resp. $u= z_1(\al')$. Since
$$
\frac{1}{\al - \al'}- \frac{z_{1{\al}}}{z_1(\al) - z_1(\al')}= \frac{\pa_\al F(z_1)(\al, \al')}{F(z_1)(\al, \al')}
$$
and $z_{1\al}\in H^2_{\beta + 1}(B_{\epsilon}(\al_*))$, the corresponding integral belongs to $H^1_{\ga + 1}(B_{\epsilon}(\al_*))$ by Lemma \ref{lem:derFphiI} cf. Appendix. In particular, it remains to consider the derivative of 
$$
\frac{1}{2\pi }\, \int_{-2\de}^{2\de} f(u)\, \frac{z_+'}{z_+- q_-}du - \frac{1}{\pi i}\,p.v.\int_{-\de}^{\de}\frac{1}{x - u} \bigg(\frac{1}{2\pi }\, \int_{-2\de}^{2\de} f(\t u)\,\frac{q_+'}{q_+-\t q_-} d \t u\bigg)d u, \quad x\in I_{\de/2},
$$
where $f\in \Lcal_{2, \ga}(I_{2\de})$ for some $\de>0$ small enough. Without loss of generality we only consider $x\in I_{\de/2}$ with $x>0$. The kernel can be written as
$$
\frac{z_+'}{z_+- q_-} = s(x,  u) + r(x,  u), \qquad s(x, u):= \frac{1}{(x - u) + i \rho(u)}
$$
where 
\be\label{eq:rxureg}
u \, \rightarrow  \, \int_{-2\de}^{2\de} f(\t u) r(u, \t u)\,d\t u \in H^1_{\ga + 1}(I_{\de}).
\ee
In fact, by Lemma \ref{lem:errorrk}, we know 
$$
|\pa_x r(x,u)| \, \lesssim \, \frac{1}{x}\bigg(\frac{1}{x}  +  \frac{\rho(u)}{(x-u)^2 + \rho(u)^2}\bigg), \quad  u\in I_c(x).
$$
On the other hand, we have
$$
\frac{d}{dx}\Big(\frac{z'_+}{z_+ - q_-}\Big) = -\Big(\frac{z'_+}{z_+ - q_-}\Big)^{2} + \frac{1}{z_+ - q_-} \, O(x^{\mu -1})
$$
and therefore, 
\be\label{eq:ders-IlIr}
\bigg|\frac{d}{dx}\Big(\frac{z'_+}{z_+ - q_-}\Big)\bigg|\,\lesssim \, \frac{1}{x^2}, \quad u\in I_l(x), \qquad \bigg|\frac{d}{dx}\Big(\frac{z'_+}{z_+ - q_-}\Big)\bigg|\,\lesssim \, \frac{1}{u^2}, \quad u\in I_r(x).
\ee
Clearly, when $u\in I_{2\de}\setminus I_c(x)$, the kernel $\pa_x s(x,u)$ also satisfies these estimates, hence  \eqref{eq:rxureg} follows as in the proof of Lemma \ref{lem:BRbasic}, we omit the details. Moreover, by Lemma \ref{lem:trHilbert} its Hilbert transform belongs to $H^1_{\ga + 1}(I_{\de/2})$  and therefore
$$
\Rcal f(x) = \Scal f(x) + H^1_{\ga + 1}(I_{\de/2}), \qquad \Lcal f(x) = \frac{1}{\pi i}\,p.v.\int_{-\de}^{\de}\frac{1}{x - u} \,\Scal f(u) du + H^1_{\ga + 1}(I_{\de/2}),
$$
where we have set
$$
\Scal f(x) := \frac{1}{2\pi}\int_{-2\de}^{2\de} f( u) s(x, u)\,d u. 
$$
We further claim 
\be\label{eq:cl3}
\Lcal f(x) = \frac{1}{\pi i}\,p.v.\int_{I_c(x)} \frac{1}{x - u} \,\Scal_c f(u) du + H^1_{\ga + 1}( I_{\de/2})
\ee
where we have set
$$
\Scal_c f(x) := \frac{1}{2\pi}\int_{I_c(x)} f(u) s(x,u)du.
$$ 
Indeed, it is not difficult to see that 
$$
x \mapsto \int_{I_{2\de}\setminus I_c(x)} f(u) s(x, u)du \in H^1_{\ga + 1}(I_{\de})
$$
(the interval of integration depends on $x$, but $|s(x, (1\pm \varepsilon)x)| = O(x^{-1})$), hence its Hilbert transform belongs to $H^1_{\ga + 1}(I_{\de/2})$ by Lemma \ref{lem:trHilbert}. We then restrict the Hilbert transform integral to $I_c(x)$. In fact, for $u\in I_l(x)$, we have
$$
\Big|\int_{I_l(x)}\frac{1}{(x - u)^{2}}\,\Scal_c f(u)du\Big| \,\lesssim\, \frac{1}{x^2}\int_{- x}^{x}\, |\Scal_c f(u)|du, 
$$
while, when $u \in I_r(x)$, we have
$$
\Big|\int_{I_r(x)}\frac{1}{(x - u)^2}\,\Scal_c f(u)du\Big| \,\lesssim\, \frac{1}{x}\int^{\de}_{x}\frac{|\Scal_c f(u)| + |\Scal_c f(-u)|}{u}\,du.
$$
Both are bounded in $\Lcal_{2,\ga+1}$. In particular, \eqref{eq:cl3} follows. At this point, we apply the variable change 
$$
x = h(\xi), \quad h'(\xi) = -\rho(h(\xi))
$$ 
cf. Appendix for details. Recall that $h^{-1}$ maps $(0, 2\de)$ to $\Rbb_+$ and that 
$$
f\in \Lcal^1_{2, \ga + 1}((0, 2\de)) \quad \Leftrightarrow \quad f\circ h \in \Lcal^1_{2, 1-\t \ga}(\Rbb_+, \t m), 
$$
where $\t \ga = \mu^{-1}(\ga + 1/2) +1/2$ and $\t m(\xi) = 1 + \xi$ (cf. proof of Lemma \ref{lem:BRbasic}). We claim
$$
(\Lcal f) \circ h(\xi) =  \frac{1}{\pi i}\int_{\wt I_c(\xi)} \frac{1}{\xi -\tau}\,\Scal_c f(\tau)d\tau  \, + \, \Lcal^1_{2, 1 -\t \ga}(\wt I_{\de/2}),
$$ 
where we use the notation $\wt I_c(\xi):=h^{-1}(I_c(x))$ and $\wt I_{\de/2} := h^{-1}((0, \de/2))$. Indeed, we have 
$$
\bigg|\int_{\wt I_c(\xi)}\pa_\xi\bigg(\frac{h'(\tau)}{h(\xi) - h(\tau)}-\frac{1}{\xi -\tau}\bigg)\Scal_c f(\tau)\, d\tau \bigg| \, \lesssim \, \frac{1}{\t m(\xi)^2}\int_{0}^{(1 + \varepsilon_+) \xi} |\Scal_c f(\tau)|\, d\tau
$$
(cf. Lemma \ref{lem:errVCHilb} and estimate \eqref{eq:veps12h} from the Appendix) with the right-hand side bounded in $\Lcal_{2, 1-\t \ga}(\wt I_{\de/2})$ by Hardy's inequality. The interval of integration also depends on $\xi$, but it is not difficult to see that terms coming from this part of the derivative also belong to $\Lcal_{2, 1 - \t \ga}(\wt I_{\de/2})$ cf. estimate following \eqref{eq:veps12h} in the Appendix.  

On the other hand, the kernel of $\Scal_c f$ transforms as 
$$
\frac{1}{(x- u) + i \rho(u)}\, du \quad \leadsto \quad  \frac{h'(\tau)}{(h(\xi) - h(\tau))- ih'(\tau)}\, d\tau = \left[\frac{1}{(\xi - \tau) - i}+ r(\xi, \tau)\right]d\tau
$$
where the remainder satisfies
\be\label{aux:rVC}
|r(\xi, \tau)| \, \lesssim \, \frac{1}{\t m(\tau)},\qquad
|\pa_\xi r(\xi, \tau)|\, \lesssim \, \bigg|\frac{1}{\t m(\tau)}\bigg(\frac{1}{\t m(\tau)}  + \frac{1}{(\xi - \tau)^2 + 1}\bigg)\bigg|
\ee 
for $\tau \in \wt I_c(\xi)$ uniformly for $\xi \in \wt I_{\de}$ (cf. Lemma \ref{lem:kernelDecomposition} in the Appendix). In particular,  
$$
\tau \mapsto \int_{\wt I_c(\tau)}f(\tau')r(\tau, \tau') \, d\tau' \in \Lcal^1_{2, 1 - \t\ga}(\t m).
$$
By the first estimate of \eqref{aux:rVC}, boundary terms also give correct estimates. Moreover, it is not difficult to see its Hilbert transform over $\wt I_c(\xi)$ also belongs to $\Lcal^1_{2, 1 - \t\ga}(\wt I_{\de/2})$, since we are integrating over a set where $\xi \sim \tau$. 

In particular, it remains to show
$$
\int_{\wt I_c(\xi)}f(\tau)\frac{1}{(\xi - \tau) - i} \, d\tau - \frac{1}{\pi i}\int_{\wt I_c(\xi)}\frac{1}{\xi -\tau}\, \int_{\wt I_c(\tau)}f(\t \tau)\frac{1}{(\tau - \t \tau) - i} \, d\t \tau \in \Lcal^1_{2, 1 - \t \ga}(\wt I_{\de/2}).
$$
Since $- \t \ga$ is not Muckenhaupt (we have $-\t \ga + 1/2 < 0$), we write
$$
\aligned
\int_{\wt I_c(\xi)} f(\tau) \,  \frac{1}{(\xi - \tau) - i} \, d\tau &= \t m(\xi)^k \int_{ \wt I_c(\xi)} \frac{f(\tau)}{\t m(\tau)^k} \, \frac{1}{(\xi - \tau) - i} \, d\tau + \Lcal^1_{2, 1-\t \ga}(\wt I_{\de})\\
&=\t m(\xi)^k \underbrace{\int_{\Rbb} \frac{f(\tau)}{\t m(\tau)^k} \, \frac{1}{(\xi - \tau) - i} \, d\tau}_{ = : G(\xi)} + \Lcal^1_{2, 1-\t \ga}(\wt I_{\de}),
\endaligned
$$
where $k\in \Nbb$ is such that $0 < (k - \t \ga) + 1/2 < 1$ and we have extended $f$ by zero to the entire real line. We have $G\in \Lcal_{2, k - \t \ga}(\Rbb_+, \t m)$, hence 
$$
G(\xi) = \frac{1}{\pi i}\, p.v. \int_\Rbb \frac{1}{\xi - \tau}\, G(\tau)d\tau,
$$
where it is not difficult to see that 
$$
G(\tau) =  \int_{\wt I_c(\tau)} \frac{f(\t \tau)}{\t m(\t \tau)^k}\frac{1}{(\tau - \t \tau) - i} \, d\t \tau + \Lcal^1_{2, k + 1-\t \ga}(\Rbb_+, \t m),
$$
where the Hilbert transform of the error term belongs to $\Lcal^1_{2, k + 1 -\t \ga}(\t m)$ by  Lemma \ref{lem:trHilbert}. Since, we also have 
$$
\int_{\Rbb \setminus \wt I_c(\xi)} \frac{1}{\xi - \tau} \int_{\wt I_c(\tau)} \frac{f(\t \tau)}{\t m(\t \tau)^k}\frac{1}{(\tau - \t \tau) - i} \, d\t \tau \in \Lcal^1_{2, k + 1 - \t \ga}(\wt I_{\de/2})
$$
the claim follows (note that we only need to consider $\Rbb_+ \setminus \wt I_c(\xi)$, since by definition the inner integral vanishes when $\tau<0$).
\end{proof}

We now give an extension of Lemma \ref{lem:BRcancel}, needed for the next section: 

\begin{lemma}\label{lem:struct1}
Let $1 - \mu < \ga + 1/2 < 1$ with $\mu \neq \ga + 1/2$ and let $f$ and $g$ be complex-valued. Then 
$$
\Acal(f, g) := \Re\left( \frac{ifg}{2} +  ifBR(z, g)^*z_s \right) - H \Re\left( \frac{fg}{2} +  f BR(z, g)^* z_s\right)
$$
satisfies 
$$
g\in H^k_{\ga + k}(m), \, \,  f\in H^{k + 2}_{\beta + 1 + k}(m) \quad \Rightarrow \quad \Acal(f, g) \in H^{k + 1}_{\ga + 1 + k}(m), \quad k=0,1.
$$

\end{lemma}

\begin{proof}
We have
$$
\aligned
 \Re\bigg(  \frac{ifg}{2}  +  i f BR(z, g)^*z_s\bigg) &= \Re f \Big(-\frac{\Im g}{2} + BR(z, g)\cdot z_s^\perp\Big) - \Im f \Big(\frac{\Re g}{2} + BR(z, g)\cdot z_s\Big)\\
 \Re\bigg( \frac{fg}{2} + f BR(z, g)^*z_s\bigg) &= \Re f \Big(\frac{\Re g}{2} + BR(z, g)\cdot z_s\Big) + \Im f \Big(- \frac{\Im g}{2} + BR(z, g)\cdot z_s^\perp\Big).
\endaligned
$$
Since 
$$
\aligned
-\frac{\Im g}{2} + BR(z, g)\cdot z_s^\perp &= BR(z, \Re g) \cdot z_s^\perp - \Big(\frac{\Im g}{2} + BR(z,\Im g)\cdot z_s\Big), \\
\frac{\Re g}{2} + BR(z, g)\cdot z_s &= \Big(\frac{\Re g}{2} + BR(z,\Re g)\cdot z_s\Big) + BR(z, \Im g) \cdot z_s^\perp,
\endaligned
$$
after some rearrangement we obtain the formula
$$
\aligned
\Acal(f&, g)= -\Re(f)\bigg[\Big(BR(z,\Im g)\cdot z_s + \frac{\Im g}{2}\Big) + H \big(BR(z,\Im g)\cdot z_s^\perp\big) \bigg] + [\Re(f), H ]BR(z,\Im g)\cdot z_s^\perp\\
&+ \Re(f)\bigg[BR(z,\Re g)\cdot z_s^\perp - H\Big(BR(z,\Re g)\cdot z_s + \frac{\Re g}{2}\Big)\bigg] + [\Re(f), H ]\bigg(BR(z,\Re g)\cdot z_s + \frac{\Re g}{2}\bigg) \\
& -\Im(f)\bigg[\Big(BR(z,\Re g)\cdot z_s + \frac{\Re g}{2}\Big) + H \big(BR(z,\Re g)\cdot z_s^\perp\big) \bigg] + [\Im(f), H ]BR(z,\Re g)\cdot z_s^\perp\\
&-\Im(f)\bigg[BR(z,\Im g)\cdot z_s^\perp - H\Big(BR(z,\Im g)\cdot z_s + \frac{\Im g}{2}\Big)\bigg] - [\Im(f), H ]\bigg(BR(z,\Im g)\cdot z_s + \frac{\Im g}{2}\bigg).
\endaligned
$$
When $g\in \Lcal_{2, \ga}(m)$, the claim follows using Lemma \ref{lem:BRcancel} on the terms in brackets and Lemma \ref{lem:derFphi} from the Appendix which deals with commutators. However, the definition of $\Acal f$ from Lemma \ref{lem:BRcancel} differs from the above by a factor of $1/|z_\al|$, hence in order to use Lemma \ref{lem:BRcancel}, we write 
$$
\left(BR(z,\Im g)\cdot z_s + \frac{\Im g}{2}\right) + H \big(BR(z,\Im g)\cdot z_s^\perp\big) = \frac{1}{|z_\al|} \Im \Acal (\Im g) + \left[\frac{1}{|z_\al|} , H\right] BR(z, \Im g)\cdot z_\al^\perp
$$
and similarly for other terms in brackets. Since $1/|z_\al|\in H^2_{\beta + 1}(m)$, the commutator clearly belongs to $H^1_{\ga + 1}(m)$ by Lemma \ref{lem:derFphi}. 

Assume now $g\in H^1_{\ga + 1}(m)$ with $f\in H^3_{\beta + 2}(m)$. First note that commutators belong to $H^2_{\ga + 2}(m)$ by Lemma \ref{lem:derFphi}. As for the terms in brackets, we need to take one derivative directly (we control only two derivatives of $1/|z_\al|$), then use Lemma \ref{lem:BRcancel}. Since $D_s g\in \Lcal_{2, \ga + 1}(m)$ is not integrable, we have
$$
\pa_\al BR(z, g)^* =  z_\al BR_{+1}(z, D_sg)^* = \frac{z_\al}{z}\, BR(z, zD_s g)^*, \quad \al\in [-\pi, \pi]
$$
(cf. definition \eqref{def:BR+k}), where $1/z(\al)$ has singularities at both $\pm \al_*$. When considering corresponding derivatives of the Hilbert transform we therefore have to proceed as in the proof of Lemma \ref{lem:BRcancel} and consider $\al \in [-\pi, \pi]\setminus B_\epsilon(\al_*)\cup B_\epsilon(-\al_*)$ and $\al \in B_\epsilon(\pm\al_*)$ separately. 

W.lo.g. we only consider $\al \in B_\epsilon(\al_*)$. Then, as in the proof of Lemma \ref{lem:trHilbert}, we have
$$
\aligned
\pa_\al H(BR(z, g)^*z_s) &= \frac{1}{\al - \al_*} \int^\pi_0 \frac{1}{\al - \al'}(\al' - \al_*)\pa_{\al'} (BR(z, g)^*z_s(\al'))\,d\al' + H^1_{\ga + 2}(B_\epsilon(\al_*)) \\
&= \frac{1}{\al - \al_*} \int^\pi_0 \frac{1}{\al - \al'} (\al' - \al_*) \frac{z_s(\al')}{z(\al')} BR(z, zD_s g)^*z_{\al'} d\al' + H^1_{\ga + 2}(B_\epsilon(\al_*)) \\
&= \frac{z_s(\al)}{z(\al)} \int^\pi_0  \frac{1}{\al- \al'} \, BR(z, z D_s g)^*z_{\al'}d\al' + H^1_{\ga + 2}(B_\epsilon(\al_*)).
\endaligned
$$
In the last step, we have used Lemma \ref{lem:derFphi} with $\phi := (\al - \al_*)\frac{z_s}{z}\in H^2_{\beta + 1}([0, \pi])$, since 
$$
\pa_\al \Big(\frac{\al - \al_*}{z(\al)}\Big) = \frac{1}{z(\al)^2}(z(\al) - z(\al_*) - z_\al(\al)(\al - \al_*)) \in H^1_{\beta + 1}([0, \pi]).
$$ 
In particular, the claim follows applying Lemma \ref{lem:BRcancel} to $\Acal(zD_s g)$. 
\end{proof}

\begin{lemma}\label{lem:BR-additional-cancelation}
	Let $0 < \ga + 1/2 < 1$ such that $\mu \neq \ga + 1/2$ and let 
	$$
	I(f) :=  \Re \bigg(\frac{1}{2\pi i}\int_\Ga f(\al')\bigg(\frac{z_s(\al)}{z_s(\al')} + 1\bigg)\frac{z_s(\al)}{z(\al) - z(\al')}\, ds_{\al'}\bigg) 
	$$
	Then, we have 
	$$
	\aligned
	f\in\Lcal_{2,\ga}(m) \quad &\Rightarrow \quad I(f) \in H^1_{\ga + 1}(m),\\
	f\in\Lcal^1_{2,\ga+1}(m) \quad &\Rightarrow \quad I(f) \in H^1_{\ga + 1 - \lambda}(m).
	\endaligned
	$$ 
	where $\lambda := 1-(\beta + 1/2)$. 
\end{lemma}

\begin{proof}
	The derivative of $I(f)$ reads
	$$
	\aligned
	\pa_s I(f) = \Re \bigg(\frac{1}{2\pi i}\int_\Ga f(\al')\bigg(\frac{z_s(\al)}{z_s(\al')} &+ 1\bigg) \, \frac{d}{ds_\al}\Big(\frac{z_s(\al)}{z(\al) - z(\al')}\Big)\, ds_{\al'}\bigg) +\\
	&+ \frac{\theta_s}{2\pi}\Re \bigg(\int_\Ga f(\al')\frac{z_s(\al)}{z_s(\al')} \, \frac{z_s(\al)}{z(\al) - z(\al')}\, ds_{\al'}\bigg) =: I_1 + I_2.
	\endaligned
	$$
	First note that $1 - \mu < \beta + 1/2$ implies $\Lcal_{2,\ga + 1 -\mu}(m)\subseteq\Lcal_{2,\ga + 1 -\lambda}(m)$. When $f\in\Lcal_{2,\ga}(m)$, we clearly have $I_2\in \Lcal_{2,\ga + 1 -\mu}(m)$ (using $\theta_s = O(m^{\mu-1})$). 
	
	To estimate $I_1$, assume without loss of generality $z = z_+\in \Ga^+$ with $q = q_\pm\in \Ga^\pm$ and let 
	$$
	I_1 = J^+ + J^- + J
	$$
	corresponding to integrals over $\Ga^\pm$ respectively the integral over $\Ga\setminus (\Ga^+\cup \Ga^-)$.
	 
	Passing over to the graph parametrization, we have $z_s(\al) \rightarrow z_s^+(x):=\frac{z_+'}{|z_+'|}$ and $z_s(\al')\rightarrow \pm q_s^\pm(u):=\pm \frac{q_\pm'}{|q_\pm'|}$. Therefore, it is not difficult to see that 
	\be\label{aux:additional-real}
	\Re \left(\frac{z^+_s}{q_s^+} + 1\right) = O(1), \quad q\in \Ga^+ , \qquad \Re \left(-\frac{z^+_s}{q_s^-} + 1\right) =  O(|\ka'(x)|^2 + |\ka'(u)|^2), \quad q\in\Ga^-
	\ee
	(note that $-\Re(z_+' (q_-')^*) = -1 + \ka'(x)\ka'(u)$) and
	\be\label{aux:additional-img}
	\Im \left(\pm\frac{z^+_s}{q_s^\pm} + 1\right) = O(|\ka'(x) \mp \ka'(u)|), \quad q\in \Ga^\pm.
	\ee
	We first consider $J^+$. We can write 
	$$
	\aligned
	J^+ = \Re \bigg(\frac{1}{2\pi i}&\int_{-2\de}^{2\de} f(u) g_+(x,u) \, \frac{d}{dx}\Big(\frac{z_+'}{z_+ - q_+} - \frac{1}{x-u}\Big)\, du\bigg) \\
	&- \frac{1}{2\pi}\int_{-2\de}^{2\de} f(u) \Im g_+(x,u)\frac{1}{(x-u)^2} \,du  +  \Lcal_{2, \ga + 1 - \mu}(m) =: J_1^+ + J_2^+ +  \Lcal_{2, \ga + 1 - \mu}(m),
	\endaligned
	$$
	where the real and the imaginary parts of $g_+(x,u)$ satisfy \eqref{aux:additional-real} and \eqref{aux:additional-img}-type estimates with $q\in \Ga^+$ respectively. Since $\Im g_+(x,u) = \frac{1}{|z_+'|^2}(\ka'(u) - \ka'(x))$, it is not difficult to see that we can write
	$$
	\Im g_+(x,u)\frac{1}{(x-u)^2} = \frac{1}{|z_+'|^2}\Big(F(\ka')_x + \frac{\ka''(x)}{x-u}\,\Big),
	$$
	hence $J^+_2\in \Lcal_{2, \ga + 1 - \lambda}(m)$ by Lemma \ref{lem:derFphi}. On the other hand, the kernel of $J^+_1$ can be written as
	$$
	\frac{d}{dx}\bigg(\frac{z_+'}{z_+ - q_+} - \frac{1}{x-u}\bigg)  = - \frac{iF(\ka)_{xx}}{1 + iF(\ka)} - \frac{F(\ka)_{x}^2}{(1 + iF(\ka))^2}
	$$
	cf. \eqref{aux:gac+}, hence using $g_+ = O(1)$, we conclude
	$$
	|J^+_1| \, \lesssim \, \int_{-2\de}^{2\de} |f(u)| \big(|F(\ka)_{xx}| + |F(\ka)_{x}|^2\big)du
	$$ 
	which can be estimated as in Lemma \ref{lem:derFphi} to show it belongs to $\Lcal_{2, \ga + 1 -\lambda}(m)$ (note that $F(\ka)_{xx}$ satisfies similar estimates as $F(\ka')_x$ does; we omit further details).
	
	We consider $J^-$ next. We can write
	$$
	J^- = \Re \bigg(\frac{1}{2\pi i}\int_{-2\de}^{2\de} f(u) g_-(x,u) \, \frac{d}{dx}\Big(\frac{z_+'}{z_+ - q_-}\Big)\, du\bigg) +  \Lcal_{2, \ga + 1 - \mu}(m),
	$$
	where the real and the imaginary parts of $g_-(x,u)$ satisfy \eqref{aux:additional-real} and \eqref{aux:additional-img} with $q\in \Ga^-$. We have 
	$$
	\bigg|\frac{d}{dx}\Big(\frac{z_+'}{z_+ - q_-}\Big) \bigg| \,\lesssim \, \frac{1}{x^2}, \quad u\in I_l(x), \qquad \bigg|\frac{d}{dx}\Big(\frac{z_+'}{z_+ - q_-}\Big) \bigg| \,\lesssim \, \frac{1}{u^2} + \frac{x^{\mu-1}}{u}, \quad u\in I_r(x).
	$$
	Since $g_-(x, u) = O(x^\mu + |u|^\mu)$, we can use Hardy inequalities to show the corresponding integrals are bounded in $\Lcal_{2,\ga + 1-\mu}(m)$. However, when $u\in I_r(x)$, the corresponding Hardy inequality is true only under the assumption that $(\ga + 1 - \mu) + 1/2 > 0$. We therefore use an improved estimate 
	$$
	|g_-(x, u)| \, \lesssim \, 2|\ka'(x)| + |\ka'(x)-\ka'(u)| \, \lesssim \, x^\mu + x^{\mu-1}|x-u|, \qquad u\in I_r(x),
	$$ 
	which implies 
	$$
	\bigg|g_-(x, u)\frac{d}{dx}\Big(\frac{z_+'}{z_+ - q_-}\Big) \bigg| \,\lesssim \, \frac{x^\mu}{u^2} + \frac{x^{\mu-1}}{u}, \qquad u\in I_r(x).
	$$
	The corresponding integral is now bounded in $\Lcal_{2,\ga + 1-\mu}(m)$ as long as $0 < \ga + 1/2$. 
	
	It remains to estimate $u\in I_c(x)$. We have 
	$$
	\frac{d}{dx}\Big(\frac{z_+'}{z_+ - q_-}\Big) = -\frac{1}{ \big((x-u) + i \rho(u)\big)^2} + O(x^{-2}), \qquad u\in I_c(x),
	$$
	(cf. Appendix, Lemma \ref{lem:errorrk}) and therefore
	$$
	J^- = -\Re \bigg(\frac{1}{2\pi i}\int_{I_c(x)} f(u) g_-(x,u) \, \bigg(\frac{1}{ (x-u) + i\rho(u)}\bigg)^2\, du\bigg) +  \Lcal_{2, \ga + 1 - \mu}(m).
	$$
	However, we have $g_-(x,u) = O(x^\mu)$ and 
	$$
	\Big|\frac{1}{ (x-u) + i\rho(u)}\Big|^2 \, \lesssim \, \frac{1}{\rho(x)}\frac{\rho(u)}{(x-u)^2 + \rho(u)^2}, \qquad u\in I_c(x),
	$$
	hence when $f\in \Lcal_{2,\ga}(m)$ this part of the integral belongs to $\Lcal_{2,\ga+1}(m)$ only (cf. proof of Lemma \ref{lem:BRbasic} for details). On the other hand, if $f\in \Lcal^1_{2,\ga + 1}(m)$, we can write
	$$
	\Big(\frac{1}{ (x-u) + i\rho(u)}\Big)^2 = \pa_u \Big(\frac{1}{ (x-u) + i\rho(u)}\Big) - \frac{\rho'(u)}{ \big((x-u) + i\rho(u)\big)^2} 
	$$
	where the integral corresponding to the second kernel on the r.h.s. clearly belongs to $\Lcal_{2,\ga + 1 - \mu}(m)$. Finally, we can integrate by parts, to obtain
	$$
	J^- = \Re \bigg(\frac{1}{2\pi i}\int_{I_c(x)} f'(u) g_-(x,u) \, \frac{1}{ (x-u) + i\rho(u)}\, du\bigg) +  \Lcal_{2, \ga + 1 - \mu}(m).
	$$
	(terms coming from the integration limits belong to $\Lcal_{2,\ga + 1 - \mu}(m)$, since the kernel is $O(x^{\mu-1})$ when $u = (1\pm\epsilon)x$; note also that $\pa_u g_- = O(u^{\mu - 1})$). We can further write  
	$$
	J^- = \frac{1}{2\pi}\int_{I_c(x)} f'(u) \Im g_-(x,u) \, \frac{x-u}{ (x-u)^2 + \rho(u)^2}\, du +  \Lcal_{2, \ga + 1 - \mu}(m),
	$$
    where $\Im g_-(x,u) = -\frac{1}{|z_+'|}(\ka'(x) + \ka'(u))$, but this is easily seen to be bounded in $\Lcal_{2,\ga + 1 -\mu}(m)$ (see the proof of Lemma \ref{lem:BRbasic}).
\end{proof}



We finish this section with a series of Lemmas which identify certain cancellations in the Birkhoff-Rott integral. We also give growth estimates valid near the singular point, provided we control a sufficient number of derivatives of $f$. These results remain true without any assumptions on the symmetry of the interface.

\begin{lemma}\label{lem:BRdiff} Let $0 < \ga + 1/2 < 1$, where $\ga + 1/2 \neq \mu$. Then,   
	$$
	\Dcal(f)(x) := \Big(\frac{if_+(x)}{2} + BR(z,f)^* iz_+' \Big) + \Big(\frac{if_-(x)}{2} - BR(z,f)^* iz_-' \Big), \quad x\in I_\de,
	$$
	where $f_\pm = f|_{\Ga^\pm}$, satisfies
	$$
	\aligned
	f\in H^1_{\ga+1}(m) \quad &\Rightarrow \quad \Dcal(f)\in \Lcal_{2, \ga - \mu}(I_\de),\\
	f \in H^1_{\ga}(m) \quad &\Rightarrow \quad \Dcal(f) = O(m^\mu) +  \Lcal_{2, \ga - \mu - 1}(I_\de).
	\endaligned
	$$
\end{lemma}

\begin{proof}
	We estimate the integral
	\be\label{eq:Dker}
	x \mapsto \frac{1}{2\pi }\int_{\Ga} f(q)\bigg(\frac{z'_+}{z_+ - q} - \frac{z'_-}{z_- - q} \bigg) du, \quad x\in I_\de,
	\ee
	which we split as follows:
	\be\label{eq:Ix>b}
	\aligned
	\int_{\Ga^+_l(x)}f(u) \bigg(\frac{z'_+}{z_+ - q_+} - &\frac{z'_-}{z_- - q_+} \bigg)  du \pm \int_{\Ga^+_c(x)} f(u) \bigg(\frac{z_\pm'}{z_\pm - q_+} - \frac{1}{x-u}\bigg) du \\
	&+ \int_{\Ga^+_r(x)} + \int_{\Ga \setminus (\Ga^+ \cup \Ga^-)} f(u)\bigg(\frac{z'_+}{z_+ - q} - \frac{z'_-}{z_- - q} \bigg) du= 2\pi \sum_{i = 1}^5 I_i(x).
	\endaligned
	\ee
	(The  integral over $q\in \Ga^-$ can be estimated in the same way as the one over $\Ga^+$ by interchanging $z_+$ and $z_-$). Moreover, for $q\in \Ga^\pm$, we define  
	\be\label{aux1:k}
	\aligned
	k(x,u) :&= \frac{z'_+}{z_+ - q} - \frac{z'_-}{z_- - q} \\
	&= \frac{z'_+ - z'_-}{z_+ - q} -  \frac{z'_-}{z_- - q}\frac{z_+ - z_-}{z_+ - q}. 
	\endaligned
	\ee
	
	Assume without loss of generality $x>0$. The kernel of $I_1(x)$ satisfies the estimate
	$$
	|k(x,u)| \, \lesssim \, \frac{\rho'(x)}{x} + \frac{\rho(x)}{x^2}, \quad u \in I_l(x); 
	$$
	recall that $\rho(x) = \Im(z_+ - z_-) = 2\ka(x)$. In particular, we have
	$$
	|I_1(x)| \, \lesssim \, x^{\mu-1}\int_{-x}^{x} |f(u)| du
	$$
	and $I_1 \in \Lcal_{2,\ga - \mu}(I_\de)$, provided $f\in \Lcal_{2,\ga}(m)$. If $f$ is bounded, then clearly $I_1(x) = O(x^\mu)$.
	A similar estimate is true for $I_2(x)$, since in this case the kernel satisfies    
	$$
	\bigg|\frac{z_+'}{z_+ - q_+} - \frac{1}{x-u}\bigg| \, \lesssim \,  |\pa_x F(\ka)(x,u)| \, \lesssim \, x^{\mu - 1} 
	$$
	cf. estimate \eqref{aux:gac+}.

	We consider $I_3(x)$ next. In order to track the sign of the correction further down, we briefly note that the corresponding contributions over $\Ga^\pm$ read
	$$
	\pm\bigg(\frac{z_\pm'}{z_\pm - q_\mp} - \frac{1}{x-u}\bigg).
	$$
	Under our symmetry assumptions these kernels are (up to the negative sign) complex conjugates, hence the imaginary parts of the most singular kernels coincide, while real parts have opposite signs (the integral go over $f_\mp$). 
	
	From now on, we concentrate only on the $\Ga^+$ contribution. Using Lemma \ref{lem:errbasic} in the Appendix up to order $k=1$, we can write
	$$
	\frac{z_-'}{z_- - q_+} - \frac{1}{z_-- q_*}  = - \frac{i\rho'(x)}{z_-- q_*}  + O\bigg( x^{\mu - 1} + x^\mu\frac{\rho(x)}{(x -u)^2 + \rho(x)^2}\bigg)
	$$
	where we have set $q_*:= u + i \ka(x)$. In particular, we have 
	\be\label{aux:skernel}
	\aligned
	-\bigg(\frac{z_-'}{z_- - q_+} - \frac{1}{x-u}\bigg) = \frac{1}{x-u}&\frac{-i\rho(x)}{(x - u) - i\rho(x)} + \frac{i\rho'(x)}{(x -u) - i\rho(x)}  \\
	&+ O\bigg( x^{\mu - 1} + x^\mu\frac{\rho(x)}{(x -u)^2 + \rho(x)^2}\bigg).
	\endaligned
	\ee
	If $f\in \Lcal_{2,\ga}(m)$, the corresponding integrals over the error terms are bounded in $\Lcal_{2, \ga - \mu}(I_\de)$ (cf. the proof of Lemma \ref{lem:BRbasic}). When $f = O(1)$, the error terms contribute $O(x^\mu)$, since 
	$$
	\int \frac{\rho(x)}{(x -u)^2 + \rho(x)^2} \, du = O(1).
	$$
	In particular, it is enough to estimate the most singular kernel
	\be\label{eq:skernelS}
	\frac{1}{z_- - q_*} - \frac{1}{x - u} = \frac{1}{x-u}\frac{i\rho(x)}{(x - u) - i\rho(x)}
	\ee
	since the remaining kernel in \eqref{aux:skernel} can be written as
	$$
	\frac{i\rho'(x)}{z_- - q_*} =  i\rho'(x)\frac{1}{x -u}  + i\rho'(x)\Big(\frac{1}{z_- - q_*} - \frac{1}{x - u}\Big).
	$$
	To estimate \eqref{eq:skernelS}, we employ the  variable change $x = h(\xi), u = h(\tau)$ in the region 
	$$
	\xi\in \wt I_\de = h^{-1}((0,\de)), \quad \tau \in \wt I_c(\xi) = h^{-1}(I_c(x)),
	$$  
	(cf. Appendix for the definition and properties of $h$) which yields
	$$
	\int_{I_c(x)}f(u)\, \frac{1}{x-u} \frac{-i\rho(x)}{(x - u) - i\rho(x)}\,du = \int_{\wt I_c(\xi)}f(\tau)\frac{h'(\xi)}{h(\xi) - h(\tau)}\frac{-ih'(\tau)}{(h(\xi) - h(\tau)) + ih'(\xi) }\, d\tau
	$$
	(note the extra minus sign due to $h'<0$). Lemma \ref{lem:errvc} together with \eqref{eq:veps12h} from the Appendix imply the transformed kernel can be written as  
	\be\label{aux:errvc}
	\frac{h'(\xi)}{h(\xi) - h(\tau)}\frac{h'(\tau)}{(h(\xi) - h(\tau)) + ih'(\xi) } = \frac{1}{(\xi - \tau)}\frac{1}{(\xi - \tau) + i} + O\bigg( \frac{1}{\t m(\xi)} \Big(\frac{1}{\t m(\tau)} + \frac{1}{(\xi - \tau)^2 + 1} \Big)\bigg),                                                                                                                                                       
	\ee
	where $\t m(\tau) = 1 + \tau$. 
	
	Assume first $f\in H^1_{\ga + 1}(m)$. Then, we have 
	$$
	f\circ h \in \Lcal^1_{2, 1 - \t \ga}(\Rbb_+, \t m),
	$$ 
	where $\t \ga = \mu^{-1}(\ga + 1/2) + 1/2$. It is then not difficult to see that the integrals over the error terms belong to $\Lcal_{2, 1 - \t \ga}(\wt I_\de)$, provided $f\circ h \in \Lcal_{2, -\t\ga}(\t m)$ (cf. the corresponding part of the proof of Lemma \ref{lem:BRbasic} and note that 
	$$
	x^{2(\ga - \mu)} dx \, \leadsto \, h(\xi)^{2(\ga - \mu)} |h'(\xi)| d\xi, \qquad h(\xi)^{\ga - \mu} |h'(\xi)|^{1/2} \sim \t m(\xi)^{1 - \t\ga}
	$$
	implies that $\Lcal_{2, 1 - \t \ga}(\t m)$ corresponds to $\Lcal_{2, \ga - \mu}(m)$ in our original variables). To estimate the main term, we extend $f \circ h$ by zero to $\Rbb$, which to simplify the notation, we also denote by $f$. Moreover, we may assume without loss of generality  that  
	$$
	0 < (1 - \t \ga) + 1/2 < 1,
	$$
	which corresponds to $0 < \ga + 1/2 < \mu$. Otherwise, snce we are integrating over a region where $\xi \sim \tau$, we can correct by $\t m(\xi)^{k-1}/\t m(\tau)^{k-1}$ where $k \in \Nbb$ is such that $0 < (k - \t \ga) + 1/2 < 1$ (cf. the proof of Lemma \ref{lem:BRbasic}). Then, 
	$$
	\int_{\wt I_c(\xi)} f(\tau) \, \frac{1}{\xi -\tau}\frac{-i}{(\xi - \tau) +  i} d\tau = \int_{\Rbb} f(\tau) \, \frac{1}{\xi -\tau}\frac{-i}{(\xi - \tau) +  i} \, d\tau + \, \Lcal_{2, 1-\t \ga}(\wt I_\de),
	$$
	since (cf. \eqref{eq:veps12h} in the Appendix) we have
	$$
	\Big|\int_{\Rbb\setminus \wt I_c(\xi)} f(\tau) \, \frac{1}{(\xi - \tau)}\frac{1}{(\xi - \tau) + i}\, d\tau\Big|  \, \lesssim \,\int^\infty_{\xi}\frac{|f(\tau)|}{\t m(\tau)}\,\frac{d\tau}{\tau} + \frac{1}{\xi}\int^\xi_{0}\frac{|f(\tau)|}{\t m(\tau)}\,d\tau
	$$
	which is bounded in $\Lcal_{2, 1 - \t \ga}(\wt I_\de)$ by Hardy inequalities. 
	
	We claim 
	\be\label{aux:Ixi}
	I(\xi):= \int_{\Rbb} f(\tau) \, \frac{1}{\xi -\tau}\frac{-i}{(\xi - \tau) + i} \, d\tau +  i\pi  f(\xi) \, \in \Lcal_{2, 1-\t \ga}(\wt I_\de).
	\ee
	Indeed we have
	$$
	\frac{1}{\xi -\tau}\frac{-i}{(\xi - \tau)  + i} = \frac{1}{(\xi - \tau) + i} - \frac{1}{\xi -\tau}.
	$$
	and the real part of $I$ is the Hilbert transform of the imaginary part. In particular, it is enough to show 
	$$
	-\Im I(\xi) = \int_{\Rbb}   f(\tau) \, \frac{1}{(\xi - \tau)^2+ 1} \, d\tau - \pi  f(\xi)  \in \Lcal_{2, 1-\t \ga}(\wt I_\de).
	$$
	Taking the Fourier transform of this convolution, we obtain modulo constant factors
	$$
	\wh {\Im I}(s) \sim \wh {f}\, \big(e^{-|s|} - 1\big) \sim \wh {Hf'}(s) \, \frac{e^{-|s|} - 1}{|s|},
	$$
	hence taking the inverse Fourier transform we conclude
	$$
	\Im I(\xi) \sim \int_\Rbb (H f')(\tau) \log\bigg(1 + \frac{1}{(\xi - \tau)^2}\bigg)d\tau,
	$$
	which is bounded in $\Lcal_{2,1 - \t\ga}(\wt I_\de)$ as required. Indeed, the corresponding integral operator is bounded on $L^2(\Rbb)$ because its Fourier transform is a bounded function. As the kernel is $o(|\xi -\tau|)$, a suitable commutator with the weight function together with the $L^2(\Rbb)$-boundedness result shows the claim when $\xi \sim \tau$. For the remaining regions, proceed as in the first part of the proof of Lemma \ref{lem:BRbasic}.  
	
	In particular, we see that
	$$
	I^+_{3}(x) =- \frac{if_+(x)}{2} + \Lcal_{2, \ga - \mu}(I_{\de})
	$$
	where we have taken into account the normalization and we have used subscript $+$ to emphasize that we are integrating over $\Ga^+(x)$. As the corresponding part of the integral over $\Ga^-(x)$ satisfies an analogous estimate with  $i f_+(x)/2$ replaced by $i f_-(x)/2$, we conclude
	$$
	I^{+}_3(x) + I^{-}_3(x) = -\frac{i(f_+(x) + f_-(x))}{2}  + \Lcal_{2, \ga - \mu}(I_{\de}).
	$$
	
	The case $f\in H^1_\ga(m)$ can be reduced to the case $f\in H^1_{\ga + 1}(m)$. In fact, we may assume $f\in \Lcal^1_{2, \ga}(m)$ (constant factors $f_\pm(0)$ just contribute terms of order $O(x^\mu)$), in which case 
	$$
	\int_{I_c(x)}f(u)\, \frac{1}{x-u} \frac{-i\rho(x)}{(x - u) - i\rho(x)}\,du = x\int_{I_c(x)}\frac{f(u)}{u}\, \frac{1}{x-u} \frac{-i\rho(x)}{(x - u) - i\rho(x)}\,du + \Lcal_{2, \ga - \mu - 1}(I_\de) 
	$$
	Since $f(u)/u \in \Lcal^1_{2, \ga + 1}(I_{\de})$, we may proceed as before. 
	
	
	We next consider $I_4(x)$. We have 
	$$
	|k(x,u)| \, \lesssim \, \frac{x^\mu}{u}, \qquad |I_4(x)| \, \lesssim \, x^{\mu}\int_{x}^{2\de} \frac{|f(u)| + |f(-u)|}{u} du 
	$$
	and therefore $I_4\in \Lcal_{2,\ga - \mu}(I_{\de})$ by Hardy's inequality, whenever $f\in \Lcal_{2,\ga}(m)$. If we assume $f=O(1)$, we can write the kernel $k(x,u)$ as
	\be\label{aux1:k1}
	\aligned
	k(x,u) &= \frac{i\rho'(x)}{z_+ - q_+} + O\Big(\frac{\rho(x)}{u^2}\Big) \\
	&= - \frac{i\rho'(x)}{q_+} + \frac{i\rho'(x)}{q_+}\frac{z_+}{z_+ - q_+} + O\Big(\frac{x^{\mu + 1}}{u^2}\Big)\\
	&=-\frac{i\rho'(x)}{q_+} + O\Big(\frac{x^{\mu + 1}}{u^2}\Big).
	\endaligned
	\ee
	The integral over the error term gives the correct estimate $O(x^\mu)$. For the remaining term, we have
	$$
	\frac{1}{q} = \frac{1}{u} + O(u^{\mu - 1})
	$$
	and therefore 
	$$
	\int_{\Ga^+_r(x)\cup \Ga^-_r(x)} \frac{f(q)}{q}\, ds_q  = f_\pm(0)\int_{I_r(x)} \frac{du}{u} + O(1).
	$$
	However, the remaining integral is also bounded since the $\log$ contributions containing $x$ cancel out (the cusps are 'smoothly' connected). In particular, we have shown that 
	$$
	I_4(x)  =  O(x^{\mu}).
	$$
	
	Finally consider $I_5(x)$. Since by assumption $|z_\pm - q|$ is bounded away from zero by a constant depending on $\de$, the kernel satisfies the estimate 
	\be\label{eq:kfaraway}
	x\in (0, \de), \,\, q\in \Ga\setminus(\Ga^+ \cup \Ga^-) \quad \Rightarrow \quad |k(x,u)|\, \lesssim \, x^\mu
	\ee
	and the claim follows.
\end{proof}

In order to isolate the lower bound for $\si$ (and also for the upper bound), we now give precise estimates on the real and the imaginary parts of $\Dcal(f)$ when we control two derivatives of $f$.  

\begin{lemma}\label{lem:BRdiff2} Let $ 1 -\mu <\beta + 1/2 < \mu$. Then, $\Dcal(f)$ as defined in Lemma \ref{lem:BRdiff} satisfies
	$$
	f\in H^2_{\beta + 1}(m) \quad \Rightarrow \quad \Dcal(f)(x) = O(m(x)^\mu). 
	$$
	Furthermore, one can say more for the real part: actually, 	$$
	\Re \Dcal(f)(x) =  - \rho'(x) \frac{f_+(x) - f_-(x)}{2} + \Im ((z_+' - z_-')b_0(f)) + O(m(x)^{\mu + \lambda}), \quad \lambda := 1 - (\beta + 1/2) 
	$$
	where $\rho(x) := Im(z_+ - z_-)$. 
	
	Under additional assumptions on zeros of $f$, the estimate for the imaginary part can be improved as well. More precisely, we have
	$$
	f\in \Lcal^2_{2, \beta + 1}(m) \quad \Rightarrow \quad \Im \Dcal(f)(x) =  -\Re ((z_+' - z_-')b_0(f))  + O(m(x)^{\mu + \lambda}).
	$$
	Similarly, when
	$$
	f \in H^2_{\beta, 0}(m)  \quad \Rightarrow \quad \Dcal(f)(x) = -i(z_+' - z_-')b_0(f) + O(m(x)^{\mu + 1}),
	$$ 
	with the real part satisfying an improved estimate 
	$$ 
	\aligned
	\Re \Dcal(f)(x) = -\rho'(x)& \frac{f_+(x) - f_-(x)}{2} - \rho(x) \, \frac{f'_+(x) - f'_-(x)}{2}  \\
	& +  \Im \Big( (z_+' - z_-')b_0(f) + (z_+' z_+- z_-' z_-)b_1(f)\Big) + O(m(x)^{1 + \mu + \lambda}).
	\endaligned
	$$ 
	In addition, if the first order derivative of $f$ vanishes at the singular point as well, then
	$$
	f \in \Lcal^2_{2,\beta}(m) \quad \Rightarrow \quad \Im \Dcal(f)(x) =  - \Re \Big( (z_+' - z_-')b_0(f) + (z_+' z_+- z_-' z_-)b_1(f)\Big) + O(m(x)^{1 + \mu + \lambda}).
	$$
	The continuous linear functionals $b_i(f)$ have been defined in \eqref{def:bi}.
\end{lemma}

\begin{proof} Keeping the same notation, we retrace our steps in Lemma \ref{lem:BRdiff} and comment on the various integrals under current assumptions on $f$. Moreover, we refine the estimates on the real part of the kernel. In particular, we consider the integral \eqref{eq:Dker} with $x\in (0, \de)$, or more precisely the sum $\sum_{i=1}^5 I_i(x)$, defined as in \eqref{eq:Ix>b}. Let $k(x,u)$ denote the kernel \eqref{aux1:k}. 
	
	We start with $I_1(x)$. It is not difficult to see, that we actually have
	$$
	|\Re k(x,u)| \, \lesssim \, x^{2\mu-1}, \qquad \bigg|\Im \Big(\frac{z_\pm'}{z_\pm - q_+}\Big)\bigg| \, \lesssim \, x^{\mu-1}, \qquad u \in I_l(x)
	$$
	and all the statements are straightforward, noting that 
	$$
	f\in \Lcal^1_{2, \beta - k}(m) \quad \Rightarrow \quad f = O(m^{\lambda + k}), \qquad k = 0,1.
	$$
	We consider $I_2(x)$ and $I_3(x)$ next, where we integrate over $u\in I_c(x)$. The kernel in $I_2(x)$ satisfies 
	$$
	\Re \Big(\frac{z_+'}{z_+ - q_+}\Big) - \frac{1}{x-u}  = O( x^{2\mu - 1}), \quad \Im \Big(\frac{z_+'}{z_+ - q_+}\Big) = O(x^{\mu - 1})
	$$
	(cf. proof of Lemma \ref{lem:BRbasic}) and we conclude as for $I_1(x)$.
	
	As for the kernel in $I_3(x)$, we start with estimate \eqref{aux:skernel}, i.e.
	$$
	- \bigg(\frac{z_-'}{z_- - q_+} - \frac{1}{x - u}\bigg) = - \bigg( \frac{1}{z_- - q_*} - \frac{1}{x - u}\bigg) + \frac{i\rho'(x)}{z_- - q_*} + r(x,u),
	$$
	where $q_*:= u + i \ka(x)$ and  the remainder satisfies 
	$$
	r(x,u) = O\bigg(x^{\mu - 1} + x^{\mu} \frac{\rho(x)}{(x - u)^2 + \rho(x)^2}\bigg).
	$$
	In order to show the required statements on $\Re\Dcal(f)$, this estimate needs to be refined. In fact, further down we show that $\Re r(x,u)$ actually contributes $O(x^{\mu +  \lambda + k})$ and therefore all non-integral terms in the statement of the Lemma come from the first two kernels. The main term gives
	\be\label{eq:Icmain}
	\aligned
	\int_{I_c(x)} f(u)\, \frac{1}{x - u} \frac{-i\rho(x)}{(x - u) - i\rho(x)}du =& - 2(if(x) + f'(x)\rho(x))\arctan\Big(\frac{x}{2\rho(x)}\Big) \\
	&- i\rho(x)\int_{I_c(x)} \pa_x F(f)(x,u) \, \frac{x - u}{(x - u) - i\rho(x)}du, 
	\endaligned
	\ee
	where we have used 
	$$
	\int_{I_c(x)} \frac{1}{(x - u) - i\rho(x)}\,du = 2i \arctan \Big(\frac{x}{2\rho(x)}\Big), \qquad p.v.\int_{I_c(x)} \frac{1}{x - u}du = 0
	$$
 	and for the purposes of this Lemma we have set $\varepsilon = 1/2$ in the definition of $I_c(x)$ (and of $I_{r,l}(x)$). The last term in \eqref{eq:Icmain} contributes $O(x^{\mu + k + \lambda})$, since
	$$
	|\pa_x F(f)(x,u)| \, \lesssim \, x^{k -(\beta + 1)}\frac{1}{|x - u|^{1/2}}, \quad u \in I_c(x), \quad f''\in \Lcal_{2, \beta + (1 - k)}(I_\de), \quad  k=0,1.
	$$
	When $k=0$, the second term in \eqref{eq:Icmain} can also be absorbed in the error since  
	$$
	f \in H^2_{\beta + 1}(I_\de) \quad \Rightarrow \quad f'(x)\rho(x) = O(x^{\mu + \lambda}).
	$$
	
	Similarly, we conclude 
	\be\label{eq:Icmain2}
	i\rho'(x)\int_{I_c(x)} f(u) \, \frac{1}{z_- - q_*} du = -2 f(x)\rho'(x) \arctan\Big(\frac{x}{2\rho(x)}\Big) + O(x^{\mu + k + \lambda}).
	\ee
	Since  
	$$
	\arctan\Big(\frac{x}{2\rho(x)}\Big) = \frac{\pi}{2} + O(x^\mu),
	$$
	the claim for $\Re\Dcal(f)$ follows from \eqref{eq:Icmain} and \eqref{eq:Icmain2}, if we can show 
	\be\label{eq:Icremainder}
	\int_{I_c(x)} f(u) \Re r(x,u) du = O(x^{\mu +k+  \lambda}).
	\ee
	In fact, we can write $\Re r(x,u)$ as 
	$$
	\Re r(x,u)  =  \ka'(x) \Re \frac{\rho(x)}{(z_- - q_*)^2} -  \pa_x F(\ka)\Im \Big(\frac{x- u}{z_- - q_*}\Big)^2  + O\bigg(x^{2\mu - 1} + x^{2\mu} \frac{\rho(x)}{(x - u)^2 + \rho(x)^2}\bigg),
	$$
	(cf. Lemma \ref{lem:errbasic} from the Appendix), where the second kernel can be further rewritten as
	$$
	\pa_x F(\ka)\Im \Big(\frac{x- u}{z_- - q_*}\Big)^2 =  - \ka''(x)\rho(x)\frac{(x-u)^3}{|z_- - q_*|^4} + \pa^2_x F(\ka)\rho(x)\frac{(x-u)^4 }{|z_- - q_*|^4}.
	$$
	Using $\ka'''(x)\in H^1_{\beta + 2}(I_\de)$, we conclude
	$$
	\pa^2_x F(\ka)\rho(x)\frac{(x-u)^4 }{|z_- - q_*|^4} =  O\big(x^{\mu + \lambda -1}\big).
	$$
	In particular, the corresponding integral gives correct estimates (using the appropriate assumptions on the zeros of $f$). On the other hand, we have
	$$
	-\ka''(x)\rho(x)\int_{I_c(x)} f(u) \, \frac{(x-u)^3}{|z_- - q_*|^4} du =  \ka''(x)\rho(x)\int_{I_c(x)} \frac{f(x) - f(u)}{x - u} \frac{(x-u)^4}{|z_- - q_*|^4} \, du, 
	$$
	hence this terms contributes $O(m(x)^{2\mu + k})$ with $k = 0,1$. Finally, it is not difficult to see we have
	$$
	\ka'(x)\rho(x)\int_{I_c(x)} f(u) \Re\frac{1}{(z_- - q_*)^2} du = O(x^{\mu + k+ \lambda}),
	$$ 
	hence \eqref{eq:Icremainder} follows. 
	
	It remains to consider the far away contributions $I_4(x)$. To show the claim for $\Re I_4(x)$ when $f = O(x)$, respectively $\Im I_4(x)$ when $f = O(x^{1 + \lambda})$, we need second order corrections for the kernel. We can write
	$$
	\aligned
	k(x,u) &=  (z_+' - z_-')\Big(\frac{1}{q} + \frac{z_+}{q^2}\Big) + (z_+ - z_-)\frac{z_-'}{q^2} + r(x,u) \\
	&=  (z_+' - z_-')\, \frac{1}{q} + (z_+'z_+ - z_-'z_-) \, \frac{1}{q^2} + r(x,u)
	\endaligned
	$$
	where 
	$$
	r(x,u) = i\rho'(x) \Big(\frac{z_+}{q}\Big)^2 \frac{1}{z_+- q} + i\rho(x)\frac{z_-'}{z_- - q} \Big(\frac{z_+}{q}\frac{1}{z_+ - q} - \frac{z_-}{q^2}\Big). 
	$$
	In particular, we have  
	$$
	\Re r(x,u) = O(x^{\mu + 2}u^{\mu - 3}), \qquad \Im r(x,u) = O(x^{\mu + 2}u^{-3})
	$$
	since
	$$
	\Re \Big(\frac{z}{q}\Big)^k = O\Big(\frac{x^k}{u^{k}}\Big), \quad \Im \Big(\frac{z^k}{q^k}\Big) = O\Big(\frac{x^k u^\mu}{u^{k}}\Big), \quad \Re\Big(\frac{1}{z - q}\Big) = O(u^{-1}), \quad \Im \Big(\frac{1}{z - q}\Big) = O(u^{\mu-1})
	$$
	hence the integral over $I_r(x)$ contributes $O(x^{k + \lambda + \mu})$ using the appropriate assumptions on the zeros of $f$. 
\end{proof}



\begin{lemma} \label{lem:BRdiff_comm} Let $0< \ga + 1/2 < 1$ and let $\Dcal f$ be defined as in Lemma \ref{lem:BRdiff}. Then, the commutator   
	$$
	[g, \Dcal] f(x) := g(x)\Dcal(f)(x) - \Dcal(gf)(x), \quad x\in I_\de
	$$
	satisfies
	$$
	g\in H^2_{\ga}(m), \, \, f\in H^1_{\beta + 1}(m) \quad \Rightarrow \quad [g, \Dcal] f(x) = O(m(x)^{\mu}).
	$$
	Moreover, we have
	$$
	g \in H^2_{\ga, 0}(m), \, \, f\in H^1_\beta(m) \quad \Rightarrow \quad [g, \Dcal] f(x) = i(z_+' - z_-') b_0(fg)(x) + O(m(x)^{\mu + 1})
	$$
	with $b_0(fg)$ as defined in \eqref{def:bi}.
\end{lemma}

\begin{proof}
	It is enough to show the required estimate for the integral
	$$
	x \rightarrow \frac{1}{2\pi }\int_{\Ga} f(u)\Delta g(x, u)\bigg(\frac{z'_+}{z_+ - q} - \frac{z'_-}{z_- - q} \bigg) du, \quad x\in(0, \de),
	$$
	where we have set $\Delta g(x,u) := g(x) - g(u)$. We divide the integral as in Lemma \ref{lem:BRdiff} and, keeping the same notation for the various parts, we quickly indicate the necessary changes. Recall that  
	$$
	k(x,u) = \frac{z_+' - z_-'}{z_+ - q_+} + \frac{z_-'}{z_- - q_+}\frac{z_+ - z_-}{z_+ - q_+},
	$$
	cf. \eqref{aux1:k}. In particular, using $g' = O(1)$, we have
	$$
	\Delta g(x,u) k(x,u) = O(x^{\mu}), \quad u \in I_l(x) 
	$$
	and therefore $I_1(x) = O(x^{\mu + k})$ provided $f\in H^1_{\beta + (1 -k)}(m)$ with $k=0,1$. 
	
	On the other hand, when $u \in I_c(x)$, we similarly have
	$$
	\Delta g(x,u)\bigg(\frac{z_+'}{z_+ - q_+} - \frac{1}{x- u}\bigg) = O(x^\mu),
	$$
	while
	$$
	\aligned
	\Delta g(x,u)\bigg(\frac{z_-'}{z_- - q_+} - \frac{1}{x- u}\bigg) &= F(g) \frac{i\rho(x)}{(x - u) - i\rho(x)}  + O(x^{\mu})\\
	& = i\rho(x)g'(x)\frac{1}{(x - u) - i\rho(x)} + O\Big(x^{\mu + 1}\frac{x^{-\ga}}{|x - u|^{1/2}} + x^\mu\Big)
	\endaligned
	$$
	where we have used $g'' \in \Lcal_{2, \ga}(m)$ in order to arrive to the weakly singular kernel in the error term (cf. formula \eqref{eq:FTaylor} in the Appendix). The corresponding integrals over the error terms satisfy correct estimates and the same is true for the remaining term subtracting 
	$$
	i\rho(x)g'(x)f(x)\int_{I_c(x)} \frac{1}{(x - u) - i\rho(x)}du 
	$$ 
	and using that 
	$$
	|F(f)(x,u)| \, \lesssim \, \frac{x^{k-1-\beta}}{|x - u|^{1/2}}, \qquad u\in I_c(x), \qquad f \in H^1_{\beta + (1- k)}(m), \qquad k=0,1.
	$$
	In particular, the corresponding integrals contribute
	$$
	I_2(x) + I_3(x) = O(x^{\mu + k}).
	$$

	Finally, let us consider the far-away contributions to the integral, i.e. $u\in I_r(x)$. When $f\in \Lcal_{2, \beta}(m)$, the claim is straightforward since
	$$
	\Delta g(x,u) k(x,u) =  O(x^{\mu}).
	$$
	On the other hand, when $f = O(1)$ and $g = O(m)$, we have $fg = O(m)$ and the proof of Lemmas \ref{lem:BRdiff} and \ref{lem:BRdiff2} imply
	$$
	I_4(x) = -\rho'(x)b_0(fg)(x) + O(g(x)x^\mu + x^{\mu + 1})
	$$
	(in this case, there is no need to use the extra cancellation coming from $\Delta g$, as opposed to the case $f\in \Lcal_{2, \beta}(m)$).
	
\end{proof}


\begin{corollary} \label{lem:BRn}
	Let $f$ be odd with respect to the $x$-axis. Then
	$$
	\aligned
	f\in H^1_{\beta + 1}(m) \quad &\Rightarrow \quad BR(z, f)\cdot z_s^\perp \in \Lcal_{2, \beta -\mu}(m), \\
	f\in H^1_{\beta}(m) \quad &\Rightarrow \quad BR(z, f)\cdot z_s^\perp = O(m^\mu) + \Lcal_{2, \beta -\mu - 1}(m)
	\endaligned
	$$
	and
	$$
	f\in H^2_{\beta + 1}(m) \quad \Rightarrow \quad BR(z, f)\cdot z_s^\perp = O(m^\mu).
	$$
	When $f$ is even with respect to the $x$-axis, analogous statements are true for $f + 2BR(z, f)\cdot z_s$.

	Moreover, when $f$ is odd with respect to both axis, we have 
	\be\label{BRzerosOdd}
	f\in H^2_{\beta}(m)  \quad \Rightarrow \quad BR(z, f)\cdot z_s^\perp = O(m^{\mu + 1}) 
	\ee
	and similarly, when $f$ is even with respect to both axes, we have 
	\be\label{BRzerosEven}
	f\in \Lcal^2_{2, \beta}(m) \quad \Rightarrow \quad f + 2BR(z, f)\cdot z_s = O(m^{\mu + 1}).
	\ee
\end{corollary}

\begin{proof}The claim follows from Lemmas \ref{lem:BRdiff} and \ref{lem:BRdiff2}. In fact, if $f$ is e.g. odd with respect to the $x$-axis, then $BR(z, f)\cdot z_s^\perp$ is even. In particular, passing to the graph parametrization and taking into account the change of orientation on the lower branch, we see that
	$$ 
	BR(z, f)\cdot z_s^\perp = \frac{1}{2|z'(x)|}\, \Re \Dcal f(x).
	$$
	If, in addition, $f$ is odd with respect to the $y$-axis, we have $f \in H^2_{\beta, 0}(m)$ (as $f(\pm \al_*) = 0$ for all times) but also $b_0(f) = 0$. 
\end{proof}

\begin{remark}\label{rem:BRt}
 Let $f\in H^2_{\beta,0}(m)$ be odd with respect to to both axes. Then
 $$
 BR(z, f) \cdot z_s^\perp = \Im(z_s z)(-f + \Re b_1(f)) + O(m^{1+\mu +\lambda})
 $$
 by Lemma \ref{lem:BRdiff2}. Moreover, it is not difficult to see that we also have  
 $$
 \frac{f}{2} + BR(z, f)\cdot z_s = -\Re(z_s z ) (-f + \Re b_1(f)) + O(m^{1+\mu}). 
 $$
\end{remark}




\section{The local existence theorem domain with cusp}
	\label{s.estEnergy}
	
	We consider the system of evolution equations   
	\be\label{eqs:sys}
	\aligned
	&(\log |z_\al|)_t =  \frac{\om_s}{2} + BR(z,\om)_s\cdot z_s, \\
	&\theta_t = \frac{\om \theta_s}{2}  + BR(z,\om)_s \cdot z_s^\perp, \\
	&\frac{\om_t}{2} + BR(z,\om)_t \cdot z_s = 0,
	\endaligned
	\ee
	for $(\log|z_\al|, \theta, \om)$, where $\om$ is odd and $\log|z_\al|$ is even with respect to both transformations $\al \leftrightarrow - \al$ and $\al_* - \al \leftrightarrow \al_* + \al$ (by symmetry, we must have $\al_* = \pi/2$), while $\theta$ satisfies \eqref{def:sym-xy-theta}.

	To recover the parametrization of the interface from $(\theta, \log|z_\al|)$, we fix the constant of integration 
	\be\label{def:z-int-constant}
	z(\al_*, t) = 0, \qquad z_t(\al_*, t) = 0
	\ee
	(as we will see these will be consistent with the regularity assumptions on the vorticity) and we set
	\be\label{def:z}
	z(\al, t) :=  \int_{\al_*}^\al |z_\al(\al',t)| e^{i\theta(\al',t)}d\al'.
	\ee
	(symmetry-wise, the components satisfy \eqref{def:sym-x}--\eqref{def:sym-y} as required).

	Let $\mu\in (1/2, 1]$ and let $\beta\in (-1/2, 1/2)$ be such that 
	\be\label{eq:betaI}
	\aligned
	\mu \in (1/2 , 2/3] \quad &\Rightarrow \quad 1 - \mu < \beta + 1/2 < \mu,\\
	\mu \in (2/3, 1] \quad &\Rightarrow \quad 1 - \mu < \beta + 1/2 < \mu/2.
	\endaligned
	\ee   
	In both cases, we have
	\be\label{eq:betaII}
	\beta + 1/2 < 1 - \mu/2.
	\ee
	For $k\geq 2$, we say $(z(\cdot, t), \om(\cdot, t))$ belong to the Banach space $\Bcal^k_{\beta, \mu}(m)$, if the above symmetry assumptions are satisfied, we have $z(\al_*, t) = z(-\al_*, t) = 0$ and 
	\be\label{regularity-z}
	\theta(\cdot, t)\in H^{k+1}_{\beta + k}(\Tbb, m), \qquad \log|z_\al(\cdot, t)| \in  H^{k + 1/2}_{\beta - \mu/2 + (k-1) + 1/2}(\Tbb, m),
	\ee
	respectively
	\be\label{regularity-om}
	\om(\cdot, t)\in H^{k + 3/2}_{\beta + (k-1) + 1/2, 0}(\Tbb, m),
	\ee
	where the weight $m$ is a fully symmetric, non-negative function, smooth everywhere except at $\al = \pm \al_*$, such that $m(\al) \sim |\al\pm \al_*|$ in the neighborhood of $\pm \al_*$ and $m(\al) \sim 1$ otherwise. 
	
	We now give additional conditions which define a particular open set $\Ocal^k_{\beta,\mu}(m)\subseteq \Bcal^k_{\beta, \mu}(m)$ in which we will construct the solutions of the system \eqref{eqs:sys}. First, we assume additional regularity on  
	\be\label{def:phi-om}
	\varphi_s = \frac{\om_s}{2} + BR(z,\om)_s\cdot z_s.
	\ee
	More precisely, we require
	\be\label{eq:varphis-reg}
	\varphi_s(\cdot, t) \in H^{k + 1/2}_{\beta - \mu/2 + (k-1) + 1/2}(\Tbb, m), 
	\ee
	i.e. the derivatives of $\varphi_s$ grow by a factor of $O(m^{\mu/2})$ slower than the corresponding derivatives of $\om_{s}$ (respectively $\om\theta_s$). For later use, we recall that     
	\be\label{eq:varphi_st}
	\varphi_{st} =- \varphi_s^2 + \theta_t^2 - \si \theta_s.
	\ee	
	
	We only consider curves which satisfy following version of the arc-chord condition: there exists some small $r>0$, such that the parametrization satisfies
	\be\label{def:ac}
	\|z\|_\Fcal := \Fcal_r(z) + \sup_{\al\in B_r(\al_*)}\frac{1}{|z_{1\al}(\al, t)|} + \sup_{\al\in B_r(\al_*)} \frac{1}{m(\al)^{-\mu}|z_{2\al}(\al, t)|} , \qquad \|z\|_\Fcal < \infty
	\ee
	(where $\Fcal_r(z)$ has been defined in \eqref{arcchord}), together with higher order estimates 
	\be\label{def:ac2}
	m^{j-\mu} \pa_s^j \theta(\cdot, t) <\infty, \qquad 1 \leq j\leq k-1.
	\ee
	Finally, we assume the normal component of the pressure gradient $\si = \pa_n P$ satisfies the Rayleigh-Taylor condition, which in the current setting reads 
	\be\label{RTcondition}
	m^{-(\mu + 1)}\si(\cdot, t) > 0,
	\ee    
	with $\si$ given by 
	\be\label{def:si}
	\si = \frac{\om \theta_t}{2} + BR(z, \om)_t \cdot z_s^\perp. 
	\ee

	We are now able to make the statement of Theorem \ref{T.main} precise; however as we are only interested in showing there exist some solutions of Euler equations that exhibit cusp singularities, we opt to restrict the class of initial data for which we show local existence. More precisely, we have 
	
	\begin{theorem}\label{thm:T.main} Let $k\geq 2$ and let $(z^0, \om^0)\in \Ocal^k_{\beta,\mu}(m)$ be such, that 
	\be\label{eq:om-add-reg}
	\om^0 \in H^{k + 3/2}_{\beta -\mu/2 + (k-1) + 1/2,0}(\Tbb, m),
	\ee
	where $\Ocal^k_{\beta,\mu}(m)$ is the open subset of $ \Bcal^k_{\beta, \mu}(m)$ defined via \eqref{eq:varphis-reg}, \eqref{def:ac}, \eqref{def:ac2} and the Rayleigh-Taylor condition \eqref{RTcondition}. Then, there exists a time $T>0$ and $(z, \om)\in\Ccal([0, T], \Ocal^k_{\beta, \mu}(m)) \cap \Ccal^1([0, T]), \Ocal^{k-1}_{\beta, \mu}(m))$ 
	solution of \eqref{eqs:sys} up to time $T$ such that $z(\cdot, 0) = z^0$ and $\om(\cdot, 0) = \om^0$.
	\end{theorem}	

	\begin{remark}
	The additional assumption \eqref{eq:om-add-reg} on the regularity of the vorticity is chosen for convenience in the regularization part of the proof. In fact, in that case mollified initial data clearly also belong to $\Ocal^k_{\beta,\mu}(m)$, which is a-priori not assured when the vorticity satisfies \eqref{regularity-om} and \eqref{eq:varphis-reg} at the same time (cf. Section \ref{ss.existence} below for more details).

	\end{remark}
	
	\begin{remark}\label{rem:initialdata}
	There are indeed initial data $(z,\om)\in\Bcal^k_{\beta,\mu}(m)$ that satisfy \eqref{RTcondition}. First, it is not difficult to verify that the Rayleigh-Taylor condition \eqref{RTcondition} holds if we initially have e.g.
	$$
	\om_s(\al_*) = 0, \qquad b_1(\om) = \frac{1}{2\pi i}\int_\Ga \frac{\om(\al')}{z(\al')^2}\, ds_{\al'}  \neq 0
	$$ 
	(where the imaginary part of $b_1(\om)$ vanishes by symmetry). In fact, the corresponding solution $\om_t\in H^2_{\beta,0}(m)$ of \eqref{eq:omt} must satisfy 
		$$
		\om_{ts}(\al_*) - \Re b_1(\om_t) = \Im b_1(\om\theta_t) - \Re b_1(z_t D_s\om) - (\Re b_1(\om))^2.
		$$
		Setting $b:=\Re b_1(\om)$, we have 
		$$
		 \aligned
		\si &= \Im(z_s z) \left(-\om_{ts}(\al_*) + \Re b_1(\om_t) + \Im b_1(\om\theta_t) - \Re b_1(z_tD_s\om)\right) - z\cdot z_s^\perp \, b^2 + O(m^{1 + \mu + \lambda}) \\ 
		&= b^2\left(\Im(z_s z) - z\cdot z_s^\perp \right) + O(m^{1 + \mu + \lambda})\\
		&= b^2 \, z \cdot z_s \,\Im (z^2_s)  + O(m^{1 + \mu + \lambda}).
		\endaligned
		$$ 
	where $\lambda = 1 - (\beta + 1/2)$. In particular, in the neighborhood of the origin, we have
	$$
	m^{-(1 + \mu)}\si \, \gtrsim \, b^2 .
	$$
	We do this calculation in more detail in Lemma \ref{lem:corner-preliminary}, Section \ref{s.Corner} when we consider domains with corners; however here we use Remark \ref{rem:BRt} and we take into account that in general $\om_{ts}(\al_*)\neq 0$. If $\om_s(\al_*)\neq 0$, a similar result can be shown using Lemma \ref{lem:BRdiff2}. We omit the details.
		\end{remark}
	
	We first show that (time derivatives of) various quantities needed for the energy estimates belong to correct weighted Sobolev spaces with corresponding norms controlled by some power of 
	$$
	\|\om\|_{H^{k + 1}_{\beta + (k-1)}(m)} + \|\varphi_s\|_{H^k_{\beta - \mu/2 + (k-1)}(m)} + \|\theta_s\|_{H^{k}_{\beta + k}(m)} + \sup_{1\leq j<k}|m^{j-\mu} \pa_s^j \theta| + \| z \|_\Fcal
	$$ 
	We then define the full energy functional and prove the corresponding a-priori energy estimate.


\subsection{Preliminary estimates in the Lagrangian parametrization}
\label{s:estDerivatives}

\begin{lemma}\label{lem:omsss} Let $\varphi_s \in H^{k}_{\beta + (k-1)}(m)$ where $k\geq 2$ and let $\om \in H^2_{\beta}(m)$. Then, we have
	$$
	\om \in H^{k + 1}_{\beta + (k-1)}(m).
	$$ 
\end{lemma}

\begin{proof} 
	The second order derivative of the vorticity $\om_{ss}$ satisfies the equation   
	$$
	\frac{\om_{ss}}{2} + BR(z,\om_{ss})\cdot z_s = F, \qquad F := BR(z,\om_{ss})\cdot z_s - (BR(z,\om)_s \cdot z_s)_s + \varphi_{ss},
	$$
	where $F\in H^1_{\beta + 1}(m)$. Indeed, $\varphi_{ss}$ clearly satisfies this condition. On the other hand, we have
	$$
	\aligned
	(BR(z,\om)_s&\cdot z_s)_s  - BR(z,\om_{ss})\cdot z_s  = \Re \bigg(\frac{1}{2\pi i}\int_\Ga \bigg(\frac{z_s(\al)^2}{z_{s}(\al')^2} - 1\bigg)\om_{ss}(\al') \, \frac{z_s(\al)}{z(\al) - z(\al')}ds_{\al'}\bigg) + \\
	&+ \Re \bigg(\frac{1}{2\pi i}\int_\Ga \frac{z_s(\al)^2}{z_{s}(\al')^2}\Big((z_{s}^2 D^2_{s}\om)(\al') -\om_{ss}(\al')\Big) \, \frac{z_s(\al)}{z(\al) - z(\al')}ds_{\al'}\bigg) + 2\theta_s BR(z,\om)_s \cdot z_s^{\perp}
	\endaligned
	$$
	hence we control the derivative of the first integral in $\Lcal_{2,\beta + 1}(m)$, using the $\Lcal_{2,\beta}(m)$-norm of $\om_{ss}$ only (cf. Lemma \ref{lem:BR-additional-cancelation}). Recall also that we use the notation $D_s\om = \pa_s(\om/z_s)$. To show the claim for the second integral note that 
	$$
	f := z_{s}^2 D^2_{s}\om -\om_{ss} = -2 \om \theta_s^2 - i \om \theta_{ss} - 3i\theta_s\om_s \in H^1_{\beta + 1}(m),
	$$
	hence, integrating by parts, we obtain
	$$
	BR(z, f)^*_s = z_s BR_{+1} (z, D_s f)^* \in \Lcal_{2, \beta + 1}(m)
	$$ 
	cf. \eqref{def:BR+k}. Finally, it is not difficult to see that we also have $\theta_s BR_s\cdot z_s^\perp\in H^1_{\beta + 1}(m)$. In particular, we can solve the Neumann problem
	$$
	\frac{\phi}{2} + BR(z,\phi)\cdot z_s = F, \qquad F\in H^1_{\beta + 1}(m)
	$$
	for $\phi\in H^1_{\beta + 1}(m)$ even with respect to both axes (the solution exists by Theorem \ref{thm:inverseOpI}) and by uniqueness, we must have $\phi \equiv \om_{ss}$ (cf. Proposition \ref{prop:injectivity}).
	
	To obtain higher order derivatives of the vorticity, we consider the following
	$$
	\frac{z_1 \om_{sss}}{2} - BR(z, z_1 \om_{sss})\cdot z_s = z_1 F_1, \quad F_1 = -((BR(z,\om_{ss})\cdot z_s)_s + \frac{1}{z_1}BR(z, z_1 \om_{sss})\cdot z_s) + F_{s} +  \varphi_{sss}.
	$$
	It is not difficult to see that $F_1\in H^1_{\beta + 2}(m)$ (where we must use $BR_{+1}$ when taking a derivative of the first $BR$ integral) and we can proceed as before to find $\phi \equiv z_1\om_{sss}\in H^1_{\beta + 1}(m)$. But then, we must have $\om_{sss} =  \phi/z_1\in H^1_{\beta + 2}(m)$.  
\end{proof}

For convenience, we state and prove all the remaining results for $k = 2$, the generalization to higher $k$ being straightforward.

\begin{lemma}\label{lem:param}
	The time derivative of the parametrization satisfies 
	\be\label{est:zt}
	z_{t} \in H^3_{\beta + 1}(m), \qquad z_{1t} = O(m), \qquad z_{2t} = O(m^{1 + \mu}) 
	\ee
	respectively
	\be\label{eq:ztdiff}
	\sup_{\al,\al'}\bigg|\frac{z_t(\al, t) - z_t(\al', t)}{z(\al, t) - z(\al', t)}\bigg| = O(1).
	\ee
	Moreover, we have $z_{ts}^* z_s = \varphi_s - i\theta_t$, where
	$$
	\theta_t\in \Lcal^2_{2,\beta + 1}(m), \qquad \theta_t  = O(m^\mu), \qquad \theta_{ts}  = O(m^{\mu - 1}).
	$$
	and 
	\be\label{eq:th_ts-Hphi_ss}
	\theta_{ts} = H(\varphi_{ss})  + \Rcal, \qquad   \Rcal\in H^2_{\beta + 2}(m).
	\ee
\end{lemma}

\begin{proof} 
	We first claim  
	\be\label{est:tht_mu}
	\theta_t = O(m^\mu).
	\ee
	Indeed, it is not difficult to see that in complex notation we can write
	\be\label{eq:complex-phi_s-th_t}
	z_{ts}^* = \left(\frac{D_s\om}{2} + z_sBR(z, D_s\om)^*\right) = \frac{\varphi_s - i\theta_t}{z_s},
	\ee
	where recall that $D_s\om := \pa_s(\frac{\om}{z_s})$. In particular, we have
	$$
	\theta_{t}= z_{1s} \Big(-\frac{\Im(D_s\om)}{2} + BR(z, D_s\om)\cdot z_s^\perp\Big) - z_{2s}\Big(\frac{\Re (D_s\om)}{2} + BR(z, D_s\om)\cdot z_s\Big).
	$$
	Since  
	$$
	-\frac{\Im(D_s\om)}{2} + BR(D_s\om)\cdot z_s^\perp = BR(\Re D_s\om)\cdot z_s^\perp - \Big(\frac{\Im D_s\om}{2} + BR(\Im D_s\om)\cdot z_s\Big) 
	$$
	is $O(m^\mu)$ by Lemma \ref{lem:BRn} (note that $\Re(D_s\om) = O(1)$ is odd, while $\Im(D_s\om) = O(m^\mu)$ is even with respect to the $x$-axis and they both belong to $H^2_{\beta +1}(m)$), the estimate \eqref{est:tht_mu} follows using that $z_{2s} = O(m^\mu)$.

In order to show \eqref{eq:th_ts-Hphi_ss}, we take a $D_s$ derivative of \eqref{eq:complex-phi_s-th_t}, which yields
\be\label{eq:complex-phi_s-th_t-der}
\frac{\varphi_{ss} - i\theta_{ts}}{z^2_s} -2i\theta_s\frac{\varphi_s-i\theta_t}{z_s^2}= \frac{D^2_s\om}{2} + z_sBR(z, D^2_s\om)^*.
\ee
In particular, it is not difficult to see that in the notation of Lemma \ref{lem:struct1} we have
$$
\Rcal := \Acal(z_s^2, D^2_s\om) - 2\theta_s \varphi_s - 2H(\theta_s \theta_t) \in H^2_{\beta + 2}(m).
$$
It remains to show the additional regularity for $\theta_{ts}$, that is
\be\label{est:thts_mu2}
\theta_{ts} = O(m^{\mu - 1}).
\ee
First note that taking the real part of \eqref{eq:complex-phi_s-th_t-der}, then using $\varphi_{ss}\in\Lcal^1_{2,\beta-\mu/2+1}(m)$ together with $\theta_t\in\Lcal^2_{2,\beta + 1}(m)$ and $z_{2s} = O(m^\mu)$, we see that 
$$
\frac{\Re(D^2_s\om)}{2} + BR(z, D^2_s\om)\cdot z_s \in \Lcal^1_{2,\beta-\mu/2 + 1}(m)
$$
(note that $\Re(D^2_s\om)$ is odd, while $\Im(D^2_s\om)$ is even with respect to~both axes and they both belong to $H^1_{\beta + 1}(m)$). In particular, we can solve  
$$
\frac{\phi}{2} + BR(z,\phi)\cdot z_s = \frac{\Re(D^2_s\om)}{2} + BR(z, D^2_s\om)\cdot z_s
$$
for $\phi\in \Lcal^1_{2,\beta - \mu/2 + 1}(m)$ odd with respect to both axes (cf. Theorem \ref{thm:inverseOpI}, Section \ref{s:inverseOp}). By uniqueness and symmetry, we must also have 
$$
\frac{\Im(D^2_s\om)}{2} - BR(z, D^2_s\om)\cdot z_s^\perp = BR(z,\phi)\cdot z_s^\perp \in  \Lcal^1_{2,\beta - \mu/2 + 1}(m),
$$
which implies the right-hand side of \eqref{eq:complex-phi_s-th_t-der} belongs to $\Lcal^1_{2,\beta - \mu/2 + 1}(m)$ and therefore by Lemma \ref{lem:sobolev} must grow as $O(m^{\mu/2-(\beta + 1/2)})$. In particular, since $\beta + 1/2$ satisfies  \eqref{eq:betaII}, we have \eqref{est:thts_mu2}.

Finally, Lemma \ref{lem:BRdersII} and Corollary \ref{lem:BRn} imply
$$
z_t \cdot z_s  = O(m), \qquad   z_t \cdot z_s^\perp = O(m^{\mu + 1}).
$$
In order to prove \eqref{eq:ztdiff}, it is enough to consider $\al \in B_\epsilon (\al_*)$ and $\al'\in B_\epsilon(-\al_*)$ (the remaining combinations are easily seen to be bounded since $z_{t\al} = O(1)$). There, we have 
\be\label{eq:zt1aux}
z_t(\al) - z_t(\al') = z_t(\al) - z_t(-\al) + O(\al - |\al'|)
\ee
where symmetry with respect to the $x$-axis implies
\be\label{eq:zt2aux}
z_{t}(\al) - z_t(-\al) =  2i z_{2t}(\al) = O(m(\al)^{\mu + 1}).
\ee
In particular, we have
$$
\frac{z_t(\al) - z_t(\al')}{z(\al) - z(\al')} = \frac{2i z_{2t}(\al)}{z(\al) - z(\al')} + O(1)
$$
The claim now follows using Lemma \ref{lem:errorrk} in the Appendix; we omit further details.
\end{proof}


We now give an auxiliary Lemma, which we will be used to estimate $\si$ and for the existence of $\om_t$.

\begin{lemma}\label{lem:BRtime}
We have
\be\label{eq:diff-BRt-BR(t)}
f\in H^{2}_{\beta, 0}(m) \quad \Rightarrow \quad BR(z, f)_t^* - BR(z,f_t)^* \in H^{2}_{\beta, 0}(m).
\ee
When $f$ is odd with respect to both axes, we have
\be\label{est:bRt_asym}
 \big[BR(z, f)_t - BR(z,f_t)\big]\cdot z_s^\perp = O(m^{\mu + 1}).
\ee
If, in addition, we have $f\in H^3_{\beta + 1}(m)$, then 
$$
[BR(z, f)_t - BR(z,f_t)]_s\cdot z_s^\perp = O(m^\mu) + \Lcal_{2, \beta - \mu - 1}(m).
$$ 
\end{lemma}

\begin{proof}
We have (neglecting the time-dependence and recalling that $\varphi_s = \frac{|z_\al|_t}{|z_\al|}$)
$$
\aligned
BR(z,f)_t^* - BR(z,f_t)^* =& - \frac{1}{2\pi i}\int f(\al') \frac{z_t(\al) - z_t(\al ')}{(z(\al) - z(\al '))^2} \,ds_{\al '} + BR(z, \varphi_s f)^*\\
=& z_t BR(z, D_{s}f)^* - BR\big(z, D_{s}(f z_t)\big)^*  + BR(z, \varphi_s f)^* . 
\endaligned
$$
Using 
$$
D_s(f z_t) = z_t D_s f +  f \frac{z_{ts}}{z_{s}}, \qquad  \frac{z_{ts}}{z_s} = \varphi_s + i \theta_t, \qquad D_s f = \pa_s\left(\frac{f}{z_s}\right)
$$
we conclude 
\be\label{eq:BRt}
BR(z,f)_t^* - BR(z,f_t)^*  =  [z_t, BR(z, \cdot)^*] D_{s}f  - i BR(z, f \theta_t)^*.
\ee
We first consider the Birkhoff-Rott integral of $f \theta_t$. Lemma \ref{lem:param} (together with $f = O(m)$) implies  
$$
f \theta_t \in \Lcal^2_{2,\beta}(m), \qquad f\theta_t  = O(m^{\mu + 1}).
$$
In particular, Lemma \ref{lem:BRdersII} combined with Corollary \ref{lem:BRn} (when $f$ is odd with respect to both axes; in which case $f\theta_t$ is even with respect to both axes) implies 
$$
BR(z, f\theta_t)^* \in H^2_\beta(m), \qquad BR(z, f\theta_t)^* = O(m), \qquad  \frac{f \theta_t}{2} + BR(z, f \theta_t )\cdot z_s = O(m^{\mu + 1}).
$$
The claim for the derivative follows from the appropriate estimates in Corollary \ref{lem:BRn} (cf. proof of Lemma \ref{lem:param}). The details are standard.

As for the commutator 
\be\label{eq:zt_comm}
[z_t, BR(z, \cdot)^*] D_{s}f = z_t BR(z, D_s f)^* - BR(z, z_t D_s f)^*
\ee
we first show it belongs to $H^2_\beta(m)$. Indeed, taking the first derivative, we get
$$
\aligned
\pa_s \big([z_t, BR(z, \cdot)^*] D_{s}f\big) &= z_s [z_t, BR(z, \cdot)^*] D^2_s f + i z_s[\theta_t, BR(z, \cdot)^*]D_s f +z_s[\varphi_s, BR(z, \cdot)^*]D_s f \\
&=z_s [z_t, BR(z, \cdot)^*] D^2_s f + H^1_{\beta}(m).
\endaligned
$$
Using \eqref{eq:ztdiff} to estimate the first term, we conclude 
$$
\pa_s \big([z_t, BR(z, \cdot)^*] D_{s}f\big) \in L^\infty.
$$
We can take one more derivative, since  
$$
 \pa_s [z_t, BR(z, \cdot)^*] D^2_s f = \frac{1}{2\pi i}\int_\Ga D_s^2f(\al')(z_t(\al) - z_t(\al')) \pa_s\Big(\frac{1}{z(\al) - z(\al')}\Big) \, ds_{\al '} + z_{ts}BR(z, D^2_s f)^*
$$
belongs to $\Lcal_{2, \beta}(m)$.  Indeed the claim is straightforward, except for the first integral when e.g. $\al \in B_\epsilon (\al_*)$ and $\al'\in B_\epsilon(-\al_*)$. However, in that case, the most singular contribution to the kernel can be estimated as in the proof of Lemma \ref{lem:BRbasic}, combining Lemma \ref{lem:errorrk} from the Appendix with estimates \eqref{eq:zt1aux}-\eqref{eq:zt2aux} for $z_t(\al) - z_t(\al')$. We omit further details.

We now show the required asymptotic estimate for the normal component of the commutator, i.e. we show
\be\label{eq:ztcommn}
\Re \big(iz_s[z_t, BR(z, \cdot)^*]D_s f \big) = O(m^{\mu + 1})
\ee
under the assumption that $f$ is odd with respect to both axes. For later use, we write out the normal component of the commutator 
\be\label{eq:ztcommwritten}
\aligned
\Re \big(iz_s[z_t, BR(z, \cdot)^*]D_s f \big) = z_{1t} &BR(z, D_s f)\cdot z_s^\perp -z_{2t}BR(z, D_sf)\cdot z_s\\
& - BR(z, z_{1t}D_s f)\cdot z_s^\perp +  BR(z, z_{2t} D_s f )\cdot z_s,  
\endaligned
\ee 
where $\Re (z_t D_s f) = O(m)$ is odd and $\Im (z_t D_s f) = O(m^{\mu + 1})$ is even w.r.t both axes. In fact, we have $z_{1t} = O(m)$ resp. $z_{2t} = O(m^{\mu + 1})$ by Lemma \ref{lem:param} (as noted previously $\Re D_s f = O(1)$ is odd with respect to the $x$-axis, but even with respect to the $y$-axis and vice-versa for $\Im D_s f = O(m^\mu)$).

If $f\in H^3_{\beta + 1,0}(m)$, we can prove \eqref{eq:ztcommn} directly for each term using growth estimates on $z_t$ and appropriate statements from Corollary \ref{lem:BRn}. The claim for the derivative follows similarly. We omit the details and prove \eqref{eq:ztcommn} when $f\in H^2_{\beta,0}(m)$. First note that the above symmetry considerations imply terms with $z_{1t}$ can be written in the form required by Lemma \ref{lem:BRdiff_comm}. In particular, we have 
$$
z_{1t} BR(z, D_s f)\cdot z_s^\perp - BR(z, z_{1t} D_s f)\cdot z_s^\perp =  O(m^{\mu + 1})
$$
and clearly also 
$$
z_{2t}BR(z, D_sf)\cdot z_s = O(m^{\mu + 1}).
$$
As for the remaining term in \eqref{eq:ztcommwritten}, we can write 
$$
D_s f = f_s(\al_*)\frac{1}{z_s} + \frac{(f_s(\al) - f_s(\al_*)) + if\theta_s}{z_s},
$$
where 
$$
\frac{1}{z_s} \in H^2_{\beta + 1}(m), \qquad \wt {D_s f}:=D_sf - f_s(\al_*)\frac{1}{z_s} \in \Lcal^1_{2,\beta}(m).
$$
Since $z_{2t}/z_s$ satisfies assumptions of Corollary \ref{lem:BRn}, a \eqref{eq:ztcommn}-type estimate follows. 

It remains to show
$$
BR(z,z_{2t} \wt {D_s f})\cdot z_s = O(m^{\mu + 1}).
$$
Indeed, we have
$$
g := z_{2t}   \wt {D_s f} \in \Lcal^1_{2, \beta - \mu - 1}(m), \qquad D_s g \in \Lcal_{2, \beta - \mu - 1}(m) 
$$
and therefore
$$ 
BR(z, g)^* = BR_{-3}(z, g)^* - \frac{1}{2\pi i} \int_{\Ga} g(\al') \Big(\frac{1}{z(\al')} + \frac{z(\al)}{z(\al')^2} + \frac{z(\al)^2}{z(\al')^3} \Big) ds_{\al'},
$$
where only the imaginary part of the second correction integral is non-zero (all the others vanish by symmetry, since real/imaginary parts of $g$ are even/odd with respect to both axes). More precisely, we have 
$$
BR(z, g)^*  = -iz(\al)\, \Im\Big(\frac{1}{2\pi i} \int_{\Ga} g(\al')\, \frac{1}{z(\al')^2}\, ds_{\al'}\Big) + \Lcal^1_{2, \beta - \mu - 1}(m),
$$
where we have used that
$$
BR_{-3}(z, g)^*\in \Lcal_{2, \beta - \mu - 2}(m), \qquad \pa_s BR_{-3}(z, g)^* = z_sBR_{-2}(z, D_s g)^* \in \Lcal_{2, \beta - \mu - 1}(m)
$$
by Corollary \ref{lem:BRcorrections}. In particular, we have $BR_{-3}(z, g)^* = O(m^{\mu + 1})$ by Lemma \ref{lem:sobolev}, and the claim follows taking the tangential component.

We can proceed similarly to prove  
$$
[z_t, BR(z, \cdot)^*]D_s f = O(m).
$$
In fact, when $f\in H^3_{\beta + 1,0}(m)$, this follows from Lemma \ref{lem:BRdersII}. If we only control two derivatives of $f$, then we can show 
$$
[z_t, BR(z, \cdot)^*] \wt {D_s f} = O(m)
$$
as above, using one correction term less as we only have $z_t \wt {D_s f} \in \Lcal^1_{2,\beta -1}(m)$. We omit the details.
 \end{proof}

\begin{lemma}\label{lem:omt} 
There exists a unique $\om_t\in H^2_{\beta}(m)$, odd with respect to both axes, solution of 
\be\label{eq:omt}
\om_t  + 2 BR(z, \om_t)\cdot z_s = F,
\ee
where the RHS is given by  
\be\label{eq:def-F}
F := - 2 \big[BR(z,\om)_t - BR(z, \om_t)\big]\cdot z_s.
\ee
\end{lemma}

\begin{proof}

We construct $\om_t$ in several steps. First note that $F\in {H^2_{\beta, 0}}(m)$ is odd with respect to both axes and $F \equiv F(z,\om)$ doesn't depend on time derivatives of $\om$ (cf. Lemma \ref{lem:BRtime}). In particular, by Theorem \ref{thm:inverseOpI} below, there exists a unique $\om_t \in \Lcal^1_{2,\beta}(m)$ odd with respect to to both axes solution of \eqref{eq:omt}. Its derivative $\om_{ts}$ must satisfy  
$$
\om_{ts} - 2BR(z, \om_{ts})\cdot z_s = F_s - F_1,
$$
where 
\be\label{def:F1}
F_1 := 2\big[(BR(z, \om_t)\cdot z_s)_s + BR(z, \om_{ts})\cdot z_s\big].
\ee

We claim that given $\om_t\in \Lcal^1_{2,\beta}(m)$, we have $F_1\in H^1_{\beta+1}(m)$. 
Indeed, we can write
\be\label{def:F1-in-detail}
\aligned
\big[BR(z, \om_t)_s + BR(z, \om_{ts})\big]\cdot z_s &= \Re \bigg(\frac{1}{2\pi i}\int \Big(\frac{z_s(\al)}{z_s(\al')} + 1\Big)\om_{ts}(\al')\frac{z_s(\al)}{z(\al)-z(\al')}ds_{\al'}\bigg)\\ 
&+\Re \bigg(\frac{1}{2\pi i}\int \frac{z_s(\al)}{z_s(\al')}\big((z_sD_s\om_t)(\al') -\om_{ts}(\al')\big)\frac{z_s(\al)}{z(\al)-z(\al')}ds_{\al'}\bigg). 
\endaligned
\ee
The second integral clearly belongs to $H^1_{\beta + 1}(m)$, since
$$
z_sD_s\om_t -\om_{ts} = -i\om_{t}\theta_s \in H^1_{\beta + 1 - \lambda}(m) \subseteq H^1_{\beta+1}(m)
$$
where $\lambda = 1 - (\beta + 1/2)$. The same is true for the first integral by Lemma \ref{lem:BR-additional-cancelation}. In particular (by Theorem \ref{thm:inverseOpI} below) there exists $\phi\in H^1_{\beta + 1}(m)$ even with respect to both axes, solution of  
\be\label{eq:omts}
\phi - 2BR(z, \phi)\cdot z_s = \wt F. 
\ee
By uniqueness, we must have $\om_{ts} = \phi$ and therefore $\om_t\in H^2_{\beta + 1}(m)$, cf. Proposition \ref{prop:injectivity} below. However, having control over an additional derivative of $\om_t$ implies $F_1\in H^1_{\beta + 1 -\lambda}(m)$, cf. Lemma \ref{lem:BR-additional-cancelation}, which in turn implies $\om_t\in H^2_{\beta + 1 - \lambda}(m)$ by Theorem \ref{thm:inverseOpI}. 

This argument can be iterated until $\ga := \beta + 1 - k \lambda$ satisfies $\ga + 1/2 < 1$ for some $k\in\Nbb$, i.e. until we have $F_1\in H^1_{\ga}(m)$. If $\ga \leq \beta$, then we are finished applying Theorem \ref{thm:inverseOpI} one more time to conclude $\om_t\in H^2_\beta(m)$. Otherwise, we only have $\om_t\in H^2_\ga(m)$. However $(\ga - \lambda) + 1/2 = (\beta + 1 - (k + 1)\lambda) + 1/2 < \beta + 1/2$ (and clearly also $(\ga  -\mu) + 1/2 < \beta + 1/2$) and we claim that $F_{1s}\in \Lcal_{2,\beta}(m)$. In fact, when $\om_{t}\in H^2_{\ga}(m)$, $\om_t=O(m)$ we have 
$$
\aligned
F_{1s} &= \Re \bigg(\frac{1}{2\pi i}\int \Big(\frac{z_s(\al)}{z_s(\al')} + 1\Big)D_s\om_{ts}(\al')\frac{z_s(\al)}{z(\al)-z(\al')}ds_{\al'}\bigg) + \Lcal_{2,\beta}(m)\\ 
&=\Re \bigg(\frac{1}{2\pi i}\int \Big(\frac{z_s(\al)}{z_s(\al')} + 1\Big)(z_{1s}\om_{tss})(\al')\frac{z_s(\al)}{z(\al)-z(\al')}ds_{\al'}\bigg) + \Lcal_{2,\beta}(m). 
\endaligned
$$
hence the claim follows as in Lemma \ref{lem:BR-additional-cancelation}. We omit the details. In particular, we conclude $\om_t\in H^2_{\beta}(m)$ by another application of Theorem \ref{thm:inverseOpI}.
\end{proof}

\begin{lemma}\label{lem:si} 
The normal derivative of the pressure $\si$ satisfies 
\be\label{est:si_asym}
\si \in \Lcal^{3}_{2, \beta + 1}(m), \qquad \si = O(m^{\mu + 1}), \qquad \si_s = O(m^\mu) + \Lcal_{2, \beta-\mu -1}(m).
\ee 
\end{lemma}

\begin{proof}
We first show the asymptotic part of \eqref{est:si_asym}. Recall that $\si$ is given by Equation \eqref{eq:si}. We have
\be\label{eq:BRtn}
BR(z,\om)_t\cdot z_s^\perp = BR(z, \om_t)\cdot z_s^\perp + \big[BR(z, \om)_t - BR(z,\om_t)\big]\cdot z_s^\perp.
\ee
Lemma \ref{lem:BRdersII} and Corollary \ref{lem:BRn} applied to $f=\om_t$ imply
$$
BR(z, \om_t)\cdot z_s^\perp\in H^2_{\beta}(m), \qquad BR(z, \om_t)\cdot z_s^\perp = O(m^{\mu + 1}),
$$
while Lemma \ref{lem:BRtime} applied to $f = \om$ implies 
$$
[BR(z, \om)_t - BR(z,\om_t)]\cdot z_s^\perp \in H^{2}_{\beta}(m), \qquad [BR(z, \om)_t - BR(z,\om_t)]\cdot z_s^\perp = O(m^{\mu + 1}).
$$
Moreover, their derivatives are $O(m^\mu) + \Lcal_{2, \beta-\mu -1}(m)$. In particular, all the statements from \eqref{est:si_asym} follow except the claim for $\pa_s^3\si$. 

As we do not control $\pa_s^3\om_t$ and $\pa_s^3 \theta_t$, we need to use the `cancellations' proven in Lemma \ref{lem:BRcancel} to gain  control over an additional derivative of $\si$. More precisely, from \eqref{eq:si} we subtract the Hilbert transform of \eqref{eq:om_t} corresponding to the tangential derivative of the pressure (cf. proof of Lemma \ref{lem:param}, where Hilbert transform of $\varphi_{s} := (\log |z_\al|)_t$ is substracted from $\theta_t$). Indeed, we have
$$
\aligned
\frac{\om\theta_t}{2} +  BR(z, \om)_t\cdot z_s^\perp  &= BR(z, \om_t)\cdot z_s^\perp + \Big(\frac{\om\theta_t}{2} + BR(z,\om\theta_t)\cdot z_s\Big) + \Re \big( i z_s[z_t, BR(z, \cdot)^*]D_s \om\big)\\
  \frac{\om_t}{2} + BR(z, \om)_t \cdot z_s &= \Big(\frac{\om_t}{2} + BR(z,\om_t)\cdot z_s\Big) - BR(z,\om\theta_t)\cdot z_s^\perp + \Re \big( z_s[z_t, BR(z, \cdot)^*]D_s \om\big)
\endaligned
$$
cf. Lemma \ref{lem:BRtime}, where
$$
[z_t, BR(z, \cdot)^*]D_s \om = z_ t \Big( \frac{D_s\om}{2z_s} + BR(z, D_s\om)^*\Big) - \Big( \frac{z_t D_s\om}{2z_s} + BR(z, z_t D_s\om)^*\Big). 
$$
Combining all the previous Lemmas of this section, we have 
$$
\om_t \in H^2_{\beta,0}(m), \quad \om \theta_t \in \Lcal^2_{2,\beta}(m), \quad z_t D_s\om \in H^2_{\beta,0}(m).
$$ 
Note that
$$
\pa_s^2\Big(\frac{f}{2z_s} + BR(z, f)^* \Big) = z^2_s \Big(\frac{D^2_s f}{2z_s} + BR(z, D^2_s f)^*\Big) + i\theta_s z_s  \Big(\frac{D_s f}{2z_s} + BR(z, D_s f)^*\Big)
$$
and therefore, when $f\in H^2_{\beta,0}(m)$, we have
$$
\pa_s^2\Big(\frac{f}{2} + BR(z, f)^* z_s \Big) = z^2_s \Big(\frac{D^2_s f}{2} + BR(z, D^2_s f)^*z_s\Big)  + H^1_{\beta + 1}(m). 
$$

Some care is needed when considering the first term of the commutator with $z_t$. First note
$$
\pa^2_s  \Big(\frac{D_s\om}{2} + z_s BR(z, D_s\om)^*\Big) = z_s^2 \Big(\frac{D^3_s \om}{2} + BR_{+1}(z, D^3_s \om)^*z_s\Big) + H^1_{\beta + 2}(m)
$$ 
and therefore 
$$
\pa^2_s  \Big(z_t \Big(\frac{D_s\om}{2} + z_s BR(z, D_s\om)^*\Big)\Big) =   z_s^2 z_t \Big(\frac{D^3_s \om}{2} + BR_{+1}(z, D^3_s \om)^*z_s\Big) + H^1_{\beta + 1}(m)
$$
cf. \eqref{def:BR+k} for the definition of $BR_{+1}(z,\om)$. In particular, the claim now follows from Lemma \ref{lem:struct1} since 
$$ 
\aligned
\si_{ss} &= \si_{ss} - H\Big(\pa_{s}^2\Big(\frac{\om_t}{2} + BR(z, \om)_t \cdot z_s\Big)\Big) \\
&=\Acal\big(z_s^2, D^2_s (\om_t - i \theta_t\om - z_t D_s \om)\big) +  \Acal\Big(\frac{z_s^2z_t}{z}, z D_s^3\om\Big) + H^1_{\beta + 1}(m) 
\endaligned
$$
where we use Lemma \ref{lem:trHilbert} from the Appendix to estimate the Hilbert transform of a function in $H^1_{\beta + 1}(m)$.
\end{proof}

\begin{lemma}\label{lem:phi_t}
The time derivative of $\varphi_s$ 
satisfies
$$
 \varphi_{st}\in L^\infty, \qquad \varphi_{sst}\in \Lcal_{2, \beta - \mu/2}(m), \qquad \varphi_{ssst}\in \Lcal_{2, \beta - \mu/2 + 1}(m).
$$ 
 \end{lemma}

\begin{proof} First note that Lemmas \ref{lem:param} and \ref{lem:si} imply $\varphi_{st}\in L^\infty$. On the other hand, we clearly have
$$
\varphi_{sst} = -3\varphi_{s} \varphi_{ss} - (\si \theta_s)_s + \theta_t \theta_{ts} \in \Lcal_{2, \beta - \mu/2}(m)
$$ 
provided $\varphi_{ss}\in \Lcal_{2, \beta - \mu/2}(m)$. Finally, we clearly have 
$$
\varphi_{ssst} = -3\varphi_{ss} \varphi_{ss} - 4 \varphi_s \varphi_{sss}  - (\si \theta_s)_{ss} + (\theta_t \theta_{ts})_s  \in \Lcal_{2, \beta - \mu/2 + 1}(m).
$$
\end{proof}

In the next few Lemmas we consider second order time derivatives of $z$ and $\om$, with the goal of showing that $\si_t$ has zeros of the same order as $\si$ at the singular point, i.e. that $\si_t = O(m^{\mu + 1})$.

\begin{lemma}\label{lem:paramt}
We have 
\be\label{eq:ztt-formula}
z_{tt} = \si z_s^\perp.
\ee
In particular, 
$$
z_{tt}\in H^3_{\beta + 1}(m), \qquad z_{1tt} = O(m^{2\mu + 1}), \qquad z_{2tt} = O(m^{\mu + 1}).
$$
Similarly,  
\be\label{eq:thtt-formula}
\theta_{tt} = \si_s - 2\varphi_s \theta_t
\ee
and therefore
$$
\theta_{tt}	\in H^2_{\beta + 1}(m), \qquad \theta_{tt} = O(m^\mu) + \Lcal_{2, \beta - \mu -1}(m).
$$
\end{lemma}

\begin{proof}
Taking a time derivative of the $z_t$ equation, we see that
$$
z_{tt} = BR(z,\om)_t  + \frac{\om_t}{2} z_s + \frac{\om \theta_t}{2} z_s^\perp = \si z_s^\perp,
$$
where we have used Equation \eqref{eq:omt}. In particular, all the results on the regularity of $z_{tt}$ follow from the corresponding results for $\si$ and properties of $\theta_s$. The same is true for $\theta_{tt}$, since
$$
z_{tts} = z_{tst} + \varphi_s z_{ts} = (\varphi_{st} + i \theta_{tt} ) z_s + (\varphi_s + i \theta_t)^2z_s 
$$  
where we have used that $z_{ts} = (\varphi_s + i\theta_t)z_s$. In particular,
$$
\theta_{tt} +  2\varphi_s \theta_t = z_{tts}\cdot z_s^\perp  = \si_s.
$$
\end{proof}

We now state an auxiliary Lemma, whose proof is postponed to the end of the present section, then give its consequences for the construction of $\om_{tt}$ and regularity of $\si_t$.

\begin{lemma}\label{lem:BRtime2}
We have 
$$
[BR(z,\om)^*_{tt} - BR(z, \om_t)^*_t] \in H^2_\beta(m).
$$
where
\be\label{est:BRtt_asym}
[BR(z, f)_{tt} - BR(z, f_t)_t]\cdot z_s = O(m), \quad  \frac{(f\theta_t)_t}{2} + [BR(z, f)_{tt} - BR(z, f_t)_t]\cdot z_s^\perp = O(m^{\mu + 1}).
\ee
\end{lemma}

\begin{lemma}\label{lem:omtt} There exists a unique $\om_{tt}\in H^2_\beta(m)$ odd with respect to both axes, solution of 
\be\label{eq:omtt}
\om_{tt}  + 2BR(z, \om_{tt}) \cdot z_s = G, 
\ee
where 
$$
G := - 2[(BR(z, \om_{t})\cdot z_s)_t - BR(z, \om_{tt})\cdot z_s] + F_t 
$$
and $F$ as in \eqref{eq:omt}. 
\end{lemma}

\begin{proof}
It suffices to show that $G \in H^2_\beta(m)$. The existence of $\om_{tt}\in H^2_\beta(m)$ then follows as in Lemma \ref{lem:omt}. 
	Indeed, Lemma \ref{lem:BRtime} applied to $f=\om_t\in H^2_{\beta}(m)$, together with Lemmas \ref{lem:BRn} and \ref{lem:param} implies 
$$
(BR(z, \om_{t})\cdot z_s)_t - BR(z, \om_{tt})\cdot z_s = \theta_t BR(z, \om_t)\cdot z_s^\perp + [BR(z, \om_t)_t - BR(z, \om_t)]\cdot z_s
$$
belongs to $H^2_\beta(m)$. On the other hand, we have
\begin{align*}
\pa_t\big([BR(z,\om)_t - BR(z, \om_t)]\cdot z_s\big) = [BR(z, \om)_{tt} - & BR(z, \om_t)_t]  \cdot z_s  \\
 & + \theta_t[BR(z, \om)_t - BR(z, \om_t)]  \cdot z^\perp_s, 
\end{align*}
hence Lemma \ref{lem:BRtime2} applied to $f=\om\in H^2_\beta(m)$ imply
$$
 \pa_t \big(BR(z,\om)_t\cdot z_s - BR(z, \om_t)\cdot z_s\big) \in H^2_\beta(m).
$$
In particular, we have $G\in H^2_\beta(m)$.
\end{proof}

\begin{lemma}\label{lem:si_t} We have 
\be\label{est:sit_asym}
\si_t \in H^2_\beta(m), \quad \quad \si_t = O(m^{\mu + 1}). 
\ee
\end{lemma}

\begin{proof}

Taking a time derivative of \eqref{eq:si} we obtain
$$
\aligned
 \si_t &= \frac{(\om\theta_t)_t}{2} + BR(z,\om)_{tt} \cdot z_s^\perp - \theta_t BR(z,\om)_t\cdot z_s    \\
 &=  \frac{(\om\theta_t)_t}{2} + \frac{\om_t \theta_{t}}{2} + BR(z,\om)_{tt} \cdot z_s^\perp  
\endaligned
$$
where we have used that $\om_t + 2BR(z, \om)_t \cdot z_s = 0$.  We can further write 
$$
BR(z,\om)_{tt}\cdot z_s^\perp=BR(z,\om_{tt})\cdot z_s^\perp + [BR(z,\om_t)_{t} - BR(z,\om_{tt})]\cdot z_s^\perp + [BR(z,\om)_{tt} - BR(z,\om_t)_t]\cdot z_s^\perp 
$$
where Corollary \ref{lem:BRn} and Lemma \ref{lem:BRtime} imply
$$
BR(z,\om_{tt})\cdot z_s^\perp = O(m^{\mu + 1}), \qquad \frac{\om_t \theta_{t}}{2} + [BR(z,\om_t)_{t} - BR(z,\om_{tt})]\cdot z_s^\perp = O(m^{\mu + 1})
$$
(Lemma \ref{lem:omtt} implies $\om_{tt}\in H^2_{\beta,0}(m)$, while $\om_t \theta_t \in \Lcal^2_{2, \beta}(m)$ by Lemmas \ref{lem:param} and \ref{lem:omt}). Moreover, Lemma \ref{lem:BRtime2} implies
$$
\frac{(\om\theta_t)_t}{2} + [BR(z,\om)_{tt} - BR(z,\om_t)_t]\cdot z_s^\perp  = O(m^{\mu + 1}).
$$
In particular, \eqref{est:sit_asym} follows. The statement for the derivatives follows from the same Lemmas.
\end{proof}

\begin{proof}[Proof of Lemma \ref{lem:BRtime2}.]
	Equation \eqref{eq:BRt} implies
	\be\label{eq:BRtt}
	[BR(z, \om)^*_{t} - BR(z, \om_t)^*]_t  = ([z_t, BR(z, \cdot)^*] D_{s}\om)_t - i BR(z, \om \theta_t)^*_t,
	\ee
	where the time derivative of the commutator reads
	$$
	\aligned
	([z_t, BR(z, \cdot)^*] D_{s}\om)_t &=[z_{tt}, BR(z, \cdot)^*] D_{s}\om + [z_t, BR(z, \cdot)^*] (D_{s}\om)_t + [z_t, BR(z, \cdot)^*](\varphi_s D_s \om) \\
	&-\frac{1}{2\pi i}\int D_{s}\om(\al') \Big(\frac{z_t(\al) - z_t(\al')}{z(\al) - z(\al')}\Big)^2 \, ds_{\al '}.
	\endaligned
	$$
	Note that we can integrate by parts to get
	$$
	\frac{1}{2\pi i}\int D_{s}\om(\al') \Big(\frac{z_t(\al) - z_t(\al')}{z(\al) - z(\al')}\Big)^2 \, ds_{\al '} =-\big[z_t, [z_t, BR(z, \cdot)^*]\big] D^2_{s}\om + 2[z_t, BR(z, \cdot)^*]\Big(\frac{z_{ts}}{z_{s}}D_{s}\om \Big)
	$$
	and
	\be\label{eq:dcomm_zt}
	\aligned
	\big[z_t, [z_t, BR(z, \cdot)^*]\big] D^2_{s}\om = \frac{1}{2\pi i}\int D^2_{s}\om(\al') \frac{(z_t(\al) - z_t(\al'))^2}{z(\al) - z(\al')} \, ds_{\al '} 
	\endaligned
	\ee
	Taking into account $(D_s \om)_t = D_s \om_t - i D_s (\om \theta_t ) - D_s (\varphi_s \om) - \varphi_s D_s \om$, we finally obtain
	$$
	([z_t, BR(z, \cdot)^*] D_s \om)_t = [z_t, [z_t, BR(z, \cdot)^*]] D^2_s \om + [z_{tt}, BR(z, \cdot)^*] D_{s}\om  + [z_{t}, BR(z, \cdot)^*] G 
	$$
	where  
	$$
	G := 2(\varphi_s + i \theta_t)D_s f + D_s( (\varphi_s + i \theta_t)f) + D_s f_t \in H^1_\beta(m).
	$$ 
	
	It is not difficult to see that 
	$$
	[z_{tt}, BR(z, \cdot)^*] D_{s}\om\in H^2_\beta(m), \qquad [z_t, BR(z, \cdot)^*] G \in H^2_\beta(m)
	$$
	both satisfying estimates as in~\eqref{est:BRtt_asym}. In fact, comparing Lemmas \ref{lem:param} and \ref{lem:paramt} we see that $z_{tt}$ satisfies all the properties required from $z_t$ when proving analogous claims for the commutator in Lemma \ref{lem:BRtime} and the same is true for $G$ and $D_s \om$ (they satisfy the same symmetry properties and we can write $G$ as a sum of a function in $H^2_{\beta + 1}(m)$ and a function in $\Lcal^1_{2, \beta}(m)$). In particular, we can proceed as in the proof of Lemma \ref{lem:BRtime}. We omit further details.

	It remains to consider the `double' commutator \eqref{eq:dcomm_zt}. Taking a derivative, we obtain
	$$
	\pa_s \big[z_t, [z_t, BR(z, \cdot)^*]\big] D^2_{s}\om =\frac{z_s}{2\pi i}\int_{\Ga} D^2_{s}\om(\al') \bigg( \frac{2z_{ts}}{z_s} \frac{z_t(\al) - z_t(\al')}{z(\al) - z(\al')} -\Big(\frac{z_t(\al) - z_t(\al')}{z(\al) - z(\al')}\Big)^2\bigg) ds_{\al'}
	$$
	which belongs to $H^1_\beta(m)$. 
	
	In order to show \eqref{est:BRtt_asym}-type estimates for the double commutator, we write 
	$$
	\big[z_t, [z_t, BR(z, \cdot)^*]\big] D^2_{s}\om  = z_t [z_t, BR(z, \cdot) ^*] D^2_s \om - [z_t, BR(z, \cdot) ^*] (z_t D^2_s \om), 
	$$
	where the second commutator can be estimated as the commutator term in Lemma \ref{lem:BRtime}. As for the first, estimate \eqref{eq:ztdiff} implies 
	$$
	[z_t, BR(z, \cdot) ^*] D^2_s \om =O(1)
	$$ 
	and it only remains to prove 
	\be\label{aux:doublecomm}
	\Re ( iz_s [z_t, BR(z, \cdot) ^*] D^2_s \om) = O(m^\mu).
	\ee
	However, terms with $z_{1t}$ satisfy the desired estimate by Lemma \ref{lem:BRdiff_comm}, i.e. we have 
	$$
	z_{1t} BR(z , D_s^2 \om) \cdot z_s^\perp - BR(z , z_{1t} D_s^2 \om) \cdot z_s^\perp = O(m^\mu) 
	$$
	(since $z_{1t}$ is even with respect to the $x$-axis, while $D_s^2 \om \in H^1_{\beta + 1}(m)$ with the real part odd and the imaginary part even with respect to both axis). On the other hand, we have $z_{2t} D_s^2 \om \in \Lcal^1_{2, \beta - \mu}(m)$, hence we can write
	$$
	BR(z, z_{2t} D_s^2 \om)^* = BR_{-1}(z, z_{2t} D_s^2 \om)^* - \frac{1}{2\pi i}\int_{\Ga} \frac{z_t D_s^2 \om(\al')}{z(\al')} \, ds_{\al'}, 
	$$
	where the real part of the correction term vanishes by symmetry. Since $D_s(z_t D_s^2 \om)$ has even real and odd imaginary part with respect to both axes, we actually have
	$$
	BR_{-1}(z, z_{2t} D_s^2 \om)^* \in \Lcal^1_{2, \beta - \mu}(m), 
	$$
	and therefore, it is not difficult to see that 
	$$
	-z_{2t} BR(z, D_s^2 \om) \cdot z_s = O(m^\mu), \quad  BR(z, z_{2t} D_s^2 \om) \cdot z_s = O(m^\mu)
	$$
	by Lemma \ref{lem:sobolev}. In particular, \eqref{aux:doublecomm} follows.

	Finally, we consider the second term in \eqref{eq:BRtt}. We can write   
	$$
	BR(z, \om\theta_t )_t^* - BR(z, (\om\theta_t)_t)^* =  -i BR(z, \theta_t^2 \om)^* + [z_t, BR(z, \cdot)^*] D_s(\om\theta_t).
	$$
	Since $D_s(\om\theta_t) \in \Lcal^1_{2, \beta}(m)$ with $\om\theta_t$ even (when $\om$ is odd) it is not difficult to see the commutator term can be estimated as the commutator term in Lemma \ref{lem:BRtime}.  On the other hand, we have $(\om\theta_t)_t$ even with respect to both axes, with Lemmas \ref{lem:param} and \ref{lem:paramt} implying $(\om\theta_t)_{t} \in \Lcal^2_{2, \beta}(m)$, which in turn implies that we can apply Corollary \ref{lem:BRn} and Lemma \ref{lem:BRdersII} to
	$$
	\frac{(\om\theta_t)_t}{2z_s} + BR(z, (\om\theta_t)_t)^*.
	$$ 
	The same is true for the Birkhoff-Rott integral of $\om\theta_t^2$. 
\end{proof}


\subsection{The a priori energy estimate}

To simplify the notation, we set 
\be\label{def:beta'}
\beta' := \beta - \mu/2 + 1.
\ee
The lower-order contributions to the energy read
\be\label{def:energy-low}
E_{l}(t)^2 :=  \|\theta\|^2_{H^2_{\beta + 1}(m)} +  \|\log |z_\al|\|^2_{H^1_{\beta'}(m)} + \|\om\|^2_{H^2_{\beta}(m)} + \|z\|^2_\Fcal, 
\ee
with higher order contributions given by 
\be\label{def:energy-high}
\aligned
E_{h}^{k}(t)^2 := &\|\sqrt{\si}\pa_s^{k+1}\theta\|^2_{2, \beta' +  (k - 3/2)} + \|\Lambda^{1/2}(m^{\lambda + (k-2)}\pa^k_s\varphi_s)\|^2_{2, \beta '-\lambda + (k -3/2)} \\ 
&+ \|\pa_s^k \varphi_s\|^2_{2,\beta'+ (k-2)} + \|\pa_s^k\log |z_\al|\|^2_{2, \beta' + (k-2)} + \sup_{\al\in[-\pi, \pi]} |m^{(k-1)-\mu}\pa_s^{k-1}\theta|,
\endaligned
\ee
For $k\geq 2$, the energy functional is defined to be
$$
E_k(t)^2 = E_l(t)^2 + \sum_{i = 2}^k E^i_h(t)^2.
$$
It generalizes the unweighted energy functional used when the interface satisfies the arc-chord condition, i.e. when $\sup_{\al,\beta\in [-\pi,\pi]}\Fcal(z)(\al, \beta) < \infty$ with $\Fcal(z)$ as defined in \eqref{eq:defFcal} (cf. \cite{AM}, \cite{CCG}). Note that we consider the interface with respect to the Lagrangian parametrization as opposed to the arc-length parametrization used there. When the Rayleigh-Taylor condition \eqref{RTcondition} is satisfied we have $\pa_s^{k+1}\theta \in \Lcal_{2, \beta + k}(m)$ as required. Also note that the norm of $\pa_s^k\om_s$ can be replaced by  the norm of the corresponding $\pa_s^k\varphi_s$ by Lemma \ref{lem:omsss}. 

Let us comment on the value of $\lambda$ in the definition of the fractional derivative. As discussed in Section \ref{ss.weightedSobolev}, we require that both $(\beta ' - \lambda) + 1/2$ and $(\beta ' - \lambda) + 1$ satisfy Muckenhaupt condition, i.e. that
$$
-1/2 < (\beta ' - \lambda) + 1/2 < 0,
$$ 
Depending on the value of $\mu$, we distinguish between two cases: 
\begin{itemize}
\item When $\mu \in (1/2, 2/3]$, bounds for $\beta + 1/2$ imply $\beta' + 1/2$ has values in the interval $(1, 4/3)$. In this case $1 < \lambda < 2$.
\item When $\mu \in (2/3, 1]$, we have $\lambda<1$, since $\beta' + 1/2$ has values in the interval $(1/2, 1)$. 
\end{itemize}

We are now ready to prove:

\begin{lemma}\label{lem:apriori} Let $(z,\om)$ be a sufficiently regular solution of \eqref{eqs:sys} so that  
$$
a(t):=\min_{\al \in [-\pi, \pi]}m(\al)^{-(\mu + 1)}\si(\al, t)>0
$$ 
and let $k\geq 2$. Then, the energy functional $E_k(t)$ satisfies the a-priori energy estimate
\be\label{eq:energy}
\frac{d}{dt} E_k(t)^p \, \lesssim \, \frac{1}{a(t)} \exp C E_k(t)^p
\ee
for some $p\in \Nbb$ and some constant $C>0$.
\end{lemma}

%
%
%
%
\begin{proof}
We prove the claim for $k=2$; the proof for higher~$k$ is completely analogous. We only need to consider time derivatives of the highest order terms in $\theta$ and $\varphi_{s}$, all of the remaining terms follow from the corresponding Lemmas in Section \ref{s:estDerivatives}. We first show 
\be\label{en1}
\frac{1}{2}\frac{d}{dt} \big\|\sqrt{\sigma}\pa_s^3\theta\big\|^2_{2, \beta ' + 1/2} = \int_\Ga m^{2(\beta ' + 1/2)} \si \pa_s^3\theta \, \pa_s H(\pa_s^2\varphi_{s}) ds +  \Rcal_1(t),
\ee
with the remainder $\Rcal_1(t)$ bounded by (at most) the exponential of some power of the energy. Indeed, we have
$$
\frac{1}{2}\frac{d}{dt} \big\|\sqrt{\sigma}\pa_s^3\theta\big\|^2_{2, \beta ' + 1/2} = \int_\Ga m^{2(\beta ' + 1/2)} \si \pa_s^3\theta \,  (\pa_s^3\theta)_{t}  \, ds + \int_\Ga m^{2(\beta ' + 1/2)} \Big(\frac{\si_t}{\si} + \varphi_s\Big) \, \si|\pa_s^3\theta|^2 ds. 
$$
The last two terms clearly satisfy the desired estimate, since
$$
\frac{\si_t}{\si} = O(1), \qquad \varphi_s = O(1).
$$
As for $(\pa_s^3\theta)_{t}$, successively interchanging $\pa_t$ with $\pa_s = \frac{1}{|z_\al|}\pa_\al$, then using Lemma \ref{lem:param}, we obtain 
$$
\aligned
(\pa_s^3\theta)_{t} = \pa^3_s\theta_{t} - \pa_s^2(\varphi_s \theta_{s}) - \pa_s(\varphi_s \pa_s^2\theta) - \varphi_s \pa_s^3\theta =\pa^2_s H\varphi_{ss} + \Lcal_{2, \beta + 2}(m).
\endaligned
$$
It remains to rewrite the Hilbert transform term. When e.g. $\mu\in (2/3,1]$, we have $0< \beta' + 1/2 < 1$ and   
\be\label{eq:HpaScomm}
\pa_s H \varphi_{ss}  =  \frac{1}{|z_\al|} H (|z_\al|\pa_s^2\varphi_{s}) =  H \pa_s^2\varphi_{s} + \frac{1}{|z_\al|} [H, |z_\al|]\, \pa_s^2\varphi_{s}.
\ee
However, by assumption $|z_\al| \in H^2_{\beta '}(m)$, hence
$$
\pa_s \bigg(\frac{1}{|z_\al|} [H, |z_\al|] \pa_s^2\varphi_{s}\bigg) \in \Lcal_{2, \beta' + 1}(m)\subseteq \Lcal_{2, \beta + 2}(m)
$$
by an argument similar to that of Lemma \ref{lem:derFphi} in the Appendix. When $\mu\in (1/2, 2/3]$, we have $0 < (\beta' - 1) + 1/2 < 1$, hence the claim follows applying Lemma \ref{lem:trHilbert} twice. We omit the details. In particular, we have \eqref{en1}.


On the other hand, we claim 
\be\label{en2}
\frac{d}{dt} \big\|\Lambda^{1/2}\big(m^{\lambda}\pa_s^2\varphi_{s}\big)\big\|^2_{2, \beta ' - \lambda + 1/2} = -\int m^{2(\beta ' + 1/2)} \si \pa_s^3\theta \, \pa_s H \pa_s^2\varphi_{s}  ds +  \Rcal_2(t)
\ee
with the remainder $\Rcal_2(t)$ bounded by at most the exponential of some power of the energy. Indeed, taking the time derivative, we have
$$
\frac{1}{2}\frac{d}{dt} \big\|\Lambda^{1/2}\big(m^{\lambda}\pa_s^2\varphi_{s}\big)\big\|^2_{2, \beta ' - \lambda + 1/2} = \int_{-\pi}^\pi m^{2(\beta ' - \lambda + 1/2)}  \Lambda^{1/2}\big( m^{\lambda}  (\pa_s^2\varphi_{s})_t\big) \Lambda^{1/2}\big(m^\lambda \pa_s^2\varphi_{s}\big)d\al   
$$
and we first claim 
\be\label{eq:varphissst}
\aligned
(\pa_s^2\varphi_{s})_t &=- (\varphi_{ss})^2 -2\varphi_{s} \pa_s^2\varphi_{s} + \pa_s^2(- \si\theta_s + \theta_t^2) \\
&= -2\varphi_{s} \pa_s^2\varphi_{s} + 2\theta_t \pa_s H\varphi_{ss} - \si \pa_s^3\theta + H^1_{\beta ' + 1}(m).
\endaligned
\ee
Indeed, by assumption, the first term clearly satisfies $(\varphi_{ss})^2 \in H^1_{\beta '}(m)$. On the other hand, Lemma \ref{lem:param} together with \eqref{eq:HpaScomm} implies
$$
\frac{1}{2}\pa_s^2(\theta_t^2) = (\theta_{ts})^2 + \theta_t \pa_s^2\theta_{t} = \theta_t \pa_s H\varphi_{ss} + H^1_{\beta' + 1}(m).
$$

Finally, using that $\si \in H^3_{\beta + 1}(m)$ by Lemma \ref{lem:si}, we conclude 
\be\label{aux:dersith}
\pa_s^2(\si \theta_s)  - \si \pa_s^3\theta = 2 \si_s \theta_{ss} + \si_{ss}\theta_s \in H^1_{\beta' + 1}(m).
\ee
Some care is needed here, though. If we had asymptotic estimates $\theta_{ss} = O(m^{\mu - 2})$ or $\si_{ss} = O(m^{\mu - 1})$ (we do have  corresponding asymptotic estimates for lower order derivatives) we could prove the right-hand side of \eqref{aux:dersith} actually belongs to $H^1_{\beta - \mu + 2}(m)$. However, we only have
$$
\theta_{ss} \in \Lcal_{2, \beta + 1}(m), \quad \si_{ss} = O(m^{-(\beta + 1/2)}) \quad \quad \Rightarrow \quad \quad \si_{ss}\theta_{ss} \in  \Lcal_{2, \beta' + 1}(m)
$$
provided 
$$
\beta + 1/2 \leq 1 - \mu/2
$$
which is indeed satisfied for all $(1/2, 1]$ (cf. \eqref{eq:betaII}). In particular, \eqref{eq:varphissst} follows. 

We claim 
\be\label{eq:varphissst2}
\Lambda^{1/2}(m^\lambda (\pa_s^2\varphi_{s})_t) =  \Lambda^{1/2}(m^\lambda\si \pa_s^3\theta )  + \Lcal_{2, \beta' - \lambda + 1/2}(m).  
\ee
Indeed, the part of $(\pa_s^2\varphi_{s})_t$ which belongs to $H^1_{\beta' + 1}(m)$ satisfies the desired estimate by Lemma \ref{lem:halfDer}. On the other hand, by Lemma \ref{lem:comm_gb}, we have 
$$
\Lambda^{1/2}(m^\lambda \varphi_{s}\pa_s^2\varphi_{s}) = \varphi_{s}\Lambda^{1/2}(m^\lambda \pa_s^2\varphi_{s}) +  \big[ \Lambda^{1/2}, \varphi_s\big]\big(m^{\lambda}  \pa_s^2\varphi_{s}  \big) \in \Lcal_{2, \beta' - \lambda + 1/2}(m),
$$
since $(\beta' - \lambda)$ satisfies Muckenhaupt condition by construction. 

As for the remaining term, when $\mu\in (2/3, 1]$, by \eqref{eq:HpaScomm} and Lemma \ref{lem:halfDer} it is enough to show 
$$
\Lambda^{1/2}(m^\lambda \theta_t H\pa_s^2\varphi_{s}) = \theta_t \, \Lambda^{1/2}(m^\lambda H\pa^2_s\varphi_{s}) + \big[\Lambda^{1/2}, \theta_t\big](m^\lambda H\pa^2_s\varphi_{s}) \in \Lcal_{2, \beta' - \lambda + 1/2}(m).
$$
The commutator is bounded in $\Lcal_{2, \beta' - \lambda + 1/2}(m)$ by Lemma \ref{lem:comm_gb} (we have $\theta_{ts} = O(m^{\mu - 1})$), while Lemma \ref{lem:derFphi} applied with $\ga = \beta'$ implies
$$
\pa_s[m^\lambda, H] \pa_s^2\varphi_{s} \in \Lcal_{2, \beta' -\lambda + 1}(m),
$$
hence 
$$
\Lambda^{1/2}(m^\lambda H\pa_s^2\varphi_{s}) = H \Lambda^{1/2}(m^\lambda \pa_s^2\varphi_{s}) + \Lambda^{1/2}([m^\lambda, H] \pa_s^2\varphi_{s})\in \Lcal_{2, \beta' -\lambda + 1/2}(m)
$$
by Lemma \ref{lem:halfDer}. In particular, \eqref{eq:varphissst2} follows. 

On the other hand, when $\mu\in (1/2, 2/3]$, we have
$$
m^{\lambda}\pa_sH\varphi_{ss} = \pa_s H(m^{\lambda}\varphi_{ss}) + H^1_{\beta'-\lambda + 1}(m)
$$
by Lemmas \ref{lem:trHilbert} and \ref{lem:derFphi} (recall that $\lambda > 1$ in this case). Since $m^{\lambda}\varphi_{ss}\in H^1_{\beta'-\lambda}(m)$ with $\beta' - \lambda$ Muckenhaupt by construction we can proceed as in \eqref{eq:HpaScomm} to conclude
$$
m^{\lambda}\pa_sH\varphi_{ss} = H((m^{\lambda}\pa_{s}^2\varphi_{s})) + H^1_{\beta'-\lambda + 1}(m).
$$
In particular, we have
$$
\Lambda^{1/2}(m^\lambda \theta_t \pa_sH\varphi_{ss}) = \theta_t\Lambda^{1/2}(H(m^{\lambda}\pa^2_s\varphi_{s})) + \Lcal_{2, \beta' - \lambda + 1/2}(m)
$$
where, as in the case $\mu\in (2/3, 1]$, we have used Lemma \ref{lem:halfDer}.

It remains to consider the term containing $\si$. We have
$$
\aligned
\int_{-\pi}^\pi m^{2(\beta' - \lambda + 1/2)}  \Lambda^{1/2}\big( m^{\lambda} \si \pa_s^3\theta  \big) &\Lambda^{1/2}\big(m^\lambda \pa_s^2\varphi_{s}\big)d\al =\\
&= \int_{-\pi}^\pi m^{\lambda} \si \pa_s^3\theta \Lambda^{1/2}\big( m^{2(\beta' - \lambda + 1/2)}  \Lambda^{1/2}\big(m^\lambda \pa_s^2\varphi_{s}\big)\big)d\al,
\endaligned
$$
where we can write
$$
\aligned
\Lambda^{1/2}\big( m^{2(\beta' - \lambda + 1/2)}  \Lambda^{1/2}\big(m^\lambda \pa_s^2\varphi_{s}\big)\big) = m^{2(\beta' -\lambda+ 1/2)}&\Lambda \big(m^\lambda \pa_s^2\varphi_{s}\big) + \\
&+\big[\Lambda^{1/2},  m^{2(\beta' - \lambda + 1/2)}\big] \Lambda^{1/2}\big(m^\lambda \pa_s^2\varphi_{s}\big).
\endaligned
$$
Since $\si \pa_s^3\theta\in \Lcal_{2, \beta - \mu + 1}(m) \subseteq \Lcal_{2, \beta'}(m)$ and Lemma \ref{lem:comm_weight} implies
$$
m^{\lambda}\big[\Lambda^{1/2},  m^{2(\beta' - \lambda + 1/2)}\big] \Lambda^{1/2}\big(m^\lambda \pa_s^2\varphi_{s}\big)\in \Lcal_{2, -\beta'}(m),
$$
we are left with   
$$
 \int_{-\pi}^\pi  m^{2(\beta' + 1/2)}\si \pa_s^3\theta (m^{-\lambda} \Lambda\big(m^\lambda \pa_s^2\varphi_{s}\big)\big)d\al.
$$
Using $\Lambda = \pa_\al H$, we can write 
$$
\aligned
m^{-\lambda}\Lambda \big(m^\lambda \pa_s^2\varphi_{s}\big) &= m^{-\lambda}\pa_\al [H, m^\lambda] \pa_s^2\varphi_{s}  +  m^{-\lambda} \pa_\al (m^{\lambda} H\pa_s^2\varphi_{s})\\
&=\pa_\al H\pa_s^2\varphi_{s} + \lambda \frac{m_\al}{m} H\pa_s^2\varphi_{s}  + m^{-\lambda}\pa_\al [H, m^\lambda] \pa_s^2\varphi_{s}, 
\endaligned
$$
where by assumption and the discussion above we know
$$
m^{-\lambda}\pa_\al [H, m^\lambda] \pa_s^2\varphi_{s}\in \Lcal_{2, \beta' + 1}(m), \qquad \frac{1}{m} H\pa_s^2\varphi_{s} \in \Lcal_{2, \beta'+ 1}(m).
$$
In particular, we have shown \eqref{en2}.
\end{proof}

\subsection{Regularization of the evolution equations}
\label{ss.regularization}

In this section, we show existence of solutions to the regularized evolution equations. We first add a 'viscosity'-term to the $\om_t$ equation, which renders the system well-posed regardless of the sign of $\si$ (cf. \cite{CCG}). More precisely, we consider the following modification of the system  \eqref{eqs:sys}:
\be\label{eqs:sys-e}
\aligned
&(\log |z_\al|)_t = \varphi_s, \\
&\theta_t = H\varphi_s + \Rcal, \\
&\frac{\om_t}{2} + BR(z,\om)_t \cdot z_s = \epsilon \, m^2 \varphi_{ss}.
\endaligned
\ee
Here $\varphi_s$ is defined via \eqref{def:phi-om} and we have set   
\be\label{def:rcal}
\Rcal := \frac{\om \theta_s}{2}  + BR(z,\om)_s \cdot z_s^\perp - H\Big(\frac{\om_s}{2} + BR(z,\om)_s\cdot z_s\Big).
\ee
It will be convenient to append the (corresponding) evolution equation for $\varphi_s$, i.e.
\be\label{eq:varphist-e}
\varphi_{st} =- \varphi_s^2 + \theta_t^2 - \si \theta_s + \epsilon (m^2 \varphi_{ss})_s
\ee
to the system of equations \eqref{eqs:sys-e}. We call \eqref{eqs:sys-e} together with \eqref{eq:varphist-e} the $(\epsilon)$-system. We assume the solutions satisfy \eqref{def:ac}--\eqref{def:ac2}, and we require an additional half-derivative on $\varphi_s$ (and therefore on $\log|z_\al|$), i.e.
\be\label{eq:idata}
\theta\in H^{k+1}_{\beta + k}(m), \quad \log|z_\al| \in H^{k+1}_{\beta' + (k-1)}(m), \quad \varphi_{s} \in H^{k+1}_{\beta' + (k-1)}(m), \quad \om\in H^2_{\beta}(m)
\ee
(recall that $\beta'$ has been defined in \eqref{def:beta'}). 

\begin{lemma}\label{lem:apriori-e} Let $k\geq 2$ and let $(z,\om)$ be a sufficiently regular solution of the $(\epsilon)$-system for some $\epsilon>0$ fixed. Then, there exists $p\in \Nbb$ such that the energy functional 
$$
\t E_k(t)^2 := \|\log|z_\al|\|^2_{H^{k+1}_{\beta' + (k-1)}} +  \|\theta\|^2_{H^{k+1}_{\beta + k}} + \|\varphi_s\|^2_{H^{k+1}_{\beta' + (k-1)}} +  \|\om\|^2_{H^{2}_{\beta}}  + \|z\|^2_\Fcal
$$ 
satisfies the a priori energy estimate 
\be\label{eqs:tenergy}
\frac{d\t E_k(t)^p}{dt} \, \leq \,  C(\epsilon) \exp  (C \t E_k^p(t)),
\ee
where $C(\epsilon) = O(\epsilon^{-1})$. 
\end{lemma}

\begin{proof}
We show the claim for $k=2$. Generalizing the results of Section \ref{s:estDerivatives} to solutions of the system \eqref{eqs:sys-e} is straightforward. Compared to our original system we control an additional derivative of $\varphi_s$, hence $\theta_t$ belongs to $H^3_{\beta + 2}(m)$ with the norm controlled in terms of $\t E_2(t)$. 

Since $m^2\pa_s \varphi_s\in H^2_{\beta}(m)$, the estimate for $\om_t$ remains valid. The only interesting part is the time derivative of $\varphi_{s}$. We have
$$
\aligned
\frac{1}{2}\frac{d}{dt}\|\pa_s^3 \varphi_s\|^2_{2, \beta' + 1} =    \epsilon\int_{-\pi}^{\pi} m^{2(\beta' + 1)} &\pa_s^3 \varphi_s \pa^{4}_s(m^2 \varphi_{ss}) ds -\int_{-\pi}^{\pi} m^{2(\beta' + 1)} \pa_s^3 \varphi_s \pa^{3}_s(\si \theta_s) ds \\
&+  \text{`bounded terms'}.
\endaligned
$$  
We denote the integrals on the r.h.s. of the above equation by $I_1$ and $I_2$ respectively. For $I_1$, it is not difficult to see that integration by parts yields
$$
I_1 = - \epsilon \int_{-\pi}^{\pi} m^{2(\beta' + 2)} |\pa_s^4 \varphi_s|^2ds + \text{`bounded terms'}.
$$
For $I_2$ on the other hand, we have
\be\label{eq:I1}
\aligned
I_2 &= \int_{-\pi}^{\pi} \pa_s(m^{2(\beta' + 1)} \pa_s^3 \varphi_s) \pa^{2}_s(\si \theta_s) ds \\
&=\int_{-\pi}^{\pi} m^{2(\beta' + 1)} \pa_s^4 \varphi_s \pa^{2}_s(\si \theta_s) ds + 2(\beta' + 1)\int_{-\pi}^{\pi} m^{2(\beta' + 1)}m' \pa_s^3 \varphi_s \pa^{2}_s(\si \theta_s) ds.
\endaligned
\ee
The second integral is clearly bounded in terms of $\t E_2(t)$, while for the first one we have the estimate 
\be\label{eq:1/eps}
\int_{-\pi}^{\pi} m^{2(\beta' + 1)} \pa_s^4 \varphi_s \pa^{2}_s(\si \theta_s) ds \, \leq \, \frac{\epsilon}{2} \int_{-\pi}^{\pi} m^{2(\beta' + 2)}| \pa_s^4 \varphi_s |^2 ds + \frac{1}{2\epsilon} \int_{-\pi}^{\pi} m^{2\beta'} |\pa^{2}_s(\si \theta_s)|^2 ds.
\ee
This concludes the proof, since the most singular terms from $I_1$ and $I_2$ cancel out. 
\end{proof}

In order to prove the local existence of solutions to the $(\epsilon)$-system when $\epsilon>0$, we introduce an additional $\de$-dependent regularization such that the resulting $(\epsilon, \de)$-system can be written as an ODE on an open set of a suitable Banach space. This is accomplished applying a variable-step convolution operator $A^\de$ to the highest order derivative terms. Then, we can use the abstract Picard theorem to find a sequence of solutions  whose subsequence converges to a solution of the $(\epsilon)$-system by compactness. The argument is standard, we therefore omit most of the details and only indicate the differences due to weights (the interpolation inequalities are to be replaced by their weighted counterparts, cf. \cite{CSLin}).

We first specify the convolution operator $A^\de$. Let $\phi$ be a positive, symmetric mollifier, i.e. a smooth function $\phi:\Rbb \rightarrow \Rbb$ such that 
$$
\phi(\al) \geq 0, \qquad  \phi(\al) = \phi(-\al), \qquad \text{supp} \, \phi \subseteq B_1(0), \qquad  \int_\Rbb \phi = 1,
$$
and let $\eta\in\Ccont^\infty(\Tbb\setminus\{\pm\al_*\})$ respect both symmetries and be strictly positive on $\Tbb\setminus\{\pm\al_*\}$ with first order zeros at~$\pm\al_*$, i.e. $\eta(\pm\al_*) = 0$ but $\eta'(\pm\al_*) \neq 0$. This condition ensures that $\eta \sim m$. Then, for a sufficiently small $\de >0$ and $\al\neq \pm\al_*$, we define 
\be\label{E.defphide}
\phi_{\de \eta(\al)}(\al') := \frac{1}{\de \eta(\al)}\phi\left(\frac{\al'}{\de \eta(\al)}\right).
\ee
We extend $\phi$ and $\eta$ periodically to $\Rbb$ and define
\be\label{E.defBde}
B_{\de}f(\al) := (\phi_{\de \eta(\al)} \ast f) (\al) := \int_\Rbb \phi_{\de \eta(\al)}(\al - \al')f(\al')d\al'.
\ee
The adjoint is 
\be\label{E.defB*de}
B_{\de }^*f(\al') :=  \int_\Rbb \phi_{\de \eta(\al)}(\al - \al') f(\al')d\al',
\ee
and the convolution operator $A^\de$ is defined as the composition of $B_\de$ and $B_\de^*$:
\be\label{E.defAde}
A^\de f := B^{*}_\de  B_\de f.
\ee
It is not difficult to see the restriction of $A^\de$ to $[-\pi, \pi]$ respects both symmetries and $A^\de(f)(\pi) = A^\de(f)(-\pi)$. 

As expected these operators have the smoothing effect away from $\pm \al_*$, however they also respect growth rates as we approach $\pm \al_*$. Moreover, for fixed $\al \neq \pm\al_*$, the corresponding interval of integration always has positive distance to $\pm\al_*$, which makes them well-adapted for regularization of functions that live in weighted Sobolev spaces $\Lcal^k_{2,\ga}(\Tbb)$ (in particular those having non-integrable singularities). Note that taking a derivative on $\phi$ results in a factor of size $O((\de \eta)^{-1})$, which is why we need $\eta \sim m$. Precise technical results are presented in the Appendix \ref{a.convolution}. 

When defining the $(\epsilon, \delta)$-system, it will be convenient to introduce an additional variable $\Phi_s$ which basically satisfies the same evolution equation as $\varphi_s$, the only difference being the convolution operator applied to the $\si\theta_s$ term. More precisely, for $\epsilon>0$, we consider   
\be\label{eqs:sys-ed}
\aligned
&(\log |z_\al|)_t = \Phi_s, \\
&\theta_t = H\Phi_s + \Rcal, \\
&\frac{\om_t}{2} + BR(z,\om)_t \cdot z_s = \epsilon \, m^2 \,\Phi_{ss}, \\
&\varphi_{st} = -\Phi_s \varphi_s + \theta_t^2 -\si \theta_s + \epsilon (m^2 A^\de (\Phi_{ss}))_s, \\
&\Phi_{st} =- \Phi_s^2 + \theta_t^2 - A^\de(\si \theta_s) + \epsilon (m^2 A^\de (\Phi_{ss}))_s,
\endaligned
\ee
where $\Rcal$ is defined by \eqref{def:rcal} and $\varphi_s$ via \eqref{def:phi-om}. For the initial data, we take 
\be\label{eq:idata-ed}
\theta(\cdot, 0) = \theta^0, \qquad \log|z_\al(\cdot, 0)| = \log|z_\al^0|, \qquad \Phi_s(\cdot, 0) = \varphi_s(\cdot, 0) = \varphi_s^0, \qquad \om(\cdot, 0) = \om^0.
\ee
When $\de = 0$, these initial data ensure $\Phi_s(\cdot, t) = \varphi_s(\cdot, t)$ and we recover solutions of the $(\epsilon)$-system.

Let $\wt{\mathcal{B}}^k_{\beta, \mu}(m)$ for $k\geq 2$, denote the Banach space of all $(\theta, \log|z_\al|, \varphi_s, \Phi_s, \om)$ satisfying the appropriate symmetry assumptions such that  
$$
(\theta, \log|z_\al|) \in H^{k+1}_{\beta + k}(m) \times H^{k+1}_{\beta' + (k-1)}(m),  
$$
respectively 
$$
(\varphi_s, \Phi_s) \in H^{k}_{\beta'+ (k-2)}(m) \times H^{k+1}_{\beta' + (k-1)}(m), \qquad \om\in H^2_{\beta}(m).
$$
Let $\t \Ocal^k_{\beta, \mu}$ be the open set of all elements of $\wt{\mathcal{B}}^k_{\beta, \mu}(m)$ which satisfy \eqref{def:ac}--\eqref{def:ac2}.

The results of Section \ref{s:estDerivatives} readily generalize to the present case. The right-hand side of \eqref{eqs:sys-ed} has values in $\wt{\mathcal{B}}^k_{\beta, \mu}(m)$ if the conditions \eqref{def:ac} on the arc-chord are satisfied. The existence of $\om\in H^{k+1}_{\beta + (k-1)}(m)$ follows from Lemma \ref{lem:omsss}. We can then show $\Rcal \in H^{k+1}_{\beta + k}(m)$ and the same is true for $\theta_t$ given $\Phi_s \in H^{k+1}_{\beta' + (k-1)}(m)$. Moreover, it is not difficult to see that higher-order asymptotic estimates (of type \eqref{def:ac2}) must hold as well. To control $k$ derivatives of $\om_t$ we need control over $k+1$ derivative of $\Phi_s$. In order 
to show $\varphi_{st} \in H^{k}_{\beta'+ (k-2)}(m)$, it is enough to control $\theta\in H^{k+1}_{\beta + k}(m)$. In particular, there is no need for a convolution operator on the corresponding term. Finally, the r.h.s. of \eqref{eqs:sys-ed} is actually Lipschitz on $\t\Ocal^k_{\beta, \mu}(m)$. The estimates would be similar to those handling the difference $z_t(\al) - z_t(\al')$. We omit the details. 

In particular, by the abstract Picard theorem, there exists $T^{\epsilon, \de}>0$ and $(z^{\epsilon, \de}, \varphi_s^{\epsilon, \de}, \Phi_s^{\epsilon, \de}, \om^{\epsilon, \de})\in \Ccal^1([0, T^{\epsilon, \de}], \, \t\Ocal^k_{\beta, \mu})$ solutions of the corresponding initial value problem \eqref{eqs:sys-ed}-\eqref{eq:idata-ed}. These can be extended until the solution leaves the open set $\t \Ocal^k_{\beta, \mu}$. 

\begin{lemma}\label{lem:Tepsilon}
Let $k\geq 2$ and let $\epsilon>0$ fixed. Then, there exists $T^\epsilon>0$ independent of~$\de$ such that $(z^{\epsilon, \de}, \varphi_s^{\epsilon, \de}, \Phi_s^{\epsilon, \de}, \om^{\epsilon, \de})$ exist on $[0,T^\epsilon]$ for all small enough $\de> 0$.
\end{lemma}

\begin{proof}
We claim the time derivative of the extended energy functional 
$$
\t E_k(t)^2 := \|\log|z_\al|\|^2_{H^{k+1}_{\beta' + (k-1)}} +  \|\theta\|^2_{H^{k+1}_{\beta + k}} + \|\Phi_s\|^2_{H^{k+1}_{\beta' + (k-1)}} +  \|\varphi_s\|^2_{H^{k}_{\beta' + (k-2)}} + \|\om\|^2_{H^2_{\beta}} + \|z\|^2_\Fcal
$$
satisfies estimate \eqref{eqs:tenergy} uniformly for all $\de> 0$ (for simplicity, we keep the same notation as in Lemma \ref{lem:apriori-e}). By integration, we then obtain an upper bound on $T^\epsilon>0$ in terms of $\t E_k(0)$.   

Without loss of generality, take $k=2$. It is enough to consider terms with $\pa_s^2\varphi_s$ and $\pa_s^3\Phi_s$ (cf. the proof of Lemma \ref{lem:apriori-e}). We have
$$
\frac{1}{2}\frac{d}{dt}\|\pa_s^2 \varphi_s\|^2_{2, \beta'} = \epsilon\int_{-\pi}^\pi m^{2\beta'} \pa_s^2 \varphi_s \pa^{3}_s(m^2 A^\de(\Phi_{ss})) ds + \text{ `bounded terms'},
$$
where here, and in the sequel we write 'bounded terms' for all terms which are uniformly bounded in $\de$. We denote this integral by $I$. Since we control two derivatives of $\si \theta_s$, there is no need to single out the corresponding integral, cf. \eqref{eq:I1}. We concentrate on the most singular contribution to~$I$. Using the results of Lemma \ref{lem:B_delta} it is not difficult to see that 
$$
\aligned
I &= \epsilon\int_{-\pi}^\pi m^{2(\beta' + 1)} \pa_s^2 \varphi_s \pa^{3}_s A^\de(\Phi_{ss}) ds + \text{`bounded terms'} \\
&= \epsilon\int_{-\pi}^\pi  B_\de \Big(m^{2(\beta' + 1)} \frac{1}{|z_\al|^3}\pa_\al\varphi_{ss}\Big) \pa_\al B_\de(\pa^2_\al\Phi_{ss}) d\al + \text{`bounded terms'}. 
\endaligned
$$
In particular, we have
\be\label{eq:est}
|I| \, \leq \,  \frac{\epsilon}{N}\int_{-\pi}^\pi m^{2(\beta' + 2)}\frac{1}{|z_\al|^4}|\pa_\al B_\de(\pa^2_\al\Phi_{ss})|^2d\al + \text{`bounded terms'},
\ee
where $N\in \Nbb$ is a sufficiently high number to be determined later. 

When estimating the time derivative of $\|\pa_s^3\Phi_s\|^2_{2,\beta' + 1}$ there are two terms $I_1$ and $I_2$ that we need to consider, cf. proof of Lemma \ref{lem:apriori-e}. We have   
$$
\aligned
I_1 &= \epsilon\int_{-\pi}^\pi m^{2(\beta' + 1)} \pa_s^3 \Phi_s \pa^{4}_s\left(m^2 A^\de(\Phi_{ss})\right) ds\\
&= \epsilon\int_{-\pi}^\pi m^{2(\beta' + 2)} \pa_s^3 \Phi_s \pa^{3}_s\left(  \frac{1}{|z_\al|}B_\de^* \pa_\al B_\de (\Phi_{ss})\right) ds + \epsilon\int_{-\pi}^\pi m^{2(\beta' + 2)} \pa_s^3 \Phi_s \pa^{3}_s\t K_\de B_\de (\Phi_{ss})d\al + \text{`ok'} \\
&= I_{11} + I_{12} + \text{`bounded'},
\endaligned
$$
where we have used $(B_\de^* f)' = B_\de^* f' + \t K_\de f$ and Lemma \ref{lem:B_delta}. We consider $I_{12}$ first. Repeatedly using Lemma \ref{lem:B_delta} we find
$$
\aligned
I_{12} &= \epsilon\int_{-\pi}^\pi m^{2(\beta' + 2)} \frac{1}{|z_s|^3} \pa_s^3 \Phi_{s} \t K_\de( \pa_\al B_\de (\pa_\al^2\Phi_{ss})) d\al + \text{`bounded terms'} \\
&= \epsilon\int_{-\pi}^\pi \t K_\de^*\left(m^{2(\beta' + 2)} \frac{1}{|z_s|^3} \pa_s^3 \Phi_{s}\right)  \pa_\al B_\de (\pa_\al^2\Phi_{ss}) d\al + \text{`bounded terms'},
\endaligned
$$
where $\t K^*_\de$ is the adjoint of $\t K_\de$. In particular, absolute value of $I_{12}$ satisfies estimate \eqref{eq:est}. As for $I_{11}$, we have
$$
\aligned
I_{11} &= \epsilon\int_{-\pi}^\pi m^{2(\beta' + 2)} \pa_s^3 \Phi_s \pa^{3}_s \left(\frac{1}{|z_\al|} B_\de^* \pa_\al B_\de (\Phi_{ss})\right)ds\\
&= -\epsilon\int_{-\pi}^\pi  \pa_\al B_\de \Big(m^{2(\beta' + 2)} \frac{1}{|z_\al|^4}\pa^2_\al\Phi_{ss}\Big)  \pa_\al B_\de(\pa^2_\al\Phi_{ss}) d\al + J +  \text{`bounded terms'},
\endaligned
$$ 
where $J$ is the sum of different terms, all of the same order as $I_{12}$, which arise as errors  when we interchange a derivative with $B_\de$ or $B_\de^*$. In particular, each of these terms satisfies estimate \eqref{eq:est}. In the remainder of the proof, we always group all such terms under the name $J$, however their number might change from line to line. Using that we control an additional derivative whenever we have a commutator with any of these convolution operators, e.g. 
$$
\pa_\al \big[B_\de, \, m^{2(\beta' + 2)}/|z_\al|^4\big](\pa^2_\al \Phi_{ss}) \in \Lcal_{2, -(\beta' + 2)}(m),
$$
by Lemma \ref{lem:B_delta}, we conclude
$$
I_1 = -\epsilon\int_{-\pi}^\pi  m^{2(\beta' + 2)} \frac{1}{|z_\al|^4} |\pa_\al B_\de(\pa^2_\al\Phi_{ss})|^2 d\al + J + \text{`bounded terms'}. 
$$
Finally, the remaining term reads
$$
\aligned
I_2 =& -\int_{-\pi}^\pi m^{2(\beta' + 1)} \pa_s^3 \Phi_s \pa^{3}_s A^\de(\si\theta_s) ds  \\
=& \int_{-\pi}^\pi \pa_\al B_\de\Big(m^{2(\beta' + 1)}\frac{1}{|z_\al|^4} \pa_\al^2 \Phi_{ss} \Big) \pa^{2}_\al B_\de(\si\theta_s)  d\al + \text{`bounded terms'}.
\endaligned
$$
Again, we have
$$
\pa_\al \big[B_\de, \, m^{2(\beta' + 1)}/|z_\al|^2\big](\pa^2_\al\Phi_{ss}) \in \Lcal_{2, -\beta'}(m)
$$
and we can conclude the argument as in Lemma \ref{lem:apriori-e}. 

Piecing everything together we find 
$$
I + I_1 + I_2 \, \leq \, -\epsilon\int_{-\pi}^\pi  m^{2(\beta' + 2)} \frac{1}{|z_\al|^4} |\pa_\al B_\de(\pa^2_\al\Phi_{ss})|^2 d\al + J + \text{`bounded terms'}
$$ 
where 
$$
|J| \, \leq \, \frac{n\epsilon}{N}\int_{-\pi}^\pi  m^{2(\beta' + 2)} \frac{1}{|z_\al|^4} |\pa_\al B_\de(\pa^2_\al\Phi_{ss})|^2 d\al
$$
for some $n\in \Nbb$. In particular, choosing $N \geq n$, we conclude that $I + I_1 + I_2$ are bounded in terms of $\t E_2(t)$ and estimate \eqref{eqs:tenergy} follows.
\end{proof}


\subsection{Existence of solutions to the original system}
\label{ss.existence}

In the previous section, we have shown how to construct $(z^{\epsilon}, \om^{\epsilon})$ 
a solution of the $(\epsilon)$-system to the initial data \eqref{eq:idata}. It remains to show these solutions actually exist on some common time interval for all $\epsilon\geq 0$. In order to do so, we consider the natural generalization $E_{k,\epsilon}(t)$ of the energy functional $E_k(t)$ to the $(\epsilon)$-system which involves     
$$
E^k_{h,\epsilon}(t)^2 := \, \epsilon \, \left(\|\pa_s^{k+1}\varphi_s\|^2_{2,\beta' +(k-1)} + \|\pa_s^{k+1}\log|z_\al|\|^2_{2,\beta' +(k-1)}\right)\, + E^k_h(t)^2, \qquad k\geq 2,
$$
with $E_h^k(t)$ as in \eqref{def:energy-high}. 

Let us specialize to the case $k=2$. To prove the estimate $\si_t = O(m^{\mu + 1})$, we require $\om_{tt}\in H^2_{\beta}(m)$ (cf. Lemma \ref{lem:si_t}), which in turn depends on $\epsilon^2 \pa_s^5\varphi_s$ and we only control $\pa_s^3\varphi_s$. We therefore mollify the original initial data, solve the corresponding initial value problem for a sufficiently high $k$, then prove an a-priori estimate for a higher order energy functional $E_{2,\epsilon}(t)$ which satisfies $E_{2,\epsilon}(0) \, \lesssim \, E_2(0)$ uniformly in $\epsilon$, with $E_{2,\epsilon}(0)$ corresponding to the mollified initial data and $E_2(0)$ to the original initial data. We set 
$$
E_{2, \epsilon}(t)^2 = E_l(t)^2 + \sum_{i=2}^{ 4} (\sqrt{\epsilon})^{(i-2)} E_{h, \epsilon}^i(t)^2.
$$
It is not difficult to see this scales correctly if we take $\epsilon^{1/4}$ as the mollification parameter, i.e. we mollify using $B_{\epsilon^{1/4}}$, see Appendix \ref{a.convolution} and \eqref{E.defBde}. Finally, as noted in the previous paragraph, we also need to mollify the initial data, part of which does not belong to the family  $\Lcal^k_{2,\ga}(m)$ we used in Lemma \ref{lem:B_delta}, Appendix \ref{a.convolution}. However, it is not difficult to adapt its proof to see that $B_{\epsilon^{1/4}}$ is bounded on $H^k_\ga(m)$ as well. We are now ready to prove:

\begin{lemma}\label{lem:existence} Let $k = 2$ and let $(z^0, \om^0)\in \Bcal^k_{\beta, \mu}(m)$ be the initial data to the system of equations \eqref{eqs:sys} as in Theorem \ref{thm:T.main}.  
Let furthermore $(z^\epsilon, \om^\epsilon)$ be the solution of the $(\epsilon)$-system \eqref{eqs:sys-e} with $k=4$ and mollified initial data $(z^{0,\epsilon}, \om^{0,\epsilon})$, which written out read
$$
(B_{\epsilon^{1/4}}(\theta^{0}_s), \, B_{\epsilon^{1/4}}(\log|z^{0}_\al|), \,   B_{\epsilon^{1/4}}(\om^{0})).
$$ 
Then, there exists a time $T>0$, independent of $\epsilon$, such that $(z^{\epsilon}, \om^{\epsilon})$ exist on $[0, T]$ for all small $\ep>0$.
\end{lemma}

\begin{proof}[Proof of Lemma \ref{lem:existence}] We need to show the inequality \eqref{eq:energy} from Lemma \ref{lem:apriori} holds for $E_{2,\epsilon}(t)$. Note that $E_{2,\epsilon}(0) \lesssim E_2(0)$ as required. Consider the time derivative of $E^2_{h, \epsilon}(t)^2$. We claim it can be estimated in terms of $E^2_{h, \epsilon}(t)^2 + E_l(t)^2$, except for the term requiring $\si_t = O(m^{\mu + 1})$, for which we need higher order $E^i_{h, \epsilon}(t)$ as noted above. Indeed, the time derivative of $\epsilon\|\pa^3_s \varphi_s\|^2_{2, \beta' + 1}$ can be estimated as in the proof of Lemma \ref{lem:apriori-e}. The extra $\epsilon$ removes the $1/\epsilon$ from the r.h.s. of \eqref{eq:1/eps}. On the other hand, to estimate the term with  $\sqrt{\si}\,\pa^3_s\theta$ we proceed as in Lemma \ref{lem:apriori}, where the most singular term will be canceled by the corresponding term from $\Lambda^{1/2}(\pa_s^2\varphi_s)$. 
	
It remains to consider  
$$
\frac{d}{dt} \|\Lambda^{1/2}(m^{\lambda}\pa^2_s\varphi_s)\|^2_{2, \beta '-\lambda + 1/2}.
$$    
Let us consider the `viscosity' term from $\pa_s^2\varphi_{st}$ first. We have
\be\label{est1}
m^{\lambda}\pa_s^3(m^2 \varphi_{ss}) = \pa_\al (g_2 \pa_\al(m^\lambda \pa_s^2 \varphi_s)) + g_1\pa_\al(m^\lambda \pa_s^2\varphi_s) + \Lcal^1_{2,(\beta'-\lambda) + 1}(m),
\ee
where $g_i = O(m^{i})$, $g_i' = O(m^{i-1})$ for $i = 1,2$. {\color{ForestGreen}}Using that 
$$
\pa_s(m^\lambda \pa_s^2\varphi_s) \in \Lcal_{2,(\beta'-\lambda) + 1}(m),
$$
it is not difficult to adapt the proof of Lemma \ref{lem:comm_gb} to conclude
$$
\Lambda^{1/2}\big(g_1\pa_\al (m^\lambda \pa_s^2\varphi_s)\big) = g_1\pa_\al \Lambda^{1/2}(m^\lambda \pa_s^2\varphi_s) \in \Lcal_{2,(\beta'-\lambda) + 1/2}(m),
$$
hence an integration by parts shows the corresponding term satifies the required estimate. Note that we use the control over $\pa_s^3\varphi_s$ in order to estimate derivatives of the lower order terms in \eqref{est1}. 

In particular, it remains to consider the most singular term
$$
\aligned
I_1 = & \,\epsilon\int_{-\pi}^{\pi} m^{2(\beta' - \lambda) + 1}\Lambda^{1/2}(m^{\lambda}\pa^2_s\varphi_s)\Lambda^{1/2}(\pa_\al (g_2 \pa_\al(m^\lambda \pa_s^2 \varphi_s)))d\al \\
=& \,\epsilon\int_{-\pi}^{\pi} m^{2(\beta' - \lambda) + 1}\Lambda^{1/2}(m^{\lambda}\pa^2_s\varphi_s)\pa_\al\Big( [\Lambda^{1/2}, g_2]\pa_\al(m^\lambda \pa_s^2 \varphi_s) + g_2\pa_\al\Lambda^{1/2}(m^\lambda \pa_s^2 \varphi_s)\Big)d\al.
\endaligned
$$
Note that 
$$
\pa_\al\Big([\Lambda^{1/2}, g_2]\pa_\al(m^\lambda \pa_s^2 \varphi_s)\Big) = -\frac{1}{2}\pa_\al g_2 \, \pa_\al\Lambda^{1/2}(m^\lambda \pa_s^2 \varphi_s) +  \Lcal_{2,(\beta'-\lambda) + 1/2}(m),
$$
where we have used $g''_2 = O(1)$. In particular, we have 
$$
\aligned
I_1 =& \,\epsilon\int_{-\pi}^{\pi} m^{2(\beta' - \lambda) + 1}\Lambda^{1/2}(m^{\lambda}\pa^2_s\varphi_s)\Big( \frac{1}{2}\pa_\al g_2\pa_\al\Lambda^{1/2}(m^\lambda \pa_s^2 \varphi_s) + g_2\pa_\al^2\Lambda^{1/2}(m^\lambda \pa_s^2 \varphi_s) \Big)d\al + \text{'ok'}\\
=&-\epsilon \int_{-\pi}^{\pi} m^{2(\beta' - \lambda) + 1} g_2|\pa_\al\Lambda^{1/2}(m^\lambda \pa_s^2 \varphi_s)|^2d\al + \text{`bounded terms'}
\endaligned
$$
where we have repeatedly integrated by parts. Since we actually have $g_2 = m^2/|z_\al|^2$, the claim follows.

Compared to the proof of Lemma \ref{lem:apriori}, there is one more interesting term. There, we used that $\pa_s^2(\si \theta_s) = \si \pa_s^3\theta + \Lcal_{2,\beta + 1}^1(m)$, which relies on the control of $\pa_s^3\si$. However, for the $(\epsilon)$-system \eqref{eqs:sys-e}, the equivalent of Lemma \ref{lem:si} implies
$$
\si_{ss} = \epsilon \, H(\pa_s^2(m^2 \varphi_{ss})) + \Lcal^1_{2, \beta + 1}(m).
$$
We can further rewrite this as in \eqref{est1}, which amounts to estimating   
$$
\aligned
I_2 &= -\epsilon\int_{-\pi}^{\pi} m^{2(\beta' - \lambda) + 1}\Lambda^{1/2}(m^{\lambda}\pa^2_s\varphi_s)\Lambda^{1/2}\big(m^{\lambda}\theta_sH(m^{2-\lambda} \pa_\al(m^\lambda \pa_s^2\varphi_{s}))\big)d\al \\
& = -\epsilon\int_{-\pi}^{\pi} m^{2(\beta' - \lambda) + 1}\theta_s H\pa_\al(m^{\lambda}\pa^2_s\varphi_s) m^\lambda H (m^{2-\lambda} \pa_\al(m^\lambda \pa_s^2\varphi_{s}))d\al + \text{`bounded terms'}. 
\endaligned
$$ 
An appropriate commutator of the type considered in Lemma \ref{lem:derFphi} implies everything is controlled by $ \epsilon\|\pa^3_s \varphi_s\|^2_{2, \beta' + 1}$. Higher order $E^i_{h, \epsilon}$ follow analogously. We omit further details.
\end{proof}



\section{The local existence theorem domain with corner}
\label{s.Corner}

In this section, we indicate how to modify the argument in Sections \ref{s.estimatesSingInt} and \ref{s.estEnergy} in order to prove the local existence theorem when outward cusps are replaced by corners of opening $2\nu>0$. As it turns out, the estimates simplify considerably.

As in Section \ref{s.estEnergy}, we consider the system of evolution equations \eqref{eqs:sys} for $(\log|z_\al|, \theta, \om)$ under the same symmetry assumptions. The tangent angle $\theta$ still satisfies \eqref{def:sym-xy-theta}, that is
$$
\theta(\al, t) = \pi - \theta(-\al, t), \qquad \theta(\al_* - \al, t) = - \theta(\al_*+ \al, t),
$$
however it is not continuous at $\pm\al_*$; at $\al = \al_*$ it has a jump of size $2\nu(t)\in (0, \frac{\pi}{2})$, i.e. 
$$
\lim_{\al \nearrow \al_*}\theta(\al, t) = -\nu(t), \quad \quad \lim_{\al \searrow \al_*}\theta(\al, t) = \nu(t)
$$  
and the size of the jump at $\al = -\al_*$ is determined by symmetry. The length $|z_\al|$ of the tangent vector is well-behaved throughout (it is even w.r.t. both axis and therefore continuous). Again, we recover the parametrization of the interface from $(\theta, \log|z_\al|)$ via integration with fixed integration constant 
$$
z(\al_*,t) = 0
$$
i.e. via \eqref{def:z-int-constant}--\eqref{def:z}.

Let $0 < \beta + \frac{1}{2} < 1$. 
Under current symmetry assumptions, we say that $(z(\cdot, t), \om(\cdot, t))$ belong to the Banach space $\Bcal^k_{\beta}(m)$ with $k\geq 2$, if
\be\label{regularity-z-corner}
\theta(\cdot, t)\in H^{k+1}_{\beta + k}(m), \qquad \log|z_\al(\cdot, t)|\in H^k_{\beta + (k-1)}(m),
\ee
respectively
\be\label{regularity-om-corner}
\om(\cdot, t)\in \Lcal^{k + 3/2}_{\beta + k - 1/2}(m)
\ee
with the weight function defined as in Section \ref{s.estEnergy}. As the interface is only piecewise smooth, recall that e.g. $H^1_\beta(m)$ is embedded in the space of piecewise continuous functions only, cf. Section \ref{ss.weightedSobolev}. For instance, by Hardy inequalities, we have 
$$
\theta(\al, t) \mp \nu(t) = \int_{\al_*}^{\al} \theta_s(\al', t)ds_{\al'} \in \Lcal^{k+1}_{2,  \beta + k}(m), \qquad \al \gtrless \al_*.
$$  
Note also how we assume $\om_s(\al_*, t) = 0$ in addition to $\om(\al_*, t) = 0$. As we will see in Lemma \ref{lem:corner-ipp} further down, this assumption is important when taking derivatives of the Birkhoff-Rott integral in this setting. 

	The parametrization is assumed to satisfy the following version of the arc-chord condition:
	\be\label{def:ac-corner}
	\|z\|_\Fcal := \Fcal_r(z) + \sup_{\al\in B_r(\al_*)} \frac{1}{|z_{1\al}(\al, t)|} + \sup_{\al\in B_r(\al_*)} \frac{1}{|z_{2\al}(\al, t)|} < \infty
	\ee
	for some small $r>0$, with $\Fcal_r(z)$ as in \eqref{arcchord}. We also assume that the normal component of the pressure gradient $\si = \pa_n P$ satisfies the following version of the Rayleigh-Taylor condition:
	\be\label{RTcondition-corner}
	\inf m(\cdot)^{-1}\si(\cdot, t) > 0.
	\ee    
	\begin{remark}\label{R.opening}
		As we will see in Lemma \ref{lem:corner-preliminary} below, under current assumptions the opening angle $\nu$ must change with time and the Rayleigh-Taylor condition \eqref{RTcondition-corner} is satisfied provided  certain continuous linear functional $b$ does not vanish. In fact, both the value of $\theta_t$ at the singular points and the lower bound in \eqref{RTcondition-corner} are directly proportional to $b$. More precisely, we have
		\be\label{eq:RT-b}
		b(t) : = \frac{1}{2\pi i} p.v.\int_{-\pi}^\pi \om(\al',t)\frac{1}{z(\al',t)^2}\,ds_{\al'}\neq 0
		\ee
		(i.e. $b \equiv b_1(\om)$ is the second moment of $\om$ in the sense of definition \eqref{def:bi}), where by symmetry we have $\Im b(t) = 0$. It is not difficult to choose a set of initial $(\om, z)$ such that $\Re b(0) \neq 0$. In fact, using  symmetry and an appropriate cut-off it is enough to specify the upper part of $\Ga$ in the vicinity of the singular point, i.e. $\Ga^+$ (as defined in \eqref{def:Ga_pm}) and then prescribe $\om$ on $\Ga^+\cap \{0 < x < 2\de\}$. Then taking e.g. $\om(x,0) = x^2$ and $z(x, 0) = x + iax$ with constant $a>0$, it is not difficult to see the corresponding integral over say $(0, \de)$ is strictly positive.    
	\end{remark}
		
	\begin{theorem}\label{thm:T.main-corner} Let $k\geq 2$ and let $(z^0, \om^0)\in \Ocal^k_{\beta}(m)$, where $\Ocal^k_{\beta}(m)\subseteq \Bcal^k_{\beta}(m)$ is the open set defined via \eqref{eq:RT-b}, the arc-chord condition \eqref{def:ac-corner} and 
	\be\label{eq:conditions-corner}
	2\nu(t)\in\left(0,\frac{\pi}{2}\right), \qquad \max\left\{3 - \frac{\pi}{2\nu(t)}, 3 - 2\frac{\pi}{2\overline{\nu}(t)}\right\} < \beta + \frac{1}{2}, \qquad \beta + \frac{1}{2} \neq 2-\frac{\pi}{2\overline{\nu}(t)}, 
	\ee
	where we have set $2\overline{\nu}(t) = \pi - 2\nu(t)$. Then, there exists a time $T>0$ and $(z, \om)\in\Ccal([0, T], \Ocal^k_{\beta}(m))\cap \Ccal^1([0, T], \Ocal^{k-1}_{\beta}(m))$ 
	solution of \eqref{eqs:sys} up to time $T$ such that $z(\cdot, 0) = z^0$ and $\om(\cdot, 0) = \om^0$.
	\end{theorem}

We first give some general properties of the Birkhoff-Rott integral when cusps ($\nu = 0$) are replaced by corners $(\nu > 0)$. These correspond to results of Section \ref{s.BR}.

\begin{lemma}\label{lem:BRbasic-corner}
	Let $0< \beta + 1/2 <1$. Then, we have 
	$$
	BR(z, \cdot )^*: \Lcal_{2,\beta}(m) \rightarrow \Lcal_{2,\beta}(m)
	$$ 
	and the same is true for $BR(z, \cdot )^*z_s$. In particular, the conclusions of Corollary \ref{lem:BRcorrections} remain true. In addition, we have  
	\be\label{eq:corner-der-on-ker}
	f\in\Lcal_{2,  \beta}(m) \quad \Rightarrow \quad BR(z,f)^*iz_s - Hf\in\Lcal^{k-1}_{2,\beta + (k-1)}(m),
	\ee
	where $k\geq 2$ is as in the regularity assumptions \eqref{regularity-z-corner} for the parametrization of the interface. In other words, we have
        $$
	BR(z, f)\cdot z_s \in \Lcal^{k-1}_{2,\beta + (k - 1)}(m), \qquad BR(z, f)\cdot z_s^\perp - Hf \in \Lcal^{k-1}_{2,\beta + (k - 1)}(m).
	$$
	and conclusions of Lemmas \ref{lem:BRcancel}--\ref{lem:struct1} remain true. 
\end{lemma}

\begin{proof}
As before, we pass to the graph parametrization $x\pm i\ka(x,t)$ in the neighborhood of the origin, where now
$$
\pa_x\ka(x, t) = \tan\theta(x,t), \qquad \pa_x\ka(x,t) = \pm\tan\nu(t) + O(|x|^\lambda), \qquad x\gtrless 0,
$$
with $\lambda = 1 - (\beta + 1/2)$. By integration, we then have 
$$
\ka(x,t) =  \tan\nu(t)|x| + O(|x|^{1 + \lambda}), \qquad |x|< 2\de.
$$
In particular, in a sufficiently small neighborhood of the origin, we have
\be\label{eq:ka-corner}
\pa_x\ka(x,t) \sim \text{sign}(x)\tan\nu(t), \qquad \ka(x,t) \sim |x|\tan\nu(t).
\ee
Assume $f\in\Lcal_{2,\beta}(m)$. We claim that $BR(z, f)^*\in \Lcal_{2,\beta}(m)$. The proof is completely analogous to the proof of Lemma \ref{lem:BRbasic}.  Keeping the same notation, for a given $z = z_+\in \Ga^+$ the only difference is the region $q = q_-\in \Ga^-$ with $u\in I_c(x)$. In fact, we have
$$
\frac{z_+'}{z_+- q_-} = s(x,  u) + r(x,  u), \qquad s(x, u)= \frac{1}{(x - u) + i \rho(u)}, \qquad r(x,u) = O(x^{-1})
$$
(cf. \eqref{eq:aux-I} and recall the notation $\rho(x) = 2\ka(x,t)$). However, using \eqref{eq:ka-corner}, we see that
$$
s(x, u) = O(x^{-1}), \qquad x\in I_\de, \, u\in I_c(x).
$$
In particular, the corresponding integral is bounded in $\Lcal_{2,\beta}(m)$ by Hardy's inequality. Moreover, we now have 
$$
\pa_x s(x, u) = O(x^{-2}), \qquad \pa_x r(x, u) = O(x^{-2}), \qquad x\in I_\de, \, u\in I_c(x),
$$  
(Lemma \ref{lem:errorrk} can be easily adapted to the present case) and we therefore lose $O(m)$ when we put a derivative on the kernel of $BR(z,f)^*iz_s - Hf$. In particular, assuming for simplicity $k = 2$, it is not difficults to see that $BR(z,f)^*iz_s - Hf \in \Lcal^1_{2,\beta + 1}(m)$ (more details can be found in the first part of the proof of Lemma \ref{lem:BRcancel} up to equation \eqref{eq:ders-IlIr}; it applies word-for-word to the present case and can easily be generalized to higher $k$). In particular, the conclusions of Lemmas \ref{lem:BRcancel} and \ref{lem:struct1} are an easy consequence of \eqref{eq:corner-der-on-ker} combined with Lemmas \ref{lem:trHilbert} and \ref{lem:derFphi} from the Appendix.
\end{proof}

\begin{lemma}\label{lem:corner-ipp}
Let $0 < \beta + 1/2 < 1$. In general, when $f\in H^1_{\beta}(m)$, we only have
$$
BR(z,f)^*\in\Lcal^1_{2,\beta + 1}(m).
$$
However, when $f\in \Lcal^1_{2,\beta}(m)$ we have
$$
BR_{-1}(z,f)^*\in\Lcal^1_{2,\beta}(m), \qquad BR(z,f)^*\in H^1_\beta(m)
$$
and, more generally, we have $BR_{-n}(z,f)^*\in \Lcal^n_{2,\beta}(m)$ provided $f\in \Lcal^n_{2,\beta}(m)$ where 
$$
\pa_s BR_{-n}(z,f)^* = z_s BR_{- n + 1}(z,D_sf)^*, \qquad n\geq 1.
$$
Similarly, when $f\in \Lcal^{n}_{2,\beta + n}(m)$, we have $BR(z, f)^*\in \Lcal^n_{2,\beta + n}(m)$.
\end{lemma}

\begin{proof}
	Integration by parts gives 
	$$
	\aligned
	BR(z, f)^*_s = z_s BR(z, D_sf)^* + \frac{z_s(\al)}{2\pi i z(\al)}\bigg(\frac{f(-\al_*^-)}{z_s(-\al_*^-)} - \frac{f(\al_*^+)}{z_s(\al_*^+)} + \frac{f(\al_*^-)}{z_s(\al_*^-)} - \frac{f(-\al_*^+)}{z_s(-\al_*^+)}\bigg),
	\endaligned
	$$
	where we have used the notation $f(\al_*^-) := \lim_{\al \nearrow \al_*} f(\al)$ resp. $f(\al_*^+) := \lim_{\al \searrow \al_*} f(\al)$ etc. 
	
	Assume e.g. that $f$ is odd w.r.t. $x$-axis. Then, we have 
	$$
	\Big(\frac{f(-\al_*^-)}{z_s(-\al_*^-)} - \frac{f(\al_*^+)}{z_s(\al_*^+)}\Big) + \Big(\frac{f(\al_*^-)}{z_s(\al_*^-)} - \frac{f(-\al_*^+)}{z_s(-\al_*^+)}\Big)  =  2i(f(\al_*^+) + f(\al_*^-))\sin \nu 
	$$
	where we have used that
	$$
	z_s(-\al_*^\pm) = -e^{\pm i\nu}, \qquad z_s(\al_*^\pm) = e^{\pm i\nu}.  
	$$
	In particular, the remaining term vanishes if and only if $f(\al_*^+) + f(\al_*^-) = 0$ or  $\nu = 0$. 
	
	In particular, when $f\in \Lcal^1_{2,  \beta}(m)$, boundary terms coming from the integration limits vanish regardless of any symmetry assumptions as $f = O(m^{1 - (\beta + 1/2)})$. The same is true when $f\in \Lcal^1_{2,  \beta + 1}(m)$, since a correction (cf.~equation \eqref{eq:correction-ipp}) is necessary when integrating by parts and  $zf\in \Lcal^1_{2,  \beta}(m)$.
\end{proof}

\begin{remark}
Assume $\om\in H^2_{\beta, 0}(m)$. Then, using symmetry and the proof of Lemma \ref{lem:corner-ipp}, it is not difficult to see that
$$
BR(z,\om)^*_{ss} =  4\om_s(\al_*) \sin 2\nu\, \frac{z^2_s}{2\pi z} + \Lcal_{2,\beta}(m),
$$
which without further cancellations only implies $\theta_{ts}\in\Lcal_{2,\beta + 1}(m)$. This is why the additional assumption $\om_s(\al_*) = 0$ (i.e.~$\om\in \Lcal^2_{2,\beta}(m)$) is necessary.
\end{remark}

The results of Section \ref{s:estDerivatives} can be summarized in the following lemma:

\begin{lemma}\label{lem:corner-preliminary} Let $k\geq 2$ and let
\be\label{eq:norm-energy}
\|\om\|_{\Lcal^{k + 1}_{2,\beta + (k-1)}(m)} + \|\theta\|_{H^{k+1}_{\beta + k}(m)} + \|\log|z_\al|\|_{H^{k}_{\beta + k-1}(m)} + \| z \|_\Fcal < \infty
\ee
	Then, we have $z_t\in H^{k+1}_{\beta + (k-1),0}(m)$ and the relation \eqref{eq:ztdiff} holds. More precisely, we can write
	$$
	z_t^* = -z b + \Lcal^{k+1}_{2,\beta + (k-1)}(m), 
	$$
	where $b$ is the continuous linear functional defined in \eqref{eq:RT-b}. In particular, we have 
	\be\label{eq:th_t-corner}
	\theta_t = \Im(z_s^2)b + \Lcal^k_{2,\beta + (k-1)}(m), \qquad \varphi_{s} = -\Re(z_s^2)b + \Lcal^k_{2, \beta + (k-1)}(m). 
	\ee
	In other words $\theta_t$ and $\varphi_s = (\log|z_\al|)_t$ belong to $H^{k}_{\beta + (k - 1)}(m)$ and we can write
	\be\label{eq:th_ts-corner}
	\theta_{ts} = H(\varphi_{ss})  + \Rcal, \quad \quad   \Rcal\in H^{k}_{\beta + k}(m).
	\ee	
	Moreover, there exists a unique $\om_t\in \Lcal^{k}_{2,\beta + (k-2)}(m)$ solution of \eqref{eq:omt}, which further implies 
	$$
	\si \in H^{k+1}_{\beta + (k-1),0}(m), \qquad \si  \sim m, 
	$$ 
	where the lower bound holds if and only if $b \neq 0$. Furthermore, we have $\varphi_{st}\in H^k_{\beta + (k-1)}(m)$ and 
	$$
	\si_t = O(m).
	$$
\end{lemma}

\begin{proof}
	Let $k = 2$. Since $\om\in\Lcal^3_{2,\beta+1}(m)$, we use formula \eqref{def:BR-k} to write 
	$$
	z^*_t = -b_0(\om) - zb_1(\om) + BR_{-2}(z, \om)^*, 
	$$
	where complex-valued continuous linear functionals $b_i(\om)$ have been defined in \eqref{def:bi}. By symmetry, we have $b_0(\om) = 0$ and $b:= b_1(\om)\in \Rbb$, hence 
	\be\label{eq:zt-corner}
	z_t^* = -zb + \Lcal^3_{2,\beta + 1}(m) 
	\ee
	as required. Since $z^*_{ts}z_s = \varphi_s - i\theta_t$, we have
	$$
	\theta_t = \Im(z_s^2)b + \Lcal^2_{2,\beta + 1}(m), \qquad \varphi_s = -\Re(z_s^2)b + \Lcal^2_{2,\beta + 1}(m). 
	$$ 
	In particular, $\theta_t$ has a jump when crossing $\al = \pm\al_*$ whenever $b \neq 0$. On the other hand, formula \eqref{eq:th_ts-corner} follows as in the proof of Lemma \ref{lem:param} (see also Lemma \ref{lem:BRbasic-corner}). 
	
	To find $\om_t\in \Lcal^2_{2,\beta}(m)$, we need to invert 
	\be\label{eq:omt-corner}
	\om_t  + 2 BR(z, \om_t)\cdot z_s = F,
	\ee 
	where recall that $F \equiv F(z,\om)$ is given by $F(z,\om) = -2[(BR(z, \om)_t - BR(z, \om_{t})]\cdot z_s$ and that
	$$
	BR(z, \om)_t^* - BR(z, \om_{t})^*=  [z_t, BR(z, \cdot)^*] D_{s}\om  - i BR(z, \om \theta_t)^*
	$$
	cf. \eqref{eq:BRt}. We claim that 	
	\be\label{eq:F-corner}
	F = c\Re(z_s z) +  \Lcal^2_{2,  \beta}(m).
	\ee
	Indeed, since $\om \theta_t\in \Lcal^2_{2,\beta}(m)$ we can write 
	\be\label{eq:BRom-tht-corner}
	\aligned
	-i z_sBR(z, \om\theta_t)^* &=  izz_s b_1(\om\theta_t) + z_sBR_{-2}(z, \om\theta_t)^*z_s\\
	 &=- zz_s \Im b_1(\om\theta_t) + \Lcal^2_{2,\beta}(m)
	\endaligned
	\ee
	(the first-order correction term $b_0(\om\theta_t)$ and the real part of $b_1(\om\theta_t)$ vanish by symmetry since $\om\theta_t$ is even w.r.t. both axes). As for the commutator, we have
	$$
	BR(z,D_s\om)^* = -b_0(D_s\om) + BR_{-1}(z, D_s\om)^*, 
	$$
	respectively
	$$
	BR(z,z_t D_s \om)^* = -b_0(z_t D_s\om) - z b_1(z_t D_s\om) + BR_{-2}(z, z_t D_s\om)^*, 
	$$ 
	where by symmetry we have $b_0(z_t D_s\om) = 0$ and $\Im b_1(z_t D_s\om) = 0$ (recall that the real part of $z_tD_s \om$ is odd, while the imaginary is even w.r.t. both axis). In particular, using \eqref{eq:zt-corner} and the fact that $b= b_1(\om) = b_0(D_s\om)$, we obtain
	\be\label{eq:zt_comm-corner}
	z_s[z_t, BR(z, \cdot)^*] D_{s}\om = z^*z_s |b_1(\om)|^2 + zz_s b_1(z_t D_s \om) +  \Lcal^{2}_{2,\beta}(m),
	\ee 
	since $\om\in \Lcal^{3}_{2,\beta + 1}(m)$. Taking the real part of \eqref{eq:BRom-tht-corner} resp. \eqref{eq:zt_comm-corner}, we obtain 
	$$
	\frac{F}{2} = \Re(z_s z) (\Im b_1(\om\theta_t) - \Re b_1(z_t D_s\om)) - z\cdot z_s |b_1(\om)|^2 + \Lcal^2_{2,\beta}(m).
	$$
	Since 
	\be\label{eq:zzs-aux}
	z_s z = (z\cdot z_s) z_s^2 + (z\cdot z^\perp_s) iz_s^2, \qquad z_s z^* = z \cdot z_s - i z\cdot z_s^\perp, \qquad z\cdot z_s^\perp \in\Lcal^2_{2,\beta}(m)
	\ee
	we can further write
	\be\label{eq:z.zs}
	z\cdot z_s =  \Re(z_s z)\frac{1}{\cos 2\nu} + \frac{z\cdot z_s}{\cos 2\nu}\, (\cos 2\nu - \Re (z_s^2)) +\Lcal^2_{2,\beta}(m). 
	\ee
	Noting that $\Re(z_s^2)$ is continuous with value $\cos 2\nu$ at $\al = \pm\al_*$, we conclude that
	$$
	\frac{c}{2} = \Im b_1(\om\theta_t) - \Re b_1(z_t D_s\om) -\frac{1}{\cos 2\nu}\, b^2
	$$
   In particular, by Theorem \ref{thm:inverseOpII} below, there exists $\om_t\in\Lcal_{2,  \beta - 2}(m)$ a unique solution of \eqref{eq:omt-corner}. Moreover, we must have
   \be\label{eq:constant}
   -\Re b_1(\om_t) = \frac{c}{2}
   \ee
   (this follows from \eqref{eq:omt-corner}, since $BR(z,\om_t)\cdot z_s = -\Re(zz_s)b_1(\om_t) + \Lcal_{2, \beta-2}(m)$). In order to construct derivatives of $\om_{t}$, we rewrite the equation \eqref{eq:omt-corner} as 
   $$
   \om_t = F(z,\om) - BR(z,\om_t)\cdot z_s
   $$
   and we claim that the r.h.s. belongs to $\Lcal^1_{2,\beta - 1}(m)$, where we can estimate the corresponding norm in terms of $\Lcal_{2,\beta - 2}(m)$-norm of $\om_t$ and \eqref{eq:norm-energy}. In fact, it is enough to show
   \be\label{eq:claim1}
   \pa_s (BR(z,\om_t)\cdot z_s) = -\Re(z_s^2 b_1(\om_t)) +  \Lcal_{2,\beta - 1}(m).
   \ee
   Then, we can conclude that $\om_{ts}$ exists and belongs to $\Lcal_{2,\beta - 1}(m)$ by \eqref{eq:F-corner} and \eqref{eq:constant}. The claim \eqref{eq:claim1} basically follows from Lemma \ref{lem:BRbasic-corner}, however some care is nedeed in order to isolate the constant term. In fact, writing as usual $BR(z,\om_t)^* = -b_0(\om_t) - zb_1(\om_t) + BR_{-2}(z, \om_t)^*$, then using the definition of $BR_{-2}(z,\om_t)^*$, we can write  
   $$
   BR_{-2}(z, \om_t)^*z_s =-i z^2 \left(BR\left(z, \frac{\om_t}{z^2}\right)^*iz_s - H\left(\frac{\om_t}{z^2}\right)\right) - i \left(z^2 H\left(\frac{\om_t}{z^2}\right) - H\om_t\right) - iH\om_t. 
   $$
   Taking the real part, the pure Hilbert transform term vanishes, while the terms in brackets belong to $\Lcal^1_{2,\beta - 1}(m)$ by Lemmas \ref{lem:BRbasic-corner} and \ref{lem:derFphi} (cf. Appendix) respectively. In particular, we have $\om_{t}\in \Lcal^1_{2,\beta - 1}(m)$. 
   
   It remains to construct $\om_{tss}$. Unfortunately, we cannot use Lemma \ref{lem:BRbasic-corner} directly once again. However, as we now have $\om_{ts}$, we can integrate by parts to conclude that actually 
   $$
   \pa_s (BR(z,\om_t)\cdot z_s) = BR(z,\om_t)_s\cdot z_s + \theta_s BR(z,\om_t)\cdot z_s^\perp, \qquad BR(z,\om_t)_s^* = z_s BR(z, D_s\om_t)^*
   $$
   and that we can write 
   $$
   \om_{ts} - 2BR(z, \om_{ts})\cdot z_s = F_s - F_1, 
   $$
   where 
   $$
   F_1 := 2(\pa_s (BR(z,\om_t)\cdot z_s) - BR(z,\om_{ts})\cdot z_s) 
   $$
   If $\pa_s F_1\in \Lcal_{2,\beta}(m)$, then we can proceed as in the construction of $\om_{ts}$ to show that $\om_{tss}\in \Lcal_{2,\beta}(m)$ exists. Indeed, let $F_1 = I_1 + I_2$, with $I_1$, $I_2$ as in \eqref{def:F1-in-detail}. We have
   $$
   \aligned
   I_1 =& \Re\left(\frac{1}{2\pi i}\int_{-\pi}^\pi\left(\frac{z_s(\al)}{z_s(\al')} + 1\right)\om_{ts}(\al')\left(\frac{z_s(\al)}{z(\al)-z(\al')}-\frac{1}{\al-\al'}\right)ds_{\al'} \right) \\
    &+  \frac{1}{2\pi}\int_{-\pi}^\pi\Im\left(\frac{z_s(\al)}{z_s(\al')} + 1\right)\frac{\om_{ts}(\al')}{\al-\al'}\, ds_{\al'}, 
   \endaligned
   $$
   where we can put a derivative on the kernel of the first term to conclude $\Lcal_{2,\beta}(m)$ provided $\om_{ts}\in \Lcal_{2,\beta - 1}(m)$. Similarly, the derivative of the second term belongs to $ \Lcal_{2,\beta}(m)$. In fact, the only problematic region is $\al\sim \al'$, but there 
   $$
   \aligned
   \Im\left(\frac{z_s(\al)}{z_s(\al')} + 1\right) &=  \Im(z_s(\al)z^*_s(\al')) =  \Im((z_s(\al)-z_s(\al'))z^*_s(\al'))\\
   & =\theta_{\al}(\al')(\al-\al') + O(m(\al)^{-(\lambda + 2)}|\al - \al'|^2) 
   \endaligned
   $$  
   where $\lambda = 1 - (\beta + 1/2)$ and we are finished. As for $I_2$, we have $z_sD_s\om_t -\om_{ts} = -i\om_t \theta_s \in \Lcal^1_{2,\beta - \lambda}(m)\subseteq \Lcal^1_{2,\beta}(m)$, hence $\pa_s I_2\in \Lcal^1_{2,\beta}(m)$. In particular, we have $F_s - F_1\in\Lcal^1_{2,\beta}(m)$. We omit further details.
   
   Once $\om_t\in \Lcal^2_{2,\beta}(m)$ has been constructed, it is not difficult to see that $\si\in H^2_{\beta,0}(m)$; in fact combining equation \eqref{eq:si} with \eqref{eq:BRt} and \eqref{eq:BRtn}, then taking the negative of the imaginary part of  \eqref{eq:BRom-tht-corner} resp. \eqref{eq:zt_comm-corner} we conclude 
   \be\label{eq:si-corner}
   \aligned
   \si &= \Im(z_s z) \left(\Re b_1(\om_t) + \Im b_1(\om\theta_t) - \Re b_1(z_tD_s\om)\right) - z\cdot z_s^\perp \, b^2 + \Lcal^2_{2,\beta}(m) \\ 
   &=b^2\left(\frac{1}{\cos 2\nu}\Im(z_s z) - z\cdot z_s^\perp \right) + \Lcal^2_{2,\beta}(m)\\
   &= b^2 \, z \cdot z_s \, \frac{\Im (z^2_s)}{\cos 2\nu} + b^2 z\cdot z^\perp_s \left(\frac{\Re (z_s^2)}{\cos 2\nu} - 1\right)  + \Lcal^2_{2,\beta}(m)
   \endaligned
   \ee 
   In particular, we have 
   $$
   \si = b^2 \, z \cdot z_s\, \frac{\Im (z^2_s)}{\cos 2\nu} + O(m^{1+\lambda}), \qquad \lambda = 1 - (\beta + 1/2).
   $$
   In order to show $\pa_s^3\si\in \Lcal_{2,\beta + 1}(m)$, we can proceed as in Lemma \ref{lem:si}. We omit the details.

   Moreover, it is a matter of straightforward calculation to show that $\varphi_{st}\in H^2_{\beta + 1}(m)$ holds (recall that the time derivative of $\varphi_s$ is given by \eqref{eq:varphi_st}).

   It remains to show that
   \be\label{eq:sit-corner}
   \si_t = O(m).
   \ee    
   First note that  
   $$
   z_{tt} = \si z_s^\perp  \in H^{3}_{\beta + 1,0}(m), \qquad \theta_{tt} = \si_s - 2\varphi_s\theta_t\in H^{2}_{\beta + 1}(m).
   $$
   (cf. proof of Lemma \ref{lem:paramt} for the details on the derivation of the formulas). In these formulas, using the equation for $\si$, it is not difficult to see that we can write 
   \be\label{eq:z_tt-asym-corner}
   z_{tt} = i z b^2 \frac{\Im(z_s^2)}{\cos2\nu} + \Lcal^2_{2,\beta}(m).  
   \ee 
   If we can show the following version of the auxiliary Lemma \ref{lem:BRtime2} 
	\be\label{eq:BRtime2-corner}
	[BR(z,\om)^*_{tt} - BR(z, \om_t)^*_t] = O(m),
	\ee
	with the tangential component continuous and of the form
	\be\label{eq:BRtime2-tangential-corner}
	[BR(z,\om)^*_{tt} - BR(z, \om_t)^*_t] \cdot z_s = c \Re(z_s z) + \Lcal^1_{2,\beta - 1}(m)
	\ee
	for some $c\in\Rbb$, Theorem \ref{thm:inverseOpII} and Lemma \ref{lem:BRbasic-corner} imply there exists $\om_{tt}\in \Lcal^1_{2,\beta - 1}(m)$, a solution of 
	$$
	\om_{tt}  + 2 BR(z, \om_{tt})\cdot z_s = G,
	$$ 
	where   
	\be\label{eq:def-G-corner}
	\aligned
	G &= - 2[(BR(z, \om_{t})\cdot z_s)_t - BR(z, \om_{tt})\cdot z_s] + F_t\\ 
	 & = -2\theta_t \si + \om\theta_t^2 -2 [BR(z, \om_t)_t - BR(z, \om_{tt})]\cdot z_s - 2[BR(z, \om)_t - BR(z, \om_t)]_t\cdot z_s \\
	&= c\Re(z_s z) + \Lcal^1_{2,\beta-1}(m)
	\endaligned
	\ee
	for some (possibly different) $c\in\Rbb$ (cf. proof of Lemma \ref{lem:omtt} and the construction of $\om_t$ and its derivatives).    
	
	Once $\om_{tt}\in\Lcal^1_{2,\beta - 1}(m)$ has been constructed, the estimate \eqref{eq:sit-corner} follows combining the proof of Lemma \ref{lem:si_t} with the estimate \eqref{eq:BRtime2-corner}, then using the formula \eqref{def:BR-k} together with symmetry assumptions on $BR(z,\om_{tt})^*$. Note that any function that belongs to $\Lcal^1_{2,\beta - 1}(m)$ must be $O(m^{1+\lambda})$. We omit further details.
	
	Finally, let us give some details on \eqref{eq:BRtime2-corner} and \eqref{eq:BRtime2-tangential-corner}. Writing the left hand-side of \eqref{eq:BRtime2-corner} as in the proof of Lemma \ref{lem:BRtime2}, it is not difficult to see that
	$$
	BR(z,\om)^*_{tt} - BR(z, \om_t)^*_t = z_{tt} BR(z, D_s\om)^*  + c_1 z^* + c_2 z + \Lcal^1_{2,\beta - 1}(m).
	$$
    In particular, the estimate \eqref{eq:BRtime2-corner} is straightforward. As for the estimate \eqref{eq:BRtime2-tangential-corner}, we proceed as in the proof of \eqref{eq:F-corner} to conclude   
    $$
    [BR(z,\om)^*_{tt} - BR(z, \om_t)^*_t]\cdot z_s = -\Re(izz_s)b^3 \frac{\Im(z_s^2)}{\cos 2\nu} + c\Re(z_s^2) + \Lcal^1_{2,\beta - 1}(m) ,
    $$  
    for some $c\in\Rbb$, where we have used \eqref{eq:z_tt-asym-corner}. However, using \eqref{eq:zzs-aux} and \eqref{eq:z.zs} we see that actually 
    $$
    -\Re(izz_s) \frac{\Im(z_s^2)}{\cos 2\nu} = \Re(z_s z) (\tan 2\nu)^2 + \Lcal^1_{2,\beta - 1}(m)
    $$    
    and we are finished.
\end{proof}

The lower-order contributions to the energy now read
$$
E_{l}(t)^2 :=  \|\theta\|^2_{H^2_{\beta + 1}(m)} +  \|\log |z_\al|\|^2_{H^1_{\beta}(m)} + \|\om\|^2_{\Lcal^2_{2,\beta}(m)} + \|z\|^2_\Fcal, 
$$
with the higher order contributions being given by 
$$
\aligned
E_{h}^{k}(t)^2 := &\|\sqrt{\si}\pa_s^{k+1}\theta\|^2_{2, \beta + (k - 1/2)} + \|\Lambda^{1/2}(m^{\lambda + (k-1)}\pa^k_s\varphi_s)\|^2_{2, \beta-\lambda + (k -1/2)} \\ 
&+ \|\pa_s^k \varphi_s\|^2_{2,\beta+ (k-1)} + \|\pa_s^k\log |z_\al|\|^2_{2, \beta + (k-1)},
\endaligned
$$
and $\lambda$ chosen in such a way that $-1/2 < (\beta - \lambda) + 1/2 < 0$.

Given a sufficiently regular solution  $(z,\om)$ of \eqref{eqs:sys} such that \eqref{eq:conditions-corner} is satisfied, the energy functional defined as
$$
E_k(t)^2 = E_l(t)^2 + \sum_{i = 2}^k E^i_h(t)^2, \qquad k\geq 2.
$$
satisfies an a-priori energy estimate of the type \eqref{eq:energy}. The proof is very similar to the proof of Lemma \ref{lem:apriori}. Note that the analog of Lemma \ref{lem:omsss} holds as well, i.e. we can replace the norm of $\pa_s^k\om_s$ by the norm of the corresponding $\pa_s^k\varphi_s$. The regularization of the evolution equations and the existence of solutions follow as in Sections \ref{ss.regularization} and \ref{ss.existence}. We omit further details.


\section{The inverse operator}\label{s:inverseOp}

In this section, we consider the invertibility properties of the singular integral operators
\be\label{eq:ops}
\phi \mapsto \phi \pm 2 BR(z, \phi) \cdot z_s 
\ee
in various weighted Sobolev spaces when the interface  has cusps or corners. These correspond to solving the interior respectively exterior Neumann problem in $\Om$. Indeed, let us introduce the notation 
$$
S\phi :=  2 BR(z, \phi) \cdot z_s
$$  
and let $V\phi$ be the single-layer potential   
$$
V\phi(z) := \frac{1}{\pi} \int_\Ga  \phi(\al) \log|z- z(\al)|^{-1}ds_\al, \quad z\in \Rbb^2\setminus \Ga.
$$
It is known that $V\phi$ is continuous when crossing the boundary, while its normal derivative exhibits a jump. More precisely, for $z\in \Ga\setminus\{z_*\}$ we have 
$$
(I \pm S)\phi = \pm \pa_{n^\pm} V\phi,
$$
where $\pa_{n^\pm}$ denotes the normal derivative of $V\phi$ approaching $\Ga$ from the interior of $\Om$ or, respectively, of its complement $\Om^c$, with the unit normal vector chosen to be  inward directed. In particular, we have 
$$
\pa_{n^+} V\phi-\pa_{n^-} V\phi = 2\phi.
$$

\subsection{Domain with a cusp singularity}

The invertibility properties of operators \eqref{eq:ops} associated to a bounded, connected domain with exterior (or interior) cusp can be found in \cite{MazSo}. No symmetry assumptions are made and, using conformal maps, the operators \eqref{eq:ops} are shown to be one-to-one and onto from $\Lcal_{p, \ga + 1}(m)$ to a strictly smaller subspace of $\Lcal_{p, \ga + 1}(m)$, where the operator $S$ is to be replaced by the correction
$$
S_{+1}\phi := 2BR_{+1}(z, \phi)\cdot z_s.
$$ 
To get some intuition about what the additional restrictions on the image look like see Lemma \ref{lem:BRdiff}. For the necessary modifications to the conformal maps in the present setting (i.e. when $\Om$ is a union of two domains each with an outward cusp smoothly connected through the common tip) see \cite{CEG2}, where the operator, corresponding to the exterior Neumann problem here, is inverted on the subspace of $\Lcal_{p, \ga}(m)$-functions which are even with respect to the $x$-axis. Here, although not in full generality, we extend these results to $H^1_{\ga + k}(m)$ where $k=0,1,2$. More precisely, we have:      

\begin{theorem}\label{thm:inverseOpI}
	Let $0< \ga + 1/2 < 1$ where $\ga + 1/2 \neq \mu$. Let $c\in\Rbb$ and let $\psi\in \Lcal^1_{2,\ga}(m)$ odd or $\psi \in H^1_{\ga}(m)$ even with respect to both axes. Then, there exists a unique $\phi\in \Lcal^1_{2,\ga}(m)$ odd respectively $\phi\in H^1_\ga(m)$ even with respect to both axes, solution of the interior resp. exterior Neumann problem, that is
	$$	
	(I \pm S)\phi = \psi.
	$$
	Similarly, given $\psi\in \Lcal^1_{2,\ga + k}(m)$ odd/even with respect to the $x$-axis with $k=1,2$, there exists a unique $\phi\in \Lcal^1_{2,\ga + k}(m)$ of the same parity, solution of $(I \pm S)\phi = \psi$, where the operator $S$ is to be replaced by $S_{+1}$, when $k=2$. As long as it makes sense (i.e. for $k=0,1$), we assume $\int_\Ga \psi ds = 0$.

	All the conclusions (under the corresponding assumptions) remain true if $\psi\in\Lcal_{2,\ga + k - 1}(m)$ (instead of $\Lcal^1_{2,\ga + k}(m)$) and if $\psi\in\Rbb \oplus \Lcal_{2,\ga  - 1}(m)$ (instead of $H^1_\ga(m)$).	
\end{theorem}

  For brevity, we omit most of the details of the proof and refer the reader to \cite{MazSo} (see also \cite{CEG2}) instead. The uniqueness part is a consequence of the following Proposition, which we state without proof: 

\begin{proposition}\label{prop:injectivity}
Let $0< \ga + 1/2 < 1$. Then, the operators $I \pm S$ are injective on $\Lcal_{2,\ga}(m)$. Similarly, $I \pm S_{+1}$ are injective on $\Lcal_{2,\ga + 1}(m)$. 
\end{proposition}


\begin{proof}[Sketch of the proof of Theorem \ref{thm:inverseOpI}]
	The existence (and continuity of the inverse) part is based on the construction of suitable harmonic functions on $\Om$ respectively $\Om^c$. More precisely, for say the exterior Neumann problem, we construct $u^i, u^e$ such that $u^e$ is harmonic at infinity and
	such that 
	$$
	\begin{cases}
		\Delta u^e = 0& \quad \text{in }\Om^c\\
		\pa_n u^e = \psi& \quad  \text{on }\Ga\setminus\{z_*\}
	\end{cases}
	\qquad \qquad
	\begin{cases}
		\Delta u^i = 0& \quad \text{in }\Om\\
		u^i =u^e& \quad  \text{on }\Ga\setminus\{z_*\}
	\end{cases}
	$$
	By considering complex conjugates, it is enough to solve the corresponding Dirichlet problem. The details of the construction can be found in Propositions \ref{prop:dirichletInt} and \ref{prop:dirichletExt} below. Once $u^i$ and $u^e$ have been constructed, we set 
  	$$
  	\phi := \frac{1}{2}(\psi + \pa_n u^i)
  	$$
  	and consider the harmonic mapping in $\Rbb^2\setminus \Ga$ defined via $V\phi - u^{e}$ in $\Om^c$ respectively $V\phi - u^i$ in $\Om$. This mapping can be harmonically extended to all of $\Rbb^2$ which proves that  
  	$$
  	 (I - S)\phi = \psi.
  	$$ 
    We omit further details.

\end{proof}

\subsubsection{Construction of harmonic functions on $\Om$}

\begin{proposition}\label{prop:dirichletInt} Let $0< \ga + 1/2 < 1$. Then, for any $\phi\in \Lcal^2_{2, \ga + k}(m)$ even with respect to the $x$-axis with $k = 0, 1, 2$, there exists a harmonic extension $U$ of $\phi$ to $\Omega$, i.e. a solution of
\be\label{eq:systemOmDirichlet}
\Delta U = 0, \quad \quad U\big|_\Ga = \phi,
\ee
such that 
\be\label{eq:estVnormal}
\|\pa_n U\|_{\Lcal^1_{2, \ga + k}(m)} \,\lesssim \, \|\phi\|_{\Lcal^2_{2, \ga + k}(m)}.
\ee 
\end{proposition}


Before we prove Proposition \ref{prop:dirichletInt}, we define diffeomorphisms which map the connected components $\Om^\pm := \Om \cap \Rbb^2_\pm$ of $\Om$ to the unbounded strip-like domain $\Gcal$, on which we subsequently solve the resulting elliptic problem (cf. \cite{MunRam}). More precisely, we construct  
$$
H^{-1}_\pm:\Om^\pm \rightarrow \Gcal, 
$$
where $\Gcal \cap \Rbb_-^2$ is of class $\Ccal^{3, \lambda}$ for some $\lambda\in(0,1)$ and
$$
\Gcal \cap \Rbb_+^2 = \{\t x + i\t y \, : \, \t x>0, \, |\t y| < 1/2\}.
$$ 
The same can be accomplished using conformal maps (cf. \cite{MazSo}); first map $\Om$ to a bounded domain via the mapping $z\rightarrow z^{-\mu}$, then go over to $\Gcal$ (basically) via an exponential map. However, these are too rigid for our purposes. Higher order boundary regularity requires control over a sufficient number of derivatives of $x^{-(\mu+1)}\ka(x)$, which our current regularity assumptions \eqref{regTheta}-\eqref{regTheta2} do not provide (in general we do not even have $\Ccal^{1,\lambda}$-boundary; to see how exactly these conditions come into play, cf. \cite{CEG2}). 

Here, we will only provide details on the construction of $H^\pm$ in the neighborhood of the cusp. For further details on the construction see e.g. \cite{MunRam}. So, by possibly making $\de$ smaller, we may assume the interface can be parametrized as $(x, \pm\ka(x))$ for $|x|<4\de$. We concentrate on the left connected component $\Om^-$ and drop the subscript $-$. Let 
$$
\Om_\de := \{(x,y)\in\Rbb^2 \, : \, -\de < x < 0, \, \,  |y| < \ka(x) \}.
$$
On $\Om_{2\de}$, we define
$$
h^{-1}(x) :=  \int_{-2\de}^x\frac{ds}{\rho(s)}, \qquad h^{-1}_2(x,y) :=  \frac{ y}{\rho(x)}
$$ 
and we set  
$$
H^{-1}(x,y) := (h^{-1}(x), \, h^{-1}_2(x,y)), \quad (x,y) \in \Om_{\de}.
$$
More on the properties of the mapping $h^{-1}$ and its inverse can be found in the Appendix following definition \eqref{eq:hb} (where it is defined for the right connected component $\Om^+$). 

It is a matter of straightforward calculation to see that under diffeomorphisms the Laplace operator transforms as
$$
(\Delta U) \circ H \quad \Rightarrow \quad \mathrm{div}(\Abb \na (U \circ H)), 
$$
where   
$$
\Abb \circ H^{-1} = \frac{1}{|\det DH^{-1}|} \, \, DH^{-1} (DH^{-1})^T
$$
is symmetric and positive definite. In fact, on $\Om_\de$ we have  
$$
DH^{-1}(x,y) = \frac{1}{\rho(x)}\left(
				    \begin{array}{cc}
				      1 & 0\\
				    -\frac{\rho'(x)}{\rho(x)} y & 1
				    \end{array}
				 \right),
\qquad  
DH^{-1}\circ H(\t x, \t y) = \frac{1}{h'(\t x)}\left(
						\begin{array}{cc}
						  1 & 0\\
						  -\frac{h''(\t x)}{h'(\t x)} \t y & 1
						\end{array}
					     \right)
$$
and therefore
$$
\Abb(\t x, \t y) = I + \frac{h''(\t x)}{h'(\t x)}\,\t y \left(
			\begin{array}{cc}
                           0 & -1\\		      
                          -1 &  \frac{h''(\t x)}{h'(\t x)}\,\t y
                        \end{array}
		  \right),
$$
where $I$ is the identity matrix and the eigenvalues read 
$$
\lambda_\pm (X(\t x, \t y)) = \left(\frac{X(\t x, \t y)}{2} \pm \sqrt{\frac{X(\t x, \t y)^2}{4} + 1}\right)^2, \quad X(\t x, \t y) := \frac{h''(\t x)}{h'(\t x)}\,\t y.
$$
Since $ |X(\t x, \t y)|=|\rho'(h(\t x))\t y|\leq X_0 := \frac{1}{2}|\rho'(-\de)|$, we have 
$$
\lambda_-(X_0) |u|^2 \leq \langle u, \Abb u\rangle \leq \lambda_+(X_0) |u|^2
$$
for any vector field $u$ on $\Om_\de$. 

\subsubsection*{Variational solutions on $\Gcal$}

We look for weak solutions of the elliptic equation
\be\label{eq:systemGpoisson}
\aligned
-\mathrm{div}(\Abb \na u) &= f \qquad \text{in }\Gcal\\
u &= 0 \qquad \text{on }\pa\Gcal
\endaligned
\ee
in suitable weighted Sobolev spaces. Let $\t x_0 \geq 1 $ be fixed and let $\t m$ be a smooth function on $\Rbb$ such that
$$
\t m(\t x) = 1 + \t x_0, \quad \t x\leq \t x_0, \qquad \qquad \t m(\t x) =  1 + \t x, \quad \t x \geq 2\t x_0
$$
and such that  
$$
\t m(\t x_0) \leq \t m(\t x), \qquad |\t m'(\t x)| \leq c, \qquad \forall \t x \in \Rbb
$$ 
where $c>0$ is some absolute constant independent of $\t x_0$. We set   
$$
\t m(\t x, \t y) := \t m(\t x), \qquad \forall (\t x,\t y)\in\Gcal,
$$
where we abuse the notation and use $\t m$ to denote the resulting weight function on $\Gcal$ as well. Let 
$$
u\in H^{1,\ga}(\Gcal) \quad :\Leftrightarrow \quad \t m^\ga u\in H^1(\Gcal), \quad \ga \in \Rbb
$$
and let $H^{1, \ga}_0(\Gcal)$ denote the closure of $\Ccal^\infty_c(\Gcal)$ in $H^{1,\ga}(\Gcal)$. Then, it is not difficult to see that the Poincare inequality holds on $H^{1,\ga}_0(\Gcal)$. More precisely, we have  
\be\label{est:Poincare}
\|\t u\|_{L^2(\Gcal)} \, \lesssim \, \|\na \t u\|_{L^2(\Gcal)}, \quad  \forall u \in H^{1,\ga}_0(\Gcal),
\ee
where we have set $\t u:= \t m^\ga u$ and the constant can be controlled in terms of $\t y_0:= \sup_{(\t x, \t y) \in \Gcal}|\t y|$. 

\begin{lemma}\label{lem:weaksolutions}Let $\ga \in \Rbb$. Then, for suitably chosen $\t x_0$ depending only on $\ga$ and the coefficients of $\Abb$, the bilinear form 
$$
\aligned
&a: H^{1,\ga}_0(\Gcal)\times H^{1,-\ga}_0(\Gcal)\rightarrow \Rbb\\
&a(u, v) := \int_\Gcal \Abb \na u\cdot \na v
\endaligned
$$
satisfies the assumptions of the Lax-Milgram Lemma (cf. Appendix, Lemma \ref{lem:weightedLaxMilgram}). In particular, for any $f \in (H^{1,-\ga}_0(\Gcal))'$, there exists a unique $u\in H^{1,\ga}_0(\Gcal)$ such that 
$$
a(u, v) \, = \,\langle f, v\rangle, \quad \forall v \in H^{1,-\ga}_0(\Gcal).
$$
Moreover, for $k = 0,1$ we have
\be\label{eqs:elliptic-reg}
\|u\|_{H^{k+2, \ga}(\Gcal)} \, \lesssim \, \|f\|_{H^{k,\ga}(\Gcal)}.
\ee
\end{lemma}

\begin{proof}
We first show that $a:H^{1,\ga}(\Gcal)\times H^{1,-\ga}(\Gcal)\rightarrow \Rbb$ is continuous. 
Indeed, setting
$$
\t u := \t m^\ga u, \qquad \t v = \t m^{-\ga}v,
$$
we have, by definition $\t u, \t v \in H^1(\Gcal)$ and  
$$
\|u\|_{H^{1,\ga}} = \|\t u\|_{H^1}, \qquad \|v\|_{H^{1,-\ga}} = \|\t v\|_{H^1}.
$$
Moreover, we have
\be\label{eq:auxA}
\Abb \na u\cdot \na v = \Abb \na \t u\cdot \na \t v + a_{11}\pa_{\t x} \t m^\ga\pa_{\t x} \t m^{-\ga}\t u \t v + \t m^{-\ga}\pa_{\t x} \t m^\ga  \, (\v a\cdot \na \t u)\t v + \t m^\ga\pa_{\t x} \t m^{-\ga} \t u \, \v a \cdot \na \t v, 
\ee
where we have set $\v a =(a_{11}, a_{12})$ and we have used that $\Abb$ is symmetric. 
Since 
$$
a_{ij} \in L^\infty, \qquad \pa_{\t x} \t m^\ga\pa_{\t x} \t m^{-\ga} = -\ga^2\, \frac{\t m_{\t x}^2}{\t m^2} \in L^\infty \qquad \t m^{\mp\ga}\pa_{\t x} \t m^{\pm\ga} = \pm\,\ga\, \frac{\t m_{\t x}}{\t m} \in L^\infty, 
$$
the estimate \eqref{eq:acontinuity} follows from Lemma \ref{lem:weightedLaxMilgram} in the Appendix.

We now show that $a:H^{1,\ga}_0(\Gcal)\times H^{1,-\ga}_0(\Gcal)\rightarrow \Rbb$ is coercive. More precisely, we only prove estimate \eqref{eq:acoercivity1}  in Lemma \ref{lem:weightedLaxMilgram}, and the other one follows analogously. 
Let $u\in H_0^{1,\ga}(\Gcal)$ be fixed and let
$$
v :=\t  m^{2\ga} u \in H_0^{1,-\ga}(\Gcal).
$$
Then $\t v = \t u$, the mixed terms in \eqref{eq:auxA} cancel out and 
$$
a(u, v) = a(\t u, \t u) - \ga^2\int_{\Gcal} \frac{1}{\t m^2} \, \t m_{\t x}^2\, a_{11}|\t u|^2.
$$
Since $\Abb$ is symmetric and positive-definite, we know 
$$
a(\t u, \t u) \, \geq \, \min_{(\t x, \t y)\in \Gcal} \lambda(\t x, \t y) \, \|\na \t u\|^2_{L^2(\Gcal)}
$$
while, by construction of the weight function, we have 
$$
\int_{\Gcal} \frac{1}{\t m^2}\, \t m_{\t x}^2\, a_{11}|\t u|^2 \, \lesssim \, \frac{1}{\t m(\t x_0)^2} \int_{\Gcal \cap \{\t x >\t x_0\}} \t u^2 \, \lesssim \, \frac{1}{\t m(\t x_0)^2}\,  \|\na \t u\|^2_{L^2(\Gcal)},
$$
where in the last step we have used the Poincare inequality on the half strip $\Gcal_0 := \Gcal \cap \{\t x > \t x_0\}$.
At this point, we finally determine $\t x_0$. In fact, choosing $\t x_0 \geq 1$ such that 
$$
\frac{\ga^2}{\t m(\t x_0)^2} \, \lesssim \, \frac{1}{2}  \min_{(\t x, \t y)\in \Gcal} \lambda(\t x, \t y),
$$
we have 
$$
a(u, v) \, \geq \, \frac{1}{2} \min_{(\t x, \t y)\in \Gcal} \lambda(\t x, \t y) \, \|\na \t u\|^2_{L^2(\Gcal)} \, \gtrsim \, \|u\|^2_{H^{1, \ga}(\Gcal)} 
$$
where in the last step we have used the Poincare inequality on $\Gcal$ once again. In particular, estimate \eqref{eq:acoercivity1} follows.

Finally \eqref{eqs:elliptic-reg} follows by standard elliptic estimates. The coefficients of $\Abb$ are two times continuously differentiable up to the boundary of $\Gcal$ (which is of class $\Ccal^{3, \lambda}$).
\end{proof}

\subsubsection*{Variational solutions of the original problem}

\begin{lemma}\label{lem:extension}
Let $\ga \in \Rbb$ and let $g \in \Lcal^2_{2,\ga}(\pa\Gcal)$ be even. Then, there exists an extension $G\in \Lcal^{2}_{2, \ga}(\Gcal)$ of $g$ such that 
\be\label{eq:divAbbG}
-\mathrm{div}(\Abb\na G) \in H^{1,\ga}(\Gcal).
\ee
\end{lemma}
 
\begin{proof} First note that it is enough to consider $\ga = 2$. Indeed, given $g\in \Lcal^2_{2,\ga}(\pa \Gcal)$ with $\ga \neq 2$, we have $\t g := \t m^{\ga - 2}g\in \Lcal^{2}_{2,2}(\pa\Gcal)$. Then, assuming $\t G\in \Lcal^2_{2,2}(\Gcal)$ is an extension of $\t g$ such that \eqref{eq:divAbbG} is satisfied with $\ga = 2$, it is not difficult to see that   
$$
G := \t m^{2 - \ga} \t G \in \Lcal^2_{2, \ga}(\Gcal)
$$
is the required extension of $g$.
	
So let $\ga = 2$. To simplify notation, we drop the tilde and denote coordinates by $(x,y)$ instead of $(\t x, \t y)$ and assume $g\equiv 0$, when $x < 0$. On the half-strip $\Gcal\cap \Rbb^2_+$, we have 
$$
-\mathrm{div} (\Abb \na G) = \pa_{x}^2 G + (1 + X^2) \pa_{y}^2 G - 2X \pa_{x}\pa_{y} G - \pa_{y} X \pa_{x} G + (\pa_{y} X^2 - \pa_{x} X) \pa_{y} G, 
$$
where recall that 
$$
X(x, y) = F(x)y, \qquad F(x) := \frac{h''( x)}{h'(x)}.
$$
The assumptions \eqref{regTheta} on the regularity of the interface imply $F \in  \Lcal^3_{2, 2 + (1/2 + \lambda/\mu)}(\Rbb_+)$, where we have set $\lambda := 1 - (\beta + 1/2)$. Moreover, we have
\be\label{eq:Xasymptotics}
F(x) = O(\langle x\rangle^{-1}), \qquad F'(x) = O(\langle x\rangle^{-2}), \qquad F''(x) = O(\langle x\rangle^{-\lambda/\mu - 2}),
\ee
where $\langle x\rangle := (1 + x^2)^{1/2}$ and we have used that $\t m(x) \sim \langle x\rangle$ when $x>0$. To see this, note that $F(x) = \rho'(h(x))$ and use properties of the mapping $h$ given in the Appendix following definition \eqref{eq:hb}. We extend $F$ by zero to the negative real line, smooth it out in a small neighborhood of zero and consider the operator    
$$
L = \pa_{x}^2  +  (1 + F(x)^2y^2) \pa_{y}^2 - 2F(x)y \, \pa_{x}\pa_{y}.
$$
as defined on the horizontal strip $\Pi = \{(x,y)\, : \, |y| < 1/2\}$. We consider all of the weighted Sobolev spaces on $\Pi$ with respect to the weight $m(x, y): = \langle x\rangle$, e.g. 
$$
G\in \Lcal^2_{2, 2}(\Pi) \quad :\Leftrightarrow\quad \langle x\rangle^{k} \nabla^k G \in L^2(\Pi), \quad k = 0,1,2.
$$
For given $g\in \Lcal^{2}_{2,2}(\Rbb)\subseteq H^2(\Rbb)$, we look for $G\in \Lcal^2_{2,2}(\Pi)$ in the form of an oscillatory integral
$$
G(x, y) := \int_\Rbb e^{ix \xi}a(y; x, \xi)\hat{g}(\xi)d\xi,   
$$
where $\hat{g}$ is the Fourier transform of $g$ and $a(\pm \frac{1}{2}; x, \xi) = 1$. We split $G$ as a sum $G_1 + G_2$, where the symbol of $G_1$ does not depend on $x$ and has compact support with respect to the $\xi$, while the symbol of $G_2$ decays exponentially as $|\xi|\rightarrow \infty$, but is identically equal to zero in the neighborhood of the origin. 

Let $\chi$ be a smooth cut-off such that $\chi(\xi) = 1$ if $|\xi| \leq 1/8$ and $\chi(\xi) = 0$ if $|\xi| \geq 1/4$. We will use the same cut-off further down for the $y$-variable. We start with $G_1$. We freeze the coefficients of $L$ at $y = 0$, then apply the Fourier transform with respect to the $x$ and solve the resulting ordinary differential equation 
$$
a_{yy} - \xi^2  a = 0, \quad \quad a(\pm 1/2; \, \xi) = \chi(\xi), \quad \quad |y|\leq 1/2.
$$ 
A short calculation yields
$$
a(y; \, \xi) := \chi(\xi)\frac{\cosh(\xi y)}{\cosh(\frac{\xi}{2})}. 
$$
Note that modulo a cut-off this is the symbol of the Laplace operator on the strip. 
We claim the corresponding oscillatory integral $G_1$ satisfies
\be\label{eq:G1}
G_1\in \Lcal^2_{2, 2}(\Pi), \quad \quad  LG_1 \in H^{1, 2}(\Pi).
\ee 
The first claim is straightforward. Let us consider a typical term 
$$
\aligned
x\pa_y G_1(x,y) &= \int_\Rbb xe^{i x\xi} \pa_y a(y; \xi) \h g(\xi) d\xi = \int_\Rbb xe^{i x\xi} \tanh(\xi y) a(y; \xi) \xi\h g(\xi) d\xi  \\
& = \int_\Rbb e^{i x\xi} \pa_\xi\left(\tanh(\xi y) a(y; \xi) \wh {g'}(\xi)\right) d\xi,
\endaligned
$$ 
but $\wh {g'}(\xi)\in L^2(\Rbb)$, $\pa_{\xi}\wh {g'}(\xi) \sim \wh {xg'}(\xi) \in L^2(\Rbb)$ and $\xi \mapsto \tanh (\xi y)a(y; \xi)$ is a compactly supported, smooth function, with derivatives uniformly bounded in $y$, for all $|y|\leq 1/2$. In particular, it belongs to $L^2(\Rbb)$ with respect to  $x$ uniformly in $y$ and the claim follows. On the other hand, we have 
$$
LG_1(x,y) = \int_\Rbb e^{ix\xi}\left(F(x)^2y^2a_{yy} - 2i\xi F(x)y \, a_y\right) \h g(\xi) d\xi
$$
and we need to show that $\langle x \rangle^2 LG_1 \in H^1(\Pi)$. Since $\langle x \rangle^2 = 1 + x^2$, it is enough to consider
$$
\aligned
x^2 LG_1(x,y) & = -\int_\Rbb e^{ix\xi}\pa^2_\xi\left((y^2F(x)^2a_{yy}(y;\xi) - 2i\xi yF(x) a_y(y;\xi)) \h g(\xi)\right)d\xi \\
& = - y^2F(x)^2\int_\Rbb e^{ix\xi}\pa^2_\xi (a(y;\xi) \wh {g''}(\xi))d\xi + 2iyF(x)\int_\Rbb e^{ix\xi}\pa^2_\xi \left(\tanh(\xi y) a(y; \xi)\wh {g''}(\xi)\right)d\xi
\endaligned
$$
where we have used that $a_{yy}(y;\xi) \h g(\xi) = a(y;\xi) \wh {g''}(\xi)$. By assumption, we control $\pa_{\xi}^2 \,\wh {g''}(\xi) \sim \wh {x^2g''}(\xi) \in L^2(\Rbb)$ and $a(y; \cdot)$ is smooth, of compact support, with derivatives uniformly bounded with respect to $y$. In particular, we can absorb any additional $\xi$ coming from the derivatives of $a(y;\xi)$ in the symbol, hence the claim follows using \eqref{eq:Xasymptotics} as before.

We now pass to the construction of $G_2$. It will be convenient to consider $G_2$ and the operator $L$ as defined on the upper half-plane $\Hbb$ (after preforming a translation $y \rightarrow y + 1/2$). More precisely, we construct  
\be\label{eq:tG}
G_3\in \Lcal^2_{2, 2}(\Hbb), \quad \quad  LG_3 \in H^{1, 2}(\Hbb), \quad \quad G_3(x, y) = 0, \quad |y|\geq 1/4,
\ee 
where setting 
$$
A(x, y) := 1 + F(x)^2 (y - 1/2)^2, \quad B(x, y) := -2 F(x) (y-1/2)
$$
(for $x\in \Rbb$ and $y\geq 0$) we have 
$$
L =  \pa_{x}^2  +  A(x,y) \pa_{y}^2 + B(x,y) \pa_{x}\pa_{y}. 
$$
To determine the symbol of $G_3$, we again freeze the coefficients at $y = 0$, apply the Fourier transform with respect to $x$ and look for a solution of 
$$
A(x, 0) \t a_{yy} + 2i\xi B(x, 0)\t a_{y} - \xi^2 \t a = 0, \quad \quad \t a(0;  x, \xi) = 1 - \chi(\xi), \quad \quad y \geq 0,
$$
which decays as $|\xi| \rightarrow \infty$. It is not difficult to see that  
$$
\t a(y; x , \xi) := (1 - \chi(\xi))\exp (b(x, \xi)y), \quad \quad b(x, \xi) := \frac{-\mathrm{sgn}(\xi) + i\frac{ F(x)}{2}}{1 + \frac{F( x)^2}{4}}\,\xi
$$
has the required properties. We define
$$
a(y; x , \xi) := \chi(y)\t a(y; x, \xi)
$$ 
and denote the corresponding oscillatory integrals by $G_3$ and $\t G_3$ respectively. In other words we have $G_3(x,y) = \chi(y) \t G_3(x,y)$. We first show $G_3\in \Lcal^2_{2, 2}(\Hbb)$. In fact, we can say more. We actually have $\langle x \rangle^2 G_3(x, y) \in H^2(\Hbb)$ by Theorem \ref{thm:PsiDO} from the Appendix. Indeed, integration by parts gives 
$$
x^2 G_3(x, y) = \chi(y)\int_\Rbb e^{i x \xi} \,  \pa^2_\xi \left(\frac{1 - \chi(\xi)}{\xi^2}e^{b(x, \xi) y}\wh{g''}(\xi)\right)d\xi.
$$
In particular, the symbol can be multiplied by up to $\xi^2$ for the estimate \eqref{eq:a-growth} required by Theorem \ref{thm:PsiDO} to remain true for every fixed $y$. Since these are uniform when $y$ ranges in a compact subset of $\Rbb_+$, the claim follows (as taking $x$ or $y$ derivative may result in multiplication with a factor of order $\xi$). Similarly, one can show $\langle x \rangle^2 \nabla G_3(x, y)\in H^1(\Hbb)$ and $\langle x \rangle^2 \nabla^2 G_3(x, y)\in L^2(\Hbb)$.


It remains to show 
$$
LG_3 = \chi(y)L\t E + A(x,y)(\chi''(y) \t G_3 + 2\chi'(y)\pa_y\t G_3) + B(x,y)\chi'(y)\pa_x \t G_3 \in H^{1,2}(\Hbb).
$$ 
By the previous step (and using \eqref{eq:Xasymptotics}) it is enough to consider
$$
L \t G_3 = \int_\Rbb e^{ix\xi}\left(\t a_{xx} + F(x)^2 y(y - 1 )\t a_{yy} -  2i\xi F(x)y \t a_{y} + i\xi \t a_{x} -  2F(x)(y-1/2) \t a_{xy}\right)\hat{g}(\xi)d\xi
$$
on a neighborhood of $y = 0$. Except for the first term which contains second order derivatives of $F$, all of the remaining terms are similar, they are of the form $O(y\xi^2/\langle x \rangle)$ (or better in terms of powers of $y$ or $1/\langle x \rangle$). Multiplying by $\langle x \rangle^2$, integrating by parts and taking $\pa_x$ or $\pa_y$, the terms which may pose problems are of the form 
$$
I(x,y) = \int_\Rbb e^{ix\xi} c(x)\t a(y;x,\xi) y\xi\, \h f(\xi) d\xi, \quad f = x^k g'', \quad k = 0,1,2,
$$
where $c(x) = O(\langle x \rangle^{-1})$ is a rational function of $F(x)$ possibly multiplied by some polynomial in $F'(x)$ and $F''(x)$ (there are also terms of order $y^2$ or $y^3$, but as we will see in a moment these actually behave even better). There are also terms not multiplied by any power of $y$, but as these appear after taking the $\pa_y$ derivative, there will not be an additional $\xi$, hence these can be estimated using Theorem \ref{thm:PsiDO} as before. To show that $I\in L^2(\Hbb)$ whenever $f\in L^2(\Rbb)$, we write
$$
\aligned
\|I\|^2_{L^2(\Hbb)} &= \int_0^\infty\int_\Rbb\Big|\int_\Rbb e^{ix\xi} c(x)\t a(y;x,\xi) y\xi\, \h{f}(\xi)d\xi\Big|^2dxdy \\
&\leq \int_\Rbb \int_\Rbb |c(x)|^2 \left(\int_0^\infty y^2 |\t a(y;x,\xi)|^2 dy\right)|\xi|^2|\h{f}(\xi)|^2d\xi dx.
\endaligned
$$
Since 
$$
|\t a(y;x,\xi)|^2 = (1-\chi(\xi))^2e^{2\Re b(x,\xi) y},
$$
we obtain
$$
\int_0^\infty y^2 e^{2\Re b(x,\xi) y}dy = -\frac{1}{4}\frac{1}{(\Re b(x,\xi))^3} = \frac{1}{4}\left(\frac{1 + \frac{F(x)^2}{4}}{|\xi|}\right)^3
$$
(where we have used that $\int^\infty_0 y^2e^{\lambda y}dy = \frac{d^2}{d\lambda^2} \int^\infty_0 e^{\lambda y}dy = -2\frac{1}{\lambda^3}$ if $\lambda < 0$). Collecting all the $x$ dependent terms in $c(x) = O(\langle x \rangle^{-1})$ and interchanging the integrals we can write 
$$
\|I\|^2_{L^2(\Hbb)} \, \leq \, \int_\Rbb \left(\int_\Rbb|c(x)|^2dx\right) \frac{1-\chi(\xi)^2}{|\xi|}|\h{f}(\xi)|^2d\xi \lesssim \,  \int_\Rbb\,\frac{1}{\langle \xi \rangle}\, |\h{f}(\xi)|^2 d\xi \lesssim \,\|f\|_2^2
$$
by Plancherel's theorem as required.

Finally, we comment on the term with maximal number of derivatives of $F$. We have $a_{xx} = (yb_{xx} + (yb_x)^2) a$, hence (by the previous step) 
$$
\aligned
\pa_x\left(\langle x\rangle^2\chi(y)\int_\Rbb e^{ix\xi}\t a_{xx}(y; x,\xi) \hat{g}(\xi)d\xi\right) &= \langle x\rangle^2\chi(y)y\int_\Rbb e^{ix\xi}\pa_x^3 b(x,\xi) \t a(y;x,\xi)\hat{g}(\xi)d\xi + L^2(\Hbb)\\
&= \langle x\rangle^2 F'''(x)\chi(y)y\int_\Rbb e^{ix\xi}c(x,\xi) \t a(y;x,\xi)\hat{g}(\xi)d\xi + L^2(\Hbb),
\endaligned
$$
where $c(x,\xi)$ is a nice bounded function of $F(x)$ and its derivatives up to order $2$ which depends linearly on $|\xi|$. In particular, the absolute value of the integral is controlled by the $L^2$-norm of $\wh {g'}$. Since $\langle x\rangle^{2+(1/2 + \lambda/\mu)}F'''\in L^2(\Rbb)$ cf. discussion preceding \eqref{eq:Xasymptotics}, we are finished.

Once we have $G_3$ on the upper half-plane $\Hbb$, we obtain $G_2$ on the horizontal strip $\Pi$ as
$$
G_2(x,y) = G_3(x, y + 1/2) + G_3(x, -y + 1/2)
$$ 
and we are finished.
\end{proof}

\begin{proof}[Proof of Proposition \eqref{prop:dirichletInt}.] Let $\phi\in \Lcal^2_{2,\ga + k}(\Ga)$ be given. As noted above, under the action of $H^{-1}$ the problem \eqref{eq:systemOmDirichlet} transforms to 
\be\label{eq:systemGdirichlet}
\aligned
\mathrm{div}(\Abb \na u) &= 0 \qquad \text{in }\Gcal,\\
u  &=  g   \qquad \text{on }\pa\Gcal,
\endaligned
\ee
where we have set $u = U\circ H$ and $g := \phi\circ h$. It is not difficult to see that
$$
\phi \in \Lcal^2_{2,\ga + k}(\Ga_\pm) \quad \Leftrightarrow \quad g := \phi\circ h \in \Lcal^2_{2,\ga'}(\Rbb_+),
$$
where 
$$
\ga':= - \t \ga + (2 - k)\mu^{-1} + 2, \qquad \t \ga = \mu^{-1}(\ga + 1/2) + 1/2
$$
(cf. the discussion following definition \eqref{eq:hb} in the Appendix). The normal derivative transforms as
$$
(\na U \cdot \v n)\circ h = \frac{1}{|DH \v t|}\, (\Abb \na u)\cdot \v n,
$$
and therefore \eqref{eq:estVnormal} is equivalent to proving 
$$
(\Abb \na u)\cdot \v n\in \Lcal^1_{2,\ga'}(\pa\Gcal)
$$ 
(note that $|DH\v t\, | = O(|h'(\t x)|)$). By Lemma \ref{lem:extension}, there exists an extension $G$ of $g$ such that   
$$
f = - \mathrm{div}(\Abb \na G) \in H^{1,\ga'}(\Gcal)
$$
and by Lemma \ref{lem:weaksolutions} there exists a solution $v$ of the boundary value problem \eqref{eq:systemGpoisson} with the r.h.s. $f = -\text{div}(\Abb \na G)$ such that $v\in H^{3,\ga'}(\Gcal)$. The solution of \eqref{eq:systemGdirichlet} is then given by $u = v + G$. 

It is now not difficult to verify that the estimate \eqref{eq:estVnormal} holds. Indeed, by construction we have $\na G \in \Lcal^1_{2,\ga'}(\pa\Gcal)$, while $v \in H^{2,\ga'}(\pa\Gcal)$ implies $\na v\in \Lcal^1_{2,\ga'}(\pa\Gcal)$. The claim now follows using \eqref{eq:Xasymptotics}, since 
$$
(\Abb \nabla u)\cdot \v n = (1 + X^2)u_y - Xu_x. 
$$ 
\end{proof}

\subsubsection{Construction of harmonic functions on $\Om^c$}

\begin{proposition}\label{prop:dirichletExt}
Let $0< \ga + 1/2 < 1$ and let $\phi\in \Lcal^2_{2, \ga + k}(m)$ where $k=0,1,2$. Then there exists a harmonic extension $U$ of $\phi$ to the exterior domain $\Om^c$ such that $U$ vanishes at infinity and such that
\be\label{eq:estUnormal}
\|\pa_n U\|_{H^1_{\ga + k}(m)}\, \lesssim \, \|\phi\|_{\Lcal^2_{2, \ga + k}(m)}.
\ee
Let $\phi$ be odd with respect to the $y$-axis. If $\phi \in H^1_{\ga,0}(m)$, the normal derivative is locally near the origin given by $\pa_n U = c_\pm \mathrm{sign}(z_2) z_{2s} + \Lcal^1_{2, \ga}(m)$ where $c_\pm = \phi_s(\pm \al_*)$. In particular, when $\phi \in \Lcal^2_{2,\ga}(m)$, we have $\pa_n U \in \Lcal^1_{2, \ga}(m)$.  
\end{proposition}

The proof closely follows \cite{CEG2}, where the harmonic extension of $\phi \in \Lcal^1_{2, \ga}(m)$ is constructed using conformal maps. We therefore leave some of the details out and only indicate the necessary changes due to the extra derivative. 

\begin{proof}
Without loss of generality we may assume $\Om^c$ is bounded. In fact, for any $a\in\Om\backslash\{0\}$, the function $z \mapsto  {1}/({z - a})$ maps $\Om^c$ to a bounded domain which contains the origin.  

We now sketch how to construct a conformal bijection from $\Om^c$ to the horizontal strip $\Pi:=\{\tau + i \nu \, : \, \tau \in \Rbb, \, \, |\nu| < \pi/2\}$. 
First, translate $\Om^c$ so that the origin be in the interior of the bounded component of $\Om$, then apply the complex square root with the branch cut chosen along the negative real axis. The singular point $z(\pm\al_*)$ is thereby mapped to some $\pm iy_0$ on the imaginary axes and the boundary of the resulting domain $\t \Om^c$ satisfies the arc-chord condition and is of class $\Ccal^{1, \lambda'}$ for some $\lambda' \in (0,1)$. By the Riemann mapping theorem, there exists a conformal bijection $\t f :\Dbb \rightarrow \t\Om^c$ which can be extended to a homeomorphism up to the boundary. Let $f: \Dbb \rightarrow \Om^c$ denote $\t f^2$ modulo an appropriate translation. Then we define $F: \Pi \rightarrow \Om^c$ to be $F(w) := f\circ (i\tanh(w/2))$.

This construction does not respect symmetry with respect to the $y$-axis. However, by the uniqueness part of the Riemann mapping theorem, choosing $\t f(0)$ so that $f(0) =  0$ ensures the resulting mapping from $\Dbb$ to the original unbounded (fully symmetric) $\Om^c$ respects both symmetries. In this case, we must have $\t f(\pm i) = \pm iy_0$.

We now give estimates on the derivatives of these conformal maps. The Kellog-Warschawski theorem implies that $\t f \in \Ccal^{1, \lambda'}(\overline{\Dbb})$ for some $\lambda' \in (0,1)$ and that
\be\label{eq:tf'}
|\t f'(\zeta)| \sim 1, \qquad \forall \zeta \in \overline{\Dbb}.
\ee 
In particular, we have  
$$
\frac{|F'(w)|}{|1 - \tanh^2(w/2)|} \sim 1, \qquad \frac{|F(w)|}{|1 - \tanh^2(w/2)|} \sim 1, \qquad \forall w\in \overline{\Pi}.
$$
In order to obtain estimates on the second order derivatives, we use the Schwarz integral formula on the holomorphic map $\log \t f'$. More precisely, we have
\be\label{eq:tf''}
\log \t f'(\zeta) = \log |\t f'(0)| + \frac{i}{2\pi} \int_0^{2\pi}\frac{e^{i\vartheta} + \zeta}{e^{i\vartheta} - \zeta}\, \arg \t f'(e^{i\vartheta})d\vartheta, \quad \forall\zeta \in \Dbb.
\ee 
where $\arg \t f'(e^{i\vartheta}) = \t \theta(\vartheta)- \vartheta - \pi/2$ and $\t \theta(\vartheta)$ is the tangent angle at $\t f(e^{i\vartheta})$ by the Lindel\"{o}f theorem, see e.g. \cite{Pomm}. It is not difficult to see that $\t \theta'$ belongs to $\Lcal_{2, \beta}(\pa\Dbb, \t m)$ with respect to $\t m(\zeta) \sim |\zeta \mp i|$ (when $\zeta$ is  near $\pm i$), hence we can take a complex derivative of the above formula and integrate by parts to obtain
$$
\zeta \frac{\t f''(\zeta)}{\t f'(\zeta)} = \frac{i}{2\pi} \int_0^{2\pi}\frac{e^{i\vartheta} + \zeta}{e^{i\vartheta} - \zeta}\, \t \theta'(\vartheta) d\vartheta, \quad \forall\zeta \in \Dbb.
$$
Thus, the angular limit of $\t f''$ on $\pa \Dbb$ exists and belongs to $\Lcal_{2, \beta}(\pa\Dbb, \t m)$, i.e. we have $\t f\in H^2_\beta(\pa\Dbb, \t m)$. 


Setting $g_\pm(\tau)  := \phi \circ F(\tau \pm i\pi/2)$, it is not difficult to see that 
$$
\phi\in \Lcal^2_{2, \ga + k}(\Ga) \quad \Leftrightarrow \quad g_\pm \in H^2_{-\t \ga + k - 1}\, (\Rbb, \h m), 
$$
where $\h \ga := 1 - (\ga + 1/2)$ and 
$$
\h m(\tau) := \frac{1}{\cosh\tau} \, \sim \, \left|1 - \tanh^2 \left(\frac{w}{2}\right)\right|, \qquad w=\tau + i \nu, \qquad |\nu| \leq \pi/2.
$$
Let $\Phi$ denote the harmonic extension of $g_\pm$ to $\Pi$. Then, the Fourier transform of its normal derivative reads (modulo constant factors in front of each term)  
\be\label{eq:fourier_phi}
\widehat{\pa_\nu\Phi}(\xi) \sim \widehat{g'_e}(\xi) \tanh(\pi\xi/2) \pm \widehat{g'_o}(\xi) \left(\coth (\pi\xi/2) - \frac{2}{\pi\xi}\right) \pm \widehat{g_o}(\xi),
\ee
where we have set $g_e := (g_+ + g_-)/2$ and $g_o := (g_+ - g_-)/2$.

When $|k - \t \ga - 1| < 1$, both $\tanh(\pi \xi/2)$ and $\coth (\pi\xi/2) - \frac{2}{\pi\xi}$ are Fourier multipliers for weighted Lebesgue spaces with respect to the  `exponential' weight $\h m(\tau)^{- \t \ga + k -1}$  (cf. \cite{CEG2} and references therein for proof), hence
$$
\pa_\nu\Phi, \, \pa_\tau\pa_\nu \Phi \in \Lcal_{2, - \t \ga + k -1}(\Rbb, \h m).
$$ 
In particular, when $k = 1$ or $k=2$, the estimate \eqref{eq:estUnormal} follows setting $U:= \Phi\circ F^{-1}$.

When $k=0$, the above condition is not satisfied. However, applying the inverse Fourier transform to \eqref{eq:fourier_phi}, we have (modulo constant factors)
$$
\pa_\nu\Phi(\xi) \sim \int_{\Rbb} g_e'(\tau) \,\frac{1}{\sinh(\xi - \tau)} \, d\tau \pm \int_{\Rbb} g_o'(\tau) \,\left(\coth(\xi - \tau) - \mathrm{sgn}(\xi - \tau)\right) d\tau \pm g_o(\xi),
$$
where for the symmetric part we can write
$$
\int_{-\infty}^\infty g'_e(\tau) \frac{1}{\sinh(\xi - \tau)}\, d\tau = \h m(\xi) \int_{-\infty}^\infty \frac{g'_e(\tau)}{\h m(\tau)} \frac{1}{\sinh(\xi - \tau)}\, d\tau -  \int_{-\infty}^\infty \frac{g'_e(\tau)}{\h m(\tau)} \frac{\h m(\xi) - \h m(\tau)}{\sinh(\xi - \tau)}\, d\tau. 
$$
The first integral belongs to $\Lcal_{2, -(1  + \t \ga)}(\Rbb, \h m)$ by the Fourier multiplier theorem. 

As for the second, we assume without loss of generality that $\xi > 0$ and claim
$$
\int_{-\infty}^\infty g'_e(\tau) k_e(\xi, \tau)\, d\tau = - e^{-\xi}\tanh \xi  \int_{-\infty}^\infty g'_e(\tau)e^{\tau}d\tau + \Lcal_{2, -(1  + \t \ga)}(\Rbb, \h m),
$$
where we have set 
$$
k_e(\xi, \tau) := \frac{1}{\h m(\tau)}\frac{\h m(\xi) - \h m(\tau)}{\sinh(\xi - \tau)}.
$$
Indeed, using the definition of $\h m$, we can write
$$
\aligned
 k_e(\xi, \tau) &= \frac{1}{\cosh\xi} \frac{(e^\tau - e^\xi) + (e^{-\tau} - e^{-\xi})}{e^{\xi - \tau} - e^{-\xi + \tau}} \\
 &=-\frac{1}{\cosh \xi} \frac{\sinh \xi + \sinh \tau}{(e^\xi + e^\tau)(e^{-\xi} + e^{-\tau})} = -\frac{e^{-\xi + \tau}}{(1 + e^{-\xi + \tau})^2} \frac{\sinh \xi + \sinh \tau}{\cosh \xi}.
\endaligned
$$
The claim now follows using the corresponding Hardy inequalities for exponential weights (cf. Theorem \ref{thm:hardy}). In fact, choosing, in the notation of Theorem \ref{thm:hardy}, the weight functions $m_1(\tau) \sim e^{\al\tau}$ and $m_2(\tau) \sim e^{\beta \tau}$ defined on $(0,\infty)$, it is not difficult to see that condition \eqref{eq:hardyI} from Theorem \ref{thm:hardy} holds if $\al\leq \beta, \al < 0$, while condition \eqref{eq:hardyII} holds if $\al\leq \beta, \beta > 0$. We have $\al = \beta = \t \ga + 1$ (since $\h m(\tau) \sim e^{-|\tau|}$). By the above, the kernel $k_e(\xi, \tau)$ satisfies the estimates  
$$
k_e(\xi, \tau) = e^{-\xi + \tau}\tanh \xi + O(e^{-2\xi + 2\tau}), \quad |\tau| < \xi, \qquad \qquad k_e(\xi, \tau) = O(1), \quad |\tau| > \xi.
$$
Since $g'_e(\tau) e^{|\tau|} \in \Lcal_{2, -\t \ga}(\Rbb, \h m)\subseteq L^1(\Rbb)$, the corresponding term contributes a constant to the normal derivative of $U:= \Phi\circ F^{-1}$. 

The anti-symmetric part does not contribute any constants. Proceeding as above, we are led to consider 
$$
\aligned
k_o(\xi, \tau) &:= \frac{\h m(\xi) - \h m(\tau)}{\h m(\tau)}\left(\coth(\xi - \tau) - \mathrm{sgn}(\xi - \tau)\right) \\
&= \frac{1}{2}k_e(\xi, \tau) \left(e^{\xi - \tau} + e^{-\xi + \tau} - \mathrm{sgn}(\xi - \tau)(e^{\xi - \tau} - e^{-\xi + \tau})\right),
\endaligned
$$
which implies
$$
k_o(\xi, \tau) =  O(e^{-2\xi + 2\tau}), \quad |\tau| < \xi, \qquad \qquad k_o(\xi, \tau) = O(1), \quad |\tau| > \xi
$$ 
and the corresponding integrals belong to $\Lcal_{2, -(1  + \t \ga)}(\Rbb, \h m)$ by Hardy inequalities as above.

In particular, when $\phi\in\Lcal^2_{2,\ga}(m)$ is odd with respect to the $y$-axis, we have $g_e = 0$ and therefore $\pa_\nu \Phi\in \Lcal_{2, -(\t \ga+1)}(\h m)$ as required. The same is clearly true for $\pa_\tau\pa_\nu \Phi$ as well. 

Finally, let $\phi\in H^2_{\ga,0}(m)$ be odd with respect to both axis. Then, locally we can write $\phi = \mp c_\pm \chi(z)\Re(z) + \Lcal^2_{2,\ga}(m)$, where $c_\pm := \phi'(\pm\al_*)$ and $\chi$ is a smooth, symmetric cut-off identically equal to $1$ on an open ball centered at the origin and identically equal to zero outside a larger ball. By the previous step, it is enough to consider $\phi = \mp c_\pm\chi(z)\Re(z)$. We extend $\phi$ to $\Om^c$ in an obvious manner. Its harmonic extension is then given by $U := V  - c_\pm\mathrm{sign}(z_2)\chi(z)\Re(z)$, with $V$ a solution of the Poisson equation
$$
\Delta V = G, \qquad G = \Delta(c_\pm \mathrm{sign}(z_2)\chi(z)\Re (z)), \qquad V\big|_\Ga  = 0  
$$ 
in $\Om^c$. If $\pa_n V\in \Lcal^1_{2,\ga}(m)$, then $\pa_n U = c_\pm\mathrm{sign}(z_2)\chi(z)  z_{2s} + \Lcal^1_{2,\ga}(m)$. As above, we map $\Om^c$ to the horizontal strip $\Pi$. Then $\t V = V \circ F$ satisfies the Poisson equation on $\Pi$ with the right hand side $\t G = G\circ F |F'|$ where $\t G$ is a bounded, smooth, compactly supported function on $\Pi$. The solution $\Phi$ is given by the Green's function on the strip, with the normal derivative at e.g. $w = \xi + i\pi/2$ given by
$$
\pa_n \Phi(w) = \frac{1}{\pi} \int_\Pi \t G(w')\coth(w-w')dw' = \frac{1}{2\pi}\int_{-\pi/2}^{\pi/2}\int_{-\infty}^\infty \t G(\tau, \nu)\frac{\sinh 2(\xi-\tau)}{\sinh^2(\xi-\tau) + \cos^2\nu} \, d\tau d\nu.
$$  
We need to show $\pa_n \Phi\in \Lcal_{2,  -(\t \ga + 1)}(\Rbb, \h m)$. Since $\t G$ is compactly supported, it is enough to consider $\xi$ large so that $\xi\gg\tau$. Then, we can write   
$$
\frac{1}{2}\frac{\sinh 2(\xi-\tau)}{\sinh^2(\xi-\tau) + \cos^2\nu} = \coth(\xi - \tau)(1 + O(e^{-2\xi + 2\tau})), \quad \coth(\xi - \tau) = 1 + O(e^{-2\xi + 2\tau})
$$
However, since $\t G$ is odd with respect to $\nu \leftrightarrow -\nu$, we have $\int_\Pi \t G = 0$ and the claim follows.

\end{proof}

\subsection{Domain with a corner singularity}

Next we consider the case when our domain exhibits an angled crest:

\begin{theorem}\label{thm:inverseOpII} 
	Let $\max\left\{3 - \frac{\pi}{2\nu}, 3 - 2\frac{\pi}{2\overline{\nu}}\right\} < \beta + \frac{1}{2} < 1$ with $\beta + \frac{1}{2} \neq 2-\frac{\pi}{2\overline{\nu}}$, where $2\overline{\nu}:= \pi - 2\nu$ and $2\nu\in (0,\frac{\pi}{2})$. Let $c\in \Rbb$ and let $\psi = c\pa_s \Re(z^2) +  \Lcal_{2,\beta-2}(m)$ odd with respect to both axis. Then there exists a unique solution $\phi\in \Lcal_{2,\beta-2}(m)$, odd with respect to both axes,  of 
	$$	
	(I + S)\phi = \psi.
	$$
\end{theorem}

Note that by symmetry, we have $\int \psi \,ds = 0$ on each component of the boundary.



\subsubsection{Construction of harmonic functions on $\Om$}

\begin{proposition}\label{prop:inner-corner} Let $2\nu\in (0,\frac{\pi}{2})$. If $\max\left\{3 - \frac{\pi}{2\nu}, 0\right\} < \beta + \frac{1}{2} < 1$ and $\phi = c\Re(z^2) + \Lcal^1_{2,  \beta-2}(m)$ for some $c\in \Rbb$, there exists a harmonic extension $U$ of $\phi$ to $\Om_\pm := \Om\cap \Rbb_\pm$ such that  
$$
\pa_n U = -c\pa_s\Im(z^2) + \Lcal_{2,  \beta-2}(m).
$$
\end{proposition}

\begin{proof}
	
Assume first $\phi = c \Re (z^2)$. Since $\Om$ is bounded, we can take $U = c\Re (z^2)$, which contributes $c\Re(z_s^\perp z) = -c\pa_s \Im(z^2)$ to the normal derivative as required. 

Let now $\phi \in \Lcal^1_{2,\beta - 2}(m)$. We can write $\phi = \chi\phi + (1-\chi)\phi =: \phi_1 + \phi_2$, where $\chi$ is a smooth cut-off identically equal to one in the neighborhood of the singular point. We construct the harmonic extensions $U$ of $\phi_1$ respectively $\phi_2$ separately and we claim $\pa_n U\in \Lcal_{2,\beta - 2}(m)$. 

Let $\Om$ temporarily denote the right component $\Om \cap \Rbb_+^2$. It is not difficult to see that $z \rightarrow z^{\frac{\pi}{2\nu}}$ maps $\Om$ to a bounded domain $\t \Om$ with $\Ccal^{1,\lambda'}$-boundary $\t \Ga$ for some $\lambda'\in (0,1)$. In fact, the corresponding tangent angle $\t \theta$ is Holder continuous and $\t \Om$ has a vertical tangent at the origin. By the Riemann mapping theorem there exists a conformal bijection $\t f:\Dbb\rightarrow \t \Om$ such that $\t f(-1) = 0$ and, by the  Kellog-Warschawski theorem, we must have $\t f \in \Ccal^{1, \lambda'}(\overline{\Dbb})$ with $\t f'$ satisfying a \eqref{eq:tf'}-type estimate. Let $f:\Dbb\rightarrow \Om$ denote the composition of these two mappings. Then 
$$
|f(\zeta)| \sim |\zeta + 1|^{\frac{2\nu}{\pi}}, \qquad |f'(\zeta)| \sim |\zeta + 1|^{\frac{2\nu}{\pi}-1}, \qquad \forall \zeta \in  \overline{\Dbb}\cap B_\de(-1)
$$
for some small $\de > 0$. First assume that we have found a harmonic extension $\Phi$ of $\phi_2\circ f$ to $\Dbb$, such that $U := \Phi \circ f^{-1}\in \Lcal_{2,\beta - 2}(\Om)$. It is not difficult to see that this implies 
$$
\pa_n U\in\Lcal_{2,\beta - 2}(m) \quad \Leftrightarrow \quad \pa_n \Phi\in \Lcal_{2, \frac{2\nu}{\pi} \beta' + \frac{1}{2}}(\t m),
$$ 
where $\beta' = (\beta + 1/2) - 3$ and $\t m (\zeta) \sim |\zeta + 1|$. However, by assumption $\phi_2 \circ f\in H^1(\Tbb)$ and we can define $\Phi$ on $\Dbb$ through the Poisson formula. Its normal derivative is then given by
$$
\pa_n \Phi(e^{i\vartheta}) = \frac{1}{2\pi}p.v.\int_{-\pi}^\pi \frac{d}{d\varkappa}\phi_2(f(e^{i\varkappa}))\cot\left(\frac{\vartheta - \varkappa}{2}\right)d\varkappa \in L^2(\Tbb).
$$ 
However, $\pa_n \Phi$ is actually smooth in the neighborhood of $\zeta = 1$, since $\phi_2 \circ f$ identically vanishes there. In particular, it is bounded near $\zeta = 1$ and we actually have $\pa_n \Phi \in \Lcal_{2, \frac{2\nu}{\pi} \beta' + \frac{1}{2}}(\Tbb, \t m)$ provided $0 < \left(\frac{2\nu}{\pi} \beta' + \frac{1}{2}\right) + \frac{1}{2}$. Since
$$
\left(\frac{2\nu}{\pi} \beta' + \frac{1}{2}\right) + \frac{1}{2} > 0 \quad \Leftrightarrow \quad \beta + \frac{1}{2} > 3 - \frac{\pi}{2\nu}
$$  
(where $2\nu < \frac{\pi}{2}$ is equivalent to $3 - \frac{\pi}{2\nu} < 1$), the claim follows.

It remains to consider $\phi_1$. Let $F$ be the conformal bijection from the upper half-plane $\Hbb$ to $\Om$ such that $F(0) = 0$ ($F$ is just the composition of $f$ with a suitable M{\"o}bius transformation). Then, near the origin 
$$
|F(\t z)| \sim |\t z|^{\frac{2\nu}{\pi}}, \qquad |F'(\t z)| \sim |\t z|^{\frac{2\nu}{\pi}-1}, \qquad \forall \t z \in  \overline{\Hbb}, \quad |\t z|\lesssim 1,
$$
while far away
$$
|F(\t z)| \sim |\t z|^{-1}, \qquad |F'(\t z)| \sim |\t z|^{-2}, \qquad \forall \t z \in  \overline{\Hbb}, \quad |\t z|>> 1
$$
As above, we have $\phi_1\circ F\in \Lcal^1_{2, \frac{2\nu}{\pi} \beta' + \frac{1}{2}}(\Rbb, \t m)$, where $\t m(\t x)\sim |\t x|$ in a neighborhood of the origin and $\t m(\t x)\sim |\t x|^{-1}$ otherwise. Let $\Phi$ denote the harmonic extension of $\phi_1\circ F$ given by the Poisson kernel for the upper half-plane, then 
$$
\pa_n \Phi(\t x) = \frac{1}{\pi}p.v.\int_{-\infty}^\infty \t\phi'(t) \,\frac{1}{\t x - t}\, dt \in \Lcal_{2, \frac{2\nu}{\pi} \beta' + \frac{1}{2}}(\Rbb, \t m)
$$
as required, since by assumption $\frac{2\nu}{\pi} \beta' + \frac{1}{2}$ satisfies the M{\"u}ckenhaupt condition.

\end{proof}

\subsubsection{Construction of harmonic functions on $\Om^c$}

\begin{proposition}\label{prop:outer-corner} Let $2\overline{\nu} \in [\frac{\pi}{2} , \pi)$ and let $0 < \beta + \frac{1}{2} < 1$. If $3 - 2\frac{\pi}{2\overline{\nu}} <\beta + \frac{1}{2}$ and $\phi = c\Im(z^2) + \Lcal^1_{2,  \beta-2}(m)$ is odd with respect to $y$-axis, the normal derivative of its harmonic extension $U$ to the exterior domain $\Om^c$ reads  
	$$
	\pa_n U  = c\Re(z^2) +  \Lcal_{2,  \beta-2}(m).
	$$ 

\end{proposition}

\begin{proof}
	By the Riemann mapping theorem there exists a conformal bijection  $f:\Om^c\rightarrow \t\Om^c$, with $\t \Om^c$ as in Proposition \ref{prop:dirichletExt} (i.e. $\t \Om$ has two connected components, each with an outward cusp and a common tip situated at the origin). By the Kellog-Warschawski theorem (and assuming that $f(0) = 0$), we must have
	$$
	f(z) \sim (\mp iz)^{\frac{\pi}{2\overline{\nu}}}, \qquad z\in \Om^c \cap B_\de(0),  \qquad \Im z \gtrless 0
	$$    
	for some small $\de > 0$. Parametrizing the boundary of $\pa\t \Om$ as $\t z(\al) = f(z(\al))$ where $\al \in [-\pi, \pi]$, it is not difficult to check that the corresponding tangent angle is Holder continuous and satisfies $\t \theta(\al_*) = 0$ resp. $\theta(-\al_*) = \pi$.
	
   Assume first that $\phi\in \Lcal^1_{2,  \beta-2}(m)$. As before, we have 
   $$
   \phi\in \Lcal^1_{2,  \beta-2}(m)\quad \Leftrightarrow \quad \phi\circ f^{-1}\in \Lcal^1_{2,\frac{2\overline{\nu}}{\pi}\beta' + \frac{1}{2}}(m), 
   $$
   where $\beta' = (\beta + 1/2)-3$ and by assumption $0 < \left(\frac{2\overline{\nu}}{\pi}\beta' + \frac{3}{2}\right) + \frac{1}{2} < 1$. In fact, we have
   $$
   0 < \left(\frac{2\overline{\nu}}{\pi}\beta' + \frac{3}{2}\right) + \frac{1}{2} < 1 \quad \Leftrightarrow \quad 3-2\frac{\pi}{2\overline{\nu}} < \beta + \frac{1}{2} < 3-\frac{\pi}{2\overline{\nu}},
   $$
   where by assumption $2\overline{\nu} \geq \pi/2$ and therefore $3 - \frac{\pi}{2\overline{\nu}} \geq 1$. In particular, by Proposition \ref{prop:dirichletExt}, there exists a harmonic extension $\t U$ of $\phi\circ f^{-1}$ to $\t \Om^c$ such that $\pa_n \t U\in \Lcal_{2, \frac{2\overline{\nu}}{\pi}\beta' + \frac{1}{2}}(m)$. The claim now follows going back to the original domain.

   Let now $\phi = c \Im (z^2)$ and let $\chi$ be a smooth symmetric cut-off identically equal to $1$ on a sufficiently large open ball centered at the origin (so that $\Omega$ is a proper subset of this ball) and identically equal to zero outside of some larger ball. The obvious extension $G := c\chi(z)\Im (z^2)$ to the exterior of $\Om$ leads to the Poisson problem
   $$
 	\Delta V = -\Delta G, \qquad V\big|_\Ga  = 0,  
   $$ 
  with $U = V + G$ the required harmonic extension of $\phi$. However, by the last part of the proof of Proposition \ref{prop:dirichletExt}, we know that $V\in \Lcal^1_{2,\beta - 2}(m)$, hence $\pa_n U = c\pa_s \Re (z^2) + \Lcal_{2,\beta-2}(m)$ as required. 

\end{proof}

\section*{Acknowledgements}

This project has received funding from the European Research Council (ERC) under the European Union's Horizon 2020 research and innovation program through the grant agreements 788250 (D.C.) and 633152 (A.E.). This work is supported in part by the ICMAT-Severo Ochoa grant CEX2019-000904-S.


\begin{appendices}

\section{Auxiliary estimates}

In the following sections, we will  estimate
\be\label{eq:F}
F(\phi)(x, u) := \frac{\phi(x)- \phi(u)}{x - u}.
\ee
Using the Taylor development of $\phi$, we can write  
\be\label{eq:FTaylor}
\pa_x^k F(\phi)(x, u) = \frac{1}{(x - u)^{k+1}}\int_x^{u}\pa_{\tau}^{k+1} \phi(\tau)\frac{(u - \tau)^k}{k!}d\tau
\ee
provided we control a sufficient number of derivatives of $\phi$. Finally, note that  
\be\label{eq:Fx-Fu}
\pa_u F(\phi) + \pa_x F(\phi) = F(\phi').
\ee

\subsection{The kernel on the singular set}

Here, we study the additional 'singular' kernel which appears in the Birkhoff-Rott integral as a consequence of the failure of the arc-chord condition at the cusp tip. More precisely, we isolate the most singular part and give estimates on the error  terms. Then we describe the variable change taking the cusp tip to infinity (separately for each cusp) and give estimates on the remaining kernel and its derivatives. 

\begin{lemma}\label{lem:errbasic}
Let $q_* = u - i \ka(x)$. For $k\in \Nbb$, we have the following recursive formula 
\be\label{eq:recursive}
\frac{z_+'}{z_+ - q_-} - \frac{1}{z_+ - q_*} = K_0(x,u)\sum_{j=0}^{k-1} K(x,u)^j + K(x,u)^k\bigg(\frac{z_+'}{z_+ - q_-} - \frac{1}{z_+ - q_*}\bigg),
\ee
where $\Delta\ka(x,u) := \ka(x)-\ka(u)$ and   
$$
K_0(x,u) : = \frac{i\ka'(x)}{z_+ - q_*} + \frac{i\Delta \ka(x, u)}{(z_+ - q_*)^2} \qquad K(x,u) := \frac{i\Delta \ka(x,u)}{z_+ - q_*}.
$$
For $x\in I_\de$, $u\in I_c(x)$, we have the estimates 
\be\label{eq:rKk}
K(x,u)^k\bigg(\frac{z_+'}{z_+ - q_-} - \frac{1}{z_+ - q_*}\bigg) = O\big(m(x)^{k\mu - 1}\big) 
\ee
respectively,
\be\label{eq:rK0}
\big|\t K_0(x,u) K(x,u)^j\big| \, \lesssim \, m(x)^{(j + 1)\mu} \Big(\frac{1}{m(x)} + \frac{\rho(x)}{(x -u)^2 + \rho(x)^2} \Big), \quad 0 \leq j\leq k-1,
\ee
where  
$$
\t K_0(x, u) := K_0(x,u) - \frac{i\rho'(x)}{z_+ - q_*}.
$$
\end{lemma}
\begin{proof}
To see that \eqref{eq:recursive} holds, we have
$$
\aligned
\frac{z_+'}{z_+ - q_-} - \frac{1}{z_+ - q_*} &= \frac{i\ka'(x)}{z_+ - q_*} + \frac{z_+'}{z_+ - q_-}\frac{i\Delta\ka(x, u)}{z_+ - q_*}\\
&=\frac{i\ka'(x)}{z_+ - q_*} + \frac{i\Delta\ka(x,u)}{(z_+ - q_*)^2} + \frac{i\Delta\ka(x, u)}{z_+ - q_*}\bigg(\frac{z_+'}{z_+ - q_-} - \frac{1}{z_+ - q_*}\bigg)
\endaligned
$$
which proves the claim for $k=1$. Formula \eqref{eq:recursive} now follows by induction. To see that \eqref{eq:rKk} holds, let e.g. $x>0$ and $u\in I_c(x)$. Then $x\sim u$ and   
$$
\ka'(x) = O(x^\mu), \qquad \frac{\Delta \ka(x,u)}{x - u} = O(x^\mu), \qquad |K(x,u)|^j = O(x^{j\mu}).
$$
The first line in the above calculation then implies
$$
\bigg|\frac{z_+'}{z_+ - q_-} - \frac{1}{z_+ - q_*}\bigg|\,\lesssim\, \frac{x^\mu}{|z_+ - q_*|} \,\lesssim\, \frac{1}{x}
$$
(using that $\rho(x) \sim x^{\mu + 1}$) and the claim follows. It remains to show \eqref{eq:rK0}. However, note that  
$$
K_0(x,u) = \frac{i\rho'(x)}{z_+ - q_*} + \ka'(x)\frac{\rho(x)}{(z_+ - q_*)^2} + i \pa_x F(\ka)(x,u) \frac{(x- u)^2}{(z_+ - q_*)^2},
$$
then use $\pa_x F(\ka)(x,u) = O(x^{\mu - 1})$.
\end{proof}

\begin{lemma}\label{lem:errorrk}
Let $z_*:= x + i\ka(u)$. Then, we can write 
\be\label{eq:expansionk}
\frac{z_+'}{z_+ - q_-} =  \frac{1}{z_* - q_-} + r(x,u), 
\ee
where the remainder $r(x,u)$ satisfies 
\be\label{est:r}
|r(x,u)| \, \lesssim \, \frac{1}{m(x)}, \qquad  |\pa_x r(x,u)| \, \lesssim \, \frac{1}{m(x)}\bigg(\frac{1}{m(x)} + \frac{\rho(u)}{(x- u)^2 + \rho(u)^2}\bigg),
\ee
when $u\in I_c(x)$ uniformly for $x\in I_\de$.
\end{lemma}

\begin{proof}
We can proceed as in the derivation of formula \eqref{eq:recursive} to obtain  
$$
\frac{z_+'}{z_+ - q_-} = \frac{1}{z_* - q_-} + K_0(x,u) + K(x,u)\bigg(\frac{z_+'}{z_+ - q_-} - \frac{1}{z_* - q_-}\bigg),
$$
where now 
$$
K(x,u) = -i F(\ka) \frac{x-u}{z_* - q_-}, \qquad K_0(x,u) = \frac{i\ka'(x)}{z_* - q_-} - iF(\ka)\frac{x-u}{(z_* - q_-)^2}
$$
(cf. \eqref{eq:F} for the definition of $F(\ka)$). W.l.o.g. let $x>0$ and $u\in I_c(x)$. Then $x\sim u$ and
$$
F(\ka)(x,u) = O(x^\mu) , \qquad \pa_x F(\ka)(x,u) = O(x^{\mu-1}).
$$
A short calculation then yields
$$
|\pa_x^m K(x,u)|\, \lesssim\,  \frac{1}{x}\, |x-u|^{1-m}, \qquad \bigg|\pa_x^m\bigg(\frac{z_+'}{z_+ - q_-} - \frac{1}{z_* - q_-}\bigg)\bigg| \,\lesssim\, \frac{1}{x}\frac{1}{|x-u|^m},
$$
for $m=0,1$. In particular, 
$$
\pa_x K(x,u)\left(\frac{z_+'}{z_+ - q_-} - \frac{1}{z_* - q_-}\right) + K(x,u) \pa_x \bigg(\frac{z_+'}{z_+ - q_-} - \frac{1}{z_* - q_-}\bigg) =  O\left( \frac{1}{x^{2}}\right).
$$

It remains to consider $K_0(x,u)$. Since 
$$
\ka'(x) = \ka'(u) + F(\ka')(x-u), \qquad F(\ka)= \ka'(u) + \pa_u F(\ka) (x - u)  
$$
we can write
$$
K_0(x,u) = \frac{i\ka'(u)}{z_* - q_-} - i\ka'(u) \frac{x-u}{(z_* - q_-)^2} + iF(\ka')\frac{x-u}{z_* - q_-} - i\pa_u F(\ka)\frac{(x-u)^2}{(z_* - q_-)^2}
$$
and we claim that 
$$
|\pa_x K_0(x,u)| \, \lesssim \, \frac{1}{x}\bigg(\frac{1}{x} + \frac{\rho(u)}{(x- u)^2 + \rho(u)^2}\bigg).
$$
Indeed, the first two terms combined give 
$$
s_1(x,u) := -\ka'(u)\frac{\rho(u)}{(z_* - q_-)^2}.
$$
Its derivative clearly satisfies the desired estimate. As for the remaining two terms note that   
$$
\pa_x \Big(\frac{x-u}{z_* - q_-}\Big) = \frac{i\rho(u)}{(z_* - q_-)^2}, \qquad |\pa_x F(\ka')(x,u) | \, \lesssim \, x^{\lambda - 2}, \qquad |\pa^2_x F(\ka)(x,u) | \, \lesssim \, x^{\lambda - 2}
$$
where we have set $\lambda := 1-(\beta + 1/2)$ (assumption $\ka'' \in H^2_{\beta + 2}(I_\de)$ implies $\ka'''=O(x^{\lambda - 2})$ by Lemma \ref{lem:sobolev}). Using this and \eqref{eq:Fx-Fu} the claim follows.
\end{proof}

\subsubsection*{The  change of variables $\tau\mapsto h(\tau)$}

We implicitly define  
\be\label{eq:hb}
h_+^{-1}(u) := \int^{2\de}_{u} \frac{d\nu}{\rho(\nu)}, \quad u\in (0,2\de), \qquad h_-^{-1}(u) := \int_{-2\de}^{u} \frac{d\nu}{\rho(\nu)}, \quad u\in (-2\de, 0).
\ee
By making $\de$ smaller if necessary, we may assume $\rho(\nu)$ is strictly monotonically increasing on $(0, 2\de)$ respectively decreasing on $(-2\de, 0)$ and therefore $h_{\pm}^{-1}$ can be inverted on their domains of definition. For simplicity we drop the subscript $\pm$ and consider $h^{-1}\equiv h^{-1}_+$ only. Its inverse mapping $h$ satisfies
$$
u=h(\tau), \quad h'(\tau) = -\rho(h(\tau)),
$$
hence using $\rho(u) \sim u^{\mu + 1}$, it is not difficult to see that $h^{-1}: (0, 2\de) \rightarrow (0, \infty)$ and
$$
h^{-1}(u) \sim u^{-\mu} - (2\de)^{-\mu}, \quad u \in (0, 2\de) \qquad \Leftrightarrow \qquad h(\tau) \sim (1 + \tau)^{-1/\mu}, \quad \tau\in \Rbb_+. 
$$
This asymptotic expansion can be differentiated three times (cf. estimate \eqref{assump:rho''}), i.e. setting  
$$
\t m(\tau) := 1 + \tau, \quad \tau\in \Rbb_+,
$$
we can write
\be\label{eq:hbders}
|\pa_\tau^jh(\tau)| \sim \t m(\tau)^{-j - 1/\mu}, \quad j=0,1 \qquad |\pa_\tau^j h(\tau)| \, \lesssim \, \t m(\tau)^{-j - 1/\mu} \quad j=2,3.
\ee
Note that we also have lower bounds on $\pa^j_\tau h$ when $j = 0,1$. We will typically apply $h$ to the Hilbert transform and the ` singular' kernel from the previous section in the region
$$
x \in (0, \de), \quad  u\in I_c(x) \quad \quad \Leftrightarrow \quad \quad \xi \in \wt I_\de := h^{-1}((0, \de)), \quad  \tau \in \wt I_c(\xi) := h^{-1}(I_c(x)) \quad x = h(\xi). 
$$ 
It is not difficult to see that there exist small $\varepsilon_\pm > 0$ such that
\be\label{eq:veps12h}
\{\tau \,:\, |\tau - \xi|<\varepsilon_- \xi \} \subseteq \wt I_c(\xi)\subseteq  \{\tau \,:\, |\tau - \xi|< \varepsilon_+ \xi \}, \quad \forall \xi \in \wt I_\de. 
\ee
Let us e.g. show the existence of $\varepsilon_+$. Setting $\xi_\pm:= h^{-1}((1\mp \varepsilon)x)$, we have, by definition
$$
\xi_+ - \xi_- = \int^{(1+ \varepsilon)x}_{(1-\varepsilon)x} \frac{d\nu}{\rho(\nu)} \, \sim \, \varepsilon h(\xi)^{-\mu}
$$
and the claim follows using \eqref{eq:hbders} (use that $\xi > h^{-1}(\de) = O(\de^{-\mu})$, when $\xi \in \wt I_\de$). Note that $ \varepsilon_+$ can be made arbitrarily small by making $\varepsilon$ in the definition of $I_c(x)$ smaller. We then control $\varepsilon^{-1}$ by some polynomial depending on the norm of $\rho$. The existence of $\varepsilon_-$ follows similarly.


\begin{lemma}\label{lem:errVCHilb}
Let $\epsilon>0$. Then, we can write  
 $$
\frac{h'(\tau)}{h(\xi) - h(\tau)} = \frac{1}{\xi - \tau} + r(\xi, \tau),
 $$ 
where the remainder satisfies
\be\label{eq:err20}
|r(\xi, \tau)| \, \lesssim \, \frac{1}{\t m(\tau)}, \qquad |\pa_\xi r(\xi, \tau)| \, \lesssim \, \frac{1}{\t m(\tau)^2}.
\ee
for all $|\tau - \xi| < \epsilon\, \xi$.
\end{lemma}

\begin{proof}
Indeed, we can write 
\be\label{aux:vc1}
\frac{h'(\tau)}{h(\xi) - h(\tau)}- \frac{1}{\xi - \tau} = \frac{1}{\xi - \tau} \bigg( \frac{h'(\tau)}{F(h)} - 1\bigg) = -\frac{\pa_\tau F(h)}{F(h)}.
\ee
The claim follows using \eqref{eq:FTaylor} and \eqref{eq:hbders} to conclude 
$$
|F(h)| \sim \t m(\tau)^{-1/\mu}, \qquad  |\pa_\tau^j F(h)| \, \lesssim \, \t m(\tau)^{-(j + 1) - 1/\mu}, \qquad j = 1,2, 
$$
whenever $\tau \sim \xi$.
\end{proof}


\begin{lemma}\label{lem:kernelDecomposition}
Let $\epsilon>0$. Then, we can write  
$$
\frac{h'(\tau)}{(h(\xi) - h(\tau)) - i h'(\tau)} = \frac{1}{(\xi - \tau) - i}+ r(\xi, \tau), 
$$
where the remainder $r(\xi, \tau)$ satisfies
\be\label{eq:rerrvcI}
|r(\xi, \tau)| \, \lesssim \, \frac{1}{\t m(\tau)}, \qquad |\pa_\xi r(\xi, \tau)| \, \lesssim \, \frac{1}{\t m(\tau)}\Big(\frac{1}{\t m(\tau)} +  \frac{1}{(\xi - \tau)^2 + 1}\Big)
\ee
for all $|\tau - \xi| < \epsilon\, \xi$.
\end{lemma}

\begin{proof}
Let 
$$
k(\xi, \tau) :=  \frac{h'(\tau)}{(h(\xi) - h(\tau)) - i h'(\tau)} = \frac{h'(\tau)}{F(h)} \frac{1}{\big((\xi - \tau) - i\big) + i\big(1 - \frac{h'(\tau)}{F(h)}\big)}.
$$
Then, it is not difficult to see that  
$$
r(\xi, \tau) = \Big(\frac{h'(\tau)}{F(h)} - 1\Big)\frac{1}{(\xi - \tau) - i}  -  \frac{i\big(1 - \frac{h'(\tau)}{F(h)}\big)}{(\xi - \tau) - i} \, k(\xi, \tau) =: r_1(\xi, \tau) + r_2(\xi, \tau). 
$$
Since $\xi \sim \tau$, we have by assumption 
$$
\frac{h'(\tau)}{F(h)} - 1 = O \Big(\frac{\xi - \tau}{\t m(\tau)}\Big)
$$
cf. \eqref{aux:vc1} and the zero order estimate is straightforward. To show the claim for $\pa_\xi r_1(\xi, \tau)$, we write
$$
\aligned
\pa_\xi r_1(\xi, \tau) &=  \pa_\xi\Big(\frac{h'(\tau)}{F(h)}\Big) \frac{1}{(\xi - \tau) - i} -  \Big(\frac{h'(\tau)}{F(h)} - 1\Big)\frac{1}{((\xi - \tau) - i)^2}\\
&=\pa_\xi\Big(\frac{h'(\tau)}{F(h)}\Big) \frac{1}{(\xi - \tau) - i} - \pa_\xi\Big(\frac{h'(\tau)}{F(h)}\Big)\Big|_{\xi = \tau}\frac{(\xi - \tau)}{((\xi - \tau) - i)^2} + O\Big(\frac{1}{\t m(\tau)^2}\Big)\\
&= \pa_\xi\Big(\frac{h'(\tau)}{F(h)}\Big) \frac{-i}{((\xi - \tau) - i)^2} + O\Big(\frac{1}{\t m(\tau)^2}\Big),
\endaligned
$$
where in the second line, we have used the Taylor development of the function 
$$
\xi \rightarrow g(\xi; \tau) := \frac{h'(\tau)}{F(h)(\xi, \tau)} 
$$
around $\xi = \tau$ up to order $2$, i.e. 
$$
g(\xi; \tau) = 1 + \pa_\xi\Big(\frac{h'(\tau)}{F(h)}\Big)\Big|_{\xi = \tau}(\xi - \tau) + O\Big(\frac{(\xi - \tau)^2}{\t m(\tau)^2}\Big), \qquad  \pa_\xi\Big(\frac{h'(\tau)}{F(h)}\Big)\Big|_{\xi = \tau} = -\frac{1}{2}\frac{h''(\tau)}{h'(\tau)}.
$$ 
This relation can be differentiated (using up to three derivatives of $h$), which is used to pass to the third line. Since we also have 
$$
\Big|\pa_\xi \Big(\frac{h'(\tau)}{F(h)}\Big)\Big| \, \lesssim \, \frac{1}{\t m(\tau)} ,
$$
we conclude $\pa_\xi r_1(\xi, \tau)$ satisfies estimate \eqref{eq:rerrvcI}. In order to show a similar estimate for $r_2(\xi, \tau)$, note that we can write
$$
r_2(\xi, \tau) = -i\frac{h'(\tau)}{F(h)}\bigg(1 - \frac{h'(\tau)}{F(h)}\bigg)\frac{1}{((\xi - \tau) - i)^{2}} -  \frac{\big(1 - \frac{h'(\tau)}{F(h)}\big)^2}{((\xi - \tau) - i)^2} \, k(\xi, \tau).
$$
Taking the derivative and using  
$$
|k(\xi, \tau)| = O(1), \qquad |(\xi - \tau)k(\xi, \tau)|=O(1), \qquad |\pa_\xi^j k(\xi, \tau)| \, \lesssim \, |k(\xi, \tau)|^2,
$$
when $\xi \sim \tau$ it is not difficult to see the claim follows.
\end{proof}

 
\begin{lemma}\label{lem:errvc}
Let $\epsilon>0$. Then, we can write
$$
\frac{h'(\tau)}{(h(\xi) - h(\tau)) - ih'(\xi)} - \frac{h'(\tau)}{h(\xi) - h(\tau)} = \frac{1}{(\xi - \tau) - i} - \frac{1}{\xi - \tau} + r(\xi, \tau),
$$
where the remainder $r(\xi, \tau)$ satisfies the estimate  
\be\label{eq:errvcreg}
|r(\xi, \tau)| \, \lesssim \, \frac{1}{\t m(\tau)}\Big(\frac{1}{\t m(\tau)} +  \frac{1}{1 + (\xi - \tau)^2}\Big)
\ee
for all $|\tau - \xi| < \epsilon\, \xi$.
\end{lemma}
\begin{proof}
Writing  
$$
\frac{h'(\tau)}{(h(\xi) - h(\tau))  - i h'(\xi)} = \frac{\frac{h'(\tau)}{F(h)}}{(\xi - \tau)  - i \frac{h'(\tau)}{F(h)}}, \qquad \frac{h'(\tau)}{h(\xi) - h(\tau)} = \frac{\frac{h'(\tau)}{F(h)}}{\xi - \tau},
$$
then taking the difference, we have 
$$
r(\xi, \tau)= \frac{1}{\xi - \tau}\frac{i\big(\frac{h'(\xi)h'(\tau)}{F(h)^2}  - 1\big)}{(\xi - \tau) - i\frac{h'(\xi)}{F(h)}} +\frac{\frac{1}{\xi - \tau}\big(\frac{h'(\xi)}{F(h)} - 1\big)}{((\xi - \tau) - i)\big((\xi - \tau) - i\frac{h'(\xi)}{F(h)}\big)} =: r_1(\xi, \tau) + r_2(\xi, \tau).
$$
The kernel $r_1(\xi, \tau)$ satisfies \eqref{eq:errvcreg}, since
$$
\frac{h'(\xi)h'(\tau)}{F(h)^2}  - 1 = O\Big(\frac{(\xi - \tau)^2}{\t m(\tau)^2}\Big),
$$
when $\xi\sim \tau$. Indeed, we have 
$$
\aligned
\Big(\frac{h'(\xi)h'(\tau)}{F(h)^2}  - 1\Big) &=  \Big(\frac{h'(\xi)}{F(h)} - 1\Big) + \Big(\frac{h'(\tau)}{F(h)} - 1 \Big) + \bigg(\frac{h'(\xi)}{F(h)} -1\bigg) \bigg(\frac{h'(\tau)}{F(h)}-1\bigg)\\
& = \Big(\frac{h''(\xi)}{h'(\xi)} - \frac{h''(\tau)}{h'(\tau)} \Big)\frac{\xi - \tau}{2} + O\Big(\frac{(\xi - \tau)^2)}{\t m(\tau)^2}\Big),
\endaligned
$$
since 
$$
\aligned
\frac{h'(\xi)}{F(h)}  &= 1 +  \frac{h''(\xi)}{h'(\xi)} \frac{(\xi - \tau)}{2} + O\Big(\frac{(\xi - \tau)^2)}{\t m(\tau)^2}\Big), \\
\frac{h'(\tau)}{F(h)} &= 1 -  \frac{h''(\tau)}{h'(\tau)} \frac{(\xi - \tau)}{2} + O\Big(\frac{(\xi - \tau)^2)}{\t m(\tau)^2}\Big),
\endaligned
$$
cf. the proof of Lemma \ref{lem:kernelDecomposition}. In particular, the claim follows since also
$$
\frac{1}{\xi - \tau}\Big(\frac{h''(\xi)}{h'(\xi)} - \frac{h''(\tau)}{h'(\tau)} \Big) = O\Big(\frac{1}{\t m(\tau)^2}\Big).
$$
On the other hand, we can rewrite $r_2(\xi, \tau)$ as  
$$
r_2(\xi, \tau) =\frac{\frac{1}{\xi - \tau}\big(\frac{h'(\xi)}{F(h)} - 1\big)}{((\xi - \tau) - i)^2} - \frac{\frac{i}{\xi - \tau}\big(\frac{h'(\xi)}{F(h)} - 1\big)^2}{((\xi - \tau) - i)^2\big((\xi - \tau) - i\frac{h'(\xi)}{F(h)}\big)}
$$
and the claim is then straightforward.
\end{proof}

\begin{subsection}{Commutator estimates}

\begin{lemma}\label{lem:trHilbert}
 Let $\phi\in H^1_{\ga + 1}(I_{2\de})$ where $0 < \ga + 1/2 < 1$. Then 
 $$
 x \mapsto p.v.\int^\de_{-\de} \,\frac{1}{x-u}\,\phi(u)\, du \in H^1_{\ga + 1}(I_{\de/2}).
 $$
\end{lemma}

\begin{proof}
Assume without loss of generality that $0 < x <\de/2$. Setting $\tau := x- u$, we can write 
$$
\aligned
\pa_x\Big(\int^{x +\de}_{\epsilon} \,\frac{\phi(x - \tau)}{\tau}\,d\tau\Big) &= - \int^{x +\de}_{\epsilon} \,\frac{\pa_\tau\phi(x - \tau)}{\tau}\,d\tau + \frac{\phi(-\de)}{x+\de}\\
 &= \int^{x +\de}_{\epsilon} \, \Big(\frac{1}{x} - \frac{1}{\tau}\Big)\pa_\tau\phi(x - \tau)d\tau + \Big(\frac{1}{x + \de} - \frac{1}{x}\Big)\phi(-\de) -\frac{\phi(x - \epsilon)}{x}
\endaligned
$$
and similarly for the part of the integral over $(-\epsilon, x - \de)$. Letting $\epsilon\rightarrow 0$, we conclude (after a variable change) 
$$
\aligned
 \pa_x\bigg( p.v. \int^\de_{-\de} \,\frac{1}{x-u}\,\phi(u)du\bigg) = \frac{1}{x}\, p.v.\int^\de_{-\de} \,\frac{1}{x-u}\,u\phi'(u)du - \frac{\de}{x}\bigg(\frac{\phi(\de)}{x-\de} - \frac{\phi(-\de)}{x+\de}\bigg).
\endaligned
 $$
 

\end{proof}

\begin{lemma}\label{lem:derFphi}
Let $\phi\in H^{k + 1}_{\beta + k}(I_\de)$ with $0 < \beta + 1/2 < 1$ for $k\geq 1$ and let $f\in \Lcal_{2, \ga}(I_\de)$ with $0< \ga + 1/2 < 1$. Then 
$$
[\phi, H] f \in H^{k}_{(\ga - \lambda) + k}(I_\de)
$$
where $\lambda:= 1 - (\beta + 1/2)$. The same is true, if we assume $\phi^{(j)} = O(m^{\lambda - j})$ for $1\leq j\leq k + 1$, with $\lambda < 1$. If $1 < \lambda < 2$, then 
$$
[\phi, H] f \in H^{k}_{(\ga - \lambda) + k}(I_\de)
$$ 
provided $0<(\ga - \lambda) + 1/2 < 1$. 
\end{lemma}

\begin{proof} 
It is enough to prove
$$
I(x) := \int_{-\de}^{\de} F(\phi)(x,u) f(u) du \in H^{1}_{(\ga - \lambda) + 1}(I_{\de}).
$$
Recall that, by Lemma \ref{lem:sobolev}, we have
$$
\phi' \in H^1_{\beta + 1}(I_\de) \quad \Rightarrow \quad |\phi'(x)| \, \lesssim \, m(x)^{\lambda - 1}
$$
where $\lambda := 1 - (\beta + 1/2)$. In order to estimate $\pa_x F(\phi)(x,u)$, we first consider the formula
\be\label{eq:Fj2B}
\pa_x F(\phi)(x,u) =  \frac{1}{(x - u)}\, \phi'(x)  - \frac{1}{(x - u)^2}\int_{x}^{u} \phi'(\tau)\, d\tau
\ee
(since $\phi'' \in \Lcal_{2, \beta + 1}(I_\de)$ is not integrable). Without loss of generality, assume $x>0$. We have  
$$
|\pa_x F(\phi)(x,u)| \,\lesssim\, \frac{1}{x^2}\,|\phi'(x)| + \frac{1}{x}\int^x_{-x}|\phi'(\tau)|d\tau \, \lesssim \, x^{\lambda - 2}, \quad u\in I_l(x),
$$
and therefore
$$
\int^\de_{0}x^{2(\ga - \lambda + 1)}\Big(\int_{-x}^{x/2} \pa_x F(\phi)(x,u) f(u)du \Big)^2 dx \,\lesssim\, \int^\de_{0}x^{2(\ga - 1)}\Big(\int_0^x |f(u)| + |f(-u)|\, du \Big)^2 dx
$$
which is clearly bounded by the $\Lcal_{2, \ga}$-norm of $f$ by Hardy's inequality. On the other hand, we have 
$$
|\pa_x F(\phi)(x,u)| \,\lesssim\, \frac{1}{u}\,\bigg(|\phi'(x)| + \frac{1}{u}\int^u_x|\phi'(\tau)|d\tau\bigg) \,\lesssim\, \frac{1}{u}\,\bigg(|\phi'(x)| + \int^\de_x\frac{|\phi'(\tau)|}{\tau}d\tau\bigg), \quad u\in I_r(x), \, u>0, 
$$
respectively 
$$
\aligned
|\pa_xF(\phi)(x,u)| \,\lesssim\, \frac{1}{u}\,\bigg(|\phi'(x)| + \frac{1}{x}\int^x_{-x}|\phi'(\tau)|d\tau + \frac{1}{|u|}\int^{-x}_{u}|\phi'(\tau)|d\tau\bigg)\\
\,\lesssim\, \frac{1}{u}\,\bigg(|\phi''(x)| + \frac{1}{x}\int^x_{-x}|\phi'(\tau)|d\tau + \int^{\de}_{x}\frac{|\phi'(-\tau)|}{\tau}d\tau\bigg)
\endaligned
,\quad u\in I_r(x), \, u<0.
$$
In particular, we conclude 
$$
|\pa_x F(\phi)(x,u)| \,\lesssim\,  \frac{x^{\lambda - 1}}{u},  \quad u\in I_r(x)
$$
and therefore 
$$
\int^\de_{0}x^{2(\ga - \lambda + 1)}\Big(\int_{I_r(x)} \pa_x F(\phi)(x,u) f(u)du \Big)^2 dx \,\lesssim\, \int^\de_{0}x^{2 \ga}\Big(\int_x^\de \frac{|f(u)| + |f(-u)|}{u}du \Big)^2 dx,
$$
which is bounded by the $\Lcal_{2, \ga}$-norm of $f$ by Hardy's inequality. 

It remains to consider the region where $u\in I_c(x)$. In this case, we  use the formula
\be\label{eq:Fj2A}
\pa_x F(\phi)(x,u) = \frac{1}{(x - u)^{2}}\int_x^{u} \phi''(\tau)\, \frac{(u - \tau)}{2}d\tau.
\ee
If we only control the weighted $L^2$-norm of $\phi ''$, in both cases we have the estimate  
$$
|\pa_x F(\phi)(x,u)| \,\lesssim\,  x^{\lambda - 1} \, \frac{x^{-1/2}}{|x - u|^{1/2}} \, \|\phi '' \|^2_{2, \beta + 1}.
$$
The corresponding integral is easily seen to be bounded in $\Lcal_{2, \ga - \lambda + 1}(I_\de)$ by Theorem \ref{thm:rieszintegral} below. If, on the other hand, we have $\phi''(x) = O(x^{\lambda - 2})$, then
$$
|\pa_x F(\phi)(x,u)| \,\lesssim\,  x^{\lambda - 2}
$$
in which case the claim follows by Hardy's inequality.
 
Let us briefly comment on the case $1 < \lambda < 2$. In this case $\phi''$ is integrable and we can use formula \eqref{eq:Fj2A} throughout. Then, it is not difficult to see that 
\be\label{eq:Fphiderint}
|\pa_x F(\phi)(x,u)|\,\lesssim\,  x^{\lambda - 2}, \quad u\in I_\de\setminus I_r(x), \qquad \qquad  |\pa_x F(\phi)(x,u)| \,\lesssim\,  |u|^{\lambda - 2},  \quad u\in I_r(x)
\ee  
and, by assumption $\Lcal_{2, \ga}(I_\de)\subseteq \Lcal_{2, \ga'}(I_\de)$ and $(\ga' - \lambda) + 1/2>0$.

For higher $k$, we proceed similarly, we use formula \eqref{eq:FTaylor} when $u\in I_c(x)$ and we consider derivatives of \eqref{eq:Fj2B} otherwise.
\end{proof}

\begin{lemma}\label{lem:derFphiI}
Let $\phi''\in H^{1}_{\beta + 1}(I_\de)$ with $|\phi'| \sim 1$ and $0 < \beta + 1/2 < 1$. When $f\in \Lcal_{2, \ga}(I_\de)$ with $0 < \ga + 1/2 < 1$, we have
$$
 x \mapsto \int_{-\de}^\de  \frac{\pa_x F(\phi)(x, u)}{F(\phi)(x,u)}\, f(u)du  \in H^1_{\ga + 1}(I_{\de/2}).
$$ 
\end{lemma}

\begin{proof}
The assumption $|\phi'|\sim 1$, implies $1/F(\phi) = O(1)$ and the claim basically follows as in Lemma \ref{lem:derFphi}. Indeed, we have 
$$
\pa_x\bigg(\frac{\pa_xF(\phi)}{F(\phi)}\bigg)  = \frac{\pa_x^2F(\phi)}{F(\phi)} - \bigg(\frac{\pa_xF(\phi)}{F(\phi)}\bigg)^2  
$$
where, in analogy to the proof of Lemma \ref{lem:derFphi}, we use \eqref{eq:FTaylor} with $k=2$ when $u\in I_c(x)$ and 
$$
\pa_x^2 F(\phi)(x,u) = - \frac{1}{(u - x)} \frac{\phi''(x)}{2} + \frac{1}{(u - x)^3}\int_{x}^{u} \phi''(\tau)\, \frac{(u - \tau)}{2}\, d\tau
$$
otherwise and we then retrace the steps of the previous Lemma, paying attention that $(\pa_xF(\phi))^2$ satisfies the required estimates. In fact, $\pa_x F(\phi)$ satisfies \eqref{eq:Fphiderint} with $\lambda = 1 - (\beta + 1/2)$. We omit further details.
\end{proof}

\subsubsection*{Estimates involving the fractional Laplacian}

We first state a classical result on the boundedness of the Riesz potential on weighted Lebesgue spaces with power weights, the proof of which can be found in \cite{SW}:

\begin{theorem}\label{thm:rieszintegral} Let $0 < \ga < 1/2$, i.e. both $\ga$ and $\ga - \frac{1}{2}$ give rise to Muckenhaupt weights. Then
$$
\int_{-\infty}^\infty \Big(\int_{-\infty}^\infty \frac{f(u)}{|x-u|^{1/2}}du\Big)^2 |x|^{2(\ga -1/2)}dx \, \lesssim \, \int_{-\infty}^\infty |f(x)|^2 |x|^{2\ga} dx.
$$
\end{theorem}

\begin{lemma}\label{lem:halfDer}
Assume $0 < \ga < 1/2$. Then, we have
\be\label{ineq:halfd2}
\|\Lambda^{1/2}f\|_{2,\ga - \frac{1}{2}}\,\lesssim \, \|f'\|_{2,\ga}.
\ee
\end{lemma}

\begin{proof}
It is enough to show 
$$
x \mapsto  p.v.\int_{-\de}^\de \frac{f(x) - f(u)}{|x - u|^{3/2}}\,du \in \Lcal_{2,\ga-\frac{1}{2}}(I_{\de/2}) 
$$
for some $\de > 0$ whenever $f\in H^1_{\ga}(I_{\de})$. Since we are considering weighted Lebesgue spaces, for a fixed $0<x<\de/2$, it will be more convenient to  define the principal value as
$$
p.v. \int \frac{f(x) - f(u)}{|x - u|^{3/2}}\,du = \lim_{\epsilon\rightarrow 0} \int_{|u-x|>\epsilon x} \frac{f(x) - f(u)}{|x - u|^{3/2}}\,du.
$$
Integration by parts leads to
$$
\aligned
2\int^{\de}_{(1+\epsilon)x} \frac{f(x) - f(u)}{|x - u|^{3/2}}\,du &= \int^{\de}_{(1+\epsilon)x} \big(f(u) - f(x)\big)\pa_u (u-x)^{-1/2}\,du \\
& = \frac{f(u) - f(x)}{(u - x)^{1/2}} \bigg|^{\de}_{(1+\epsilon)x} -\int^{\de}_{(1+\epsilon)x}\frac{f'(u)}{|x - u|^{1/2}}\,du \\  
\endaligned
$$
where 
$$
\frac{f(u) - f(x)}{(u - x)^{1/2}} \bigg|^{\de}_{(1+\epsilon)x} = \frac{f(\de) - f(x)}{(\de - x)^{1/2}} - \frac{f((1+\epsilon)x) - f(x)}{\sqrt{\epsilon x}}. 
$$
Since $f'\in \Lcal_{2,\ga}(I_\de)$ with $\ga$ Muckenhaupt, we have
$$
\aligned
\bigg|\frac{f((1+\epsilon)x) - f(x)}{\sqrt{\epsilon x}}\bigg| &\, \lesssim \, \frac{1}{\sqrt{\epsilon x}}\int_x^{(1+\epsilon)x}|f'(s)|ds \, \lesssim \, \|f'\|_{2,\ga} \, \bigg(\frac{1}{\epsilon}\int_x^{(1+\epsilon)x}m(s)^{-2\ga}ds\bigg)^{1/2}\\
&\, \lesssim \,   x^{-\ga}\|f'\|_{\Lcal_{2,\ga}((x, (1+\epsilon)x))}
\endaligned
$$
which vanishes as $\epsilon\rightarrow 0$. We proceed analogously on $[-\de, x(1-\epsilon))$. In particular, we obtain
$$
\bigg| p.v. \int^\de_{-\de} \frac{f(x) - f(u)}{|x - u|^{3/2}}\,du \,\bigg| \,\lesssim \, \big(\, x^{-\ga + \frac{1}{2}} + 1 \big)\, \|f'\|_{2,\ga} + \int^\de_{-\de}\frac{|f'(u)|}{|x-u|^{1/2}}\, du
$$
and the claim follows from Theorem \ref{thm:rieszintegral}. 
\end{proof}

In the next two Lemmas, we prove some estimates on the commutators with the fractional Laplacian.

\begin{lemma}\label{lem:comm_gb}
Let $0 < \ga + 1/2 < 1$ and let $\lambda \leq 1$. Then
$$
f \in \Lcal_{2,\ga-1/2}(I_{2\de}), \quad g = O(1), \quad g' = O(m^{-\lambda})  \quad \Rightarrow \quad [\Lambda^{1/2}, g]f\in  \Lcal_{2, \ga}(I_{\de}).
$$
\end{lemma}

\begin{proof} 
Assume without loss of generality that $x>0$. We have the estimates  
$$
\aligned
&\Big|\int_{I_l(x)} \frac{g(x) - g(u)}{|x -u|^{3/2}} \, f(u) d u\Big|   \, \lesssim \,  x^{- 3/2} \int_{-x}^{x} |f(u)| du, \\
&\Big|\int_{I_r(x)} \frac{g(x) - g(u)}{|x -u|^{3/2}} f(u) d u\Big| \, \lesssim \, \int_{x}^{\de} \frac{|f(u)| + |f(-u)|}{u^{3/2}}\, du. 
\endaligned
$$
Both are bounded in the desired space by Hardy inequalities. When $u\in I_c(x)$ (and we assume e.g. $-1/2< \ga < 0$) we have
$$
\Big|\, p.v\int_{I_c(x)} \frac{g(x) - g(u)}{|x -u|^{3/2}} \, f(u) d u\Big| \, \lesssim \, \int_{I_c(x)} \frac{|f(u)|u^{-1}}{|x - u|^{1/2}} \, d u,
$$
where we have used $g'(x) =O( x^{-\lambda} ) = O( x^{-1})$ since by assumption $\lambda \leq 1$. The result now follows from Theorem \ref{thm:rieszintegral} (we set $f$ identically equal to zero on the complement of say $(0, \de]$). In the case $0 <\ga< 1/2$ replace $x^\ga u^{-1}$ by $x^{\ga - 1/2}u^{-1/2}$.
\end{proof}

\begin{lemma}\label{lem:comm_weight}
Let $0 < \ga + 1/2 < 1$. Then, 
$$
f\in \Lcal_{2,\ga}(I_{2\de}) \quad \Rightarrow \quad [\Lambda^{1/2}, m^{2\ga}]f\in  \Lcal_{2, 1/2-\ga}(I_\de)
$$ 
\end{lemma}

\begin{proof} Without loss of generality, assume $x>0$ and let $u\in I_l(x)$ first. Then 
$$
x^{1/2-\ga}\bigg|\int_{I_l(x)} \frac{|x|^{2\ga} - |u|^{2\ga}}{|x -u|^{3/2}} \, f(u) d u\bigg| \, \lesssim \,  x^{\ga - 1}\int_{-x}^{x} |f(u)| d u + x^{-\ga-1}\int_{-x}^{x} |u|^{2\ga}|f(u)| d u
$$
with both terms square integrable over $(0, \de)$ by Hardy's inequality. Note that $-\ga$ satisfies the Muckenhaupt condition as well. 

On the other hand, when $u\in I_r(x)$, we have
$$
\aligned
x^{1/2-\ga}\bigg|\int_{I_r(x)} \frac{|x|^{2\ga} - |u|^{2\ga}}{|x -u|^{3/2}} f(u) d u\bigg| \, \lesssim \, x^{\ga}&\int_{x}^{\de} \frac{|f(u)| + |f(-u)|}{u}\, du \\
&+ x^{-\ga}\int_{x}^{\de} u^{2\ga}\frac{|f(u)| + |f(-u)|}{u}\, du.
\endaligned
$$
Again, both are square integrable over $(0, \de)$ by Hardy's inequality.

Finally, when $u\in I_c(x)$ (and e.g. $-1/2 < \ga < 0$, i.e. both $\ga$ and $\ga + 1/2$ satisfy the Muckenhaupt condition), we have
$$
x^{1/2 - \ga}\Big|\, p.v\int_{I_c(x)} \frac{|x|^{2\ga} - |u|^{2\ga}}{|x -u|^{3/2}} \, f(u) d u\Big| \, \lesssim \, x^{\ga}\int_{I_c(x)} \frac{1}{|x - u|^{1/2}} \, |f(u)| u^{-1/2}d u,
$$
which is square integrable over $(0, \de)$ by Theorem \ref{thm:rieszintegral}. The same is true if $0< \ga < 1/2$, however in that case we do not need the extra $u^{-1/2}$ in the integral. 
 
\end{proof}

\end{subsection}

\subsection{Miscellaneous}

For completeness we state the general form of Hardy inequalities (see e.g. \cite{KufOpic}):

\begin{theorem}\label{thm:hardy}
	Let $-\infty \leq a < b \leq \infty$ and let $m_1, m_2$ be measurable, positive functions on $(a,b)$. Then, we have 
	\be\label{eq:hardyI}
	\int_a^b\Big(\int^x_a f(s)ds\Big)^2m_1(x)^2dx \,\leq \, C \int_a^b|f(x)|^2 m_2(x)^2dx
	\ee
	if and only if 
	$$
	\sup_{a<x<b}\Big(\int_x^b m_1(s)^2ds\Big)\Big(\int^x_a m_2(s)^{-2}ds\Big) \, < \, \infty.
	$$
	Similarly, we have the dual inequality 
	\be\label{eq:hardyII}
\int_a^b\Big(\int^b_x f(s)ds\Big)^2m_1(x)^2dx  \,\leq \, C \int_a^b|f(x)|^2 m_2(x)^2dx
	\ee
	if and only if
	$$
	\sup_{a<x<b}\Big(\int_a^x m_1(s)^2ds\Big)\Big(\int^b_x m_2(s)^{-2}ds\Big) \, < \, \infty.
	$$	
\end{theorem}

It is then not difficult to verify that for power weights we have:

\begin{lemma}\label{lem:hardy}
	Let $p>1$ and let $I_\de:= (0,\de)$, with $\de < \infty$.
	\begin{enumerate}
		\item\label{itm:hardy1} If $\ga + 1/2 < 1$, then  
		$$
		f\mapsto \int_0^x f(t)dt \,:\, \Lcal_{2,\beta}(I_\de) \longrightarrow \Lcal_{2,\ga}(I_\de) 
		$$
		is continuous for all $\beta \geq \ga -1$. 
		
		\item\label{itm:hardy2} If $\ga + 1/2 >0$, then 
		$$
		f\mapsto \int_x^\de f(t)dt \,:\, \Lcal_{2,\ga}(I_\de) \longrightarrow \Lcal_{2,\beta}(I_\de)
		$$
		is continuous for all $\ga \leq \beta + 1$. 
	\end{enumerate}
\end{lemma}

We will need the following version of the Lax-Milgram Lemma, which can be found in \cite{Kufner}:

\begin{lemma}\label{lem:weightedLaxMilgram}
 Let $H_1, H_2$ be two Hilbert spaces. Let $a(u,v)$ be a bilinear form defined on the Cartesian product $H_1\times H_2$, and let there exist positive constants $c_1$, $c_2$, $c_3$ such that 
 \begin{enumerate}[leftmargin=*, label=(\roman*)]
\item for all $u\in H_1$ and $v\in H_2$ we have 
\be\label{eq:acontinuity}
|a(u,v)| \, \leq \, c_1 \|u\|_{H_1} \|v\|_{H_2};
\ee
\item for all $u\in H_1$ we have 
\be\label{eq:acoercivity1}
\sup_{\|v\|_{H_2}\leq 1}|a(u,v)| \, \geq \, c_2 \|u\|_{H_1} ;
\ee
\item for all $v\in H_2$ we have 
\be\label{eq:acoercivity2}
\sup_{\|u\|_{H_1}\leq 1}|a(u,v)| \, \geq \, c_3 \|v\|_{H_2}.
\ee
\end{enumerate}
 Let $f$ be a continuous linear functional on $H_2$. Then there exists precisely one element $u\in H_1$ such that
 $$
 a(u, v) = \langle f, v\rangle, \quad \forall v\in H_2
 $$
and there is a positive constant $c$ such that 
$$
\|u\|_{H_1} \, \leq \, c\|f\|_{H_2^*}.
$$
\end{lemma}

We state the following result on pseudodifferential operators which can be found in e.g. \cite[Proposition 9.10]{TaylorIII}: 

\begin{theorem}\label{thm:PsiDO} Let $r\in(0,\infty)$ and let $a(x,\xi)\in \Ccal^rS^0$, that is
	\be\label{eq:a-growth}
	\|\pa_\xi^k a(\cdot, \xi) \|_{\Ccal^r(\Rbb)} \, \leq \, C_k\langle\xi\rangle^{-k}, \quad \forall k\geq0. 
	\ee
	Then, $a(x, D) : W^{s,p}(\Rbb) \rightarrow W^{s,p}(\Rbb)$ for all $|s|<r$ and all $p\in(1,\infty)$, where 
	$$
	a(x,D)g = \int_{\Rbb} e^{ix\xi} a(x,\xi) \h g(\xi) d\xi 
	$$
\end{theorem}

\subsection{Variable-step convolution operator}
\label{a.convolution}

To finish, we discuss some properties of the variable-step convolution operators $B_\de, B_\de^*$ and $A^\de$ cf.\eqref{E.defphide}--\eqref{E.defAde} that were used in Sections \ref{ss.regularization} and \ref{ss.existence}.

Recall that for $\de > 0$, we have
$$
B_\de f(x) := (\phi_{\de \eta(x)} \ast f)(x) := \int_\Rbb  \frac{1}{\de \eta(x)}\phi\left(\frac{x-y}{\de \eta(x)}\right) f(y)dy
$$
where $\phi$ is (periodically extended) standard mollifier and $\eta$ is a non-negative, periodic function which we define as follows. We first extend $\al \rightarrow ||\al| - \al_*|$ for $\al \in [-\pi,\pi]$ periodically to all of $\Rbb$, then we smooth it out in a small neighborhood of $n\pi$ for $n\in\Zbb$ (in such a way that $\eta(\al) \leq ||\al| - \al_*|$ when e.g. $\al\in[-\pi,\pi]$). In particular, $\eta$ is smooth on $\Rbb\setminus\{\al_* + n\pi: n\in\Zbb\}$ and we have $\eta \sim m$. See \cite{HPR} and references therein for a more general discussion on these type of operators. 

Since $\Rbb\setminus\{\al_* + n\pi, n\in\Zbb\}$ is a countable union of disjoint open intervals $I_n$ and that given $x\in I_n$ we actually integrate over a subinterval of $I_n$, it is enough to fix one, say $I := (-\al_*, \al_*)$.    

Indeed, given $x\in I$, the interval of integration $B_{\de\eta(x)}(x) = (x-\de \eta(x), x + \de \eta(x))$ is contained in $I$ and has positive distance to $\pa I$ for all sufficiently small $\de > 0$. In particular, the integral is well-defined for $f\in L^1_{loc}(I)$ and $m(x) \sim m(y)$ uniformly for $y\in B_{\de\eta(x)}(x)$.

Similarly, in the case of the adjoint, 
$$
B_{\de}^*g(y) := \int \phi_{\de \eta(x)}(x - y)g(x)dx,
$$
for fixed $y\in I$, the integral runs over $[x^-(y), x^+(y)]\subseteq I$, where $x^\pm(y)$ are the solutions of $x \mp \de \eta(x) = y$, respectively. These are well defined for all $\de<1 / \sup |\eta'| $. It is not difficult to see that  $x^\pm(y) \rightarrow \pa I$ as $y\rightarrow \pa I$ (and they have the same limit).

In the next Lemma, we give properties of $B_\de$, $B_\de^*$ and $A^\de$ for functions defined in weighted Sobolev spaces $\Lcal^k_{2,\ga}(I)$ with weight given by the restriction $m|_I$. The generalization to the periodic setting is straightforward.

\begin{lemma}\label{lem:B_delta} 
	For any $j\geq0$, all~$\ga\in\Rbb $ and all $\de>0$, the operator $B_\de : \Lcal^j_{2, \ga + j}(I) \to  \Lcal^j_{2, \ga + j}(I)$ is continuous. Furthermore, it preserves growth rates near~$\pa I$ in the sense that, for any $f\in \Lcal^j_{2, \ga + j}(I)$,
	$$
	B_\de f\in \Ccont^\infty(I) \cap \Lcal^{j'}_{2,\ga' + j'}(I)
	$$
	for all $\ga'\in\Rbb$ and all $j'\in\Nbb $, and it approximates~$f$ as
	$$
	\lim_{\de\rightarrow 0}\|B_{\de} f -  f\|_{\Lcal^j_{2, \ga + j}(I)} = 0.
	$$
	If $f\in \Lcal^1_{2,\ga + 1}(I)$, then
	$$
	(B_\de f)' = B_\de(f') + K_\de f,
	$$
	where $K_\de: \Lcal_{2,\ga}(I) \rightarrow \Lcal_{2,\ga + 1}(I)$ is bounded, and one has the quantitative estimate
	\be\label{eq: BdeII}
	\|B_{\de} f -  f\|_{\Lcal_{2, \ga}(I)} \ \lesssim \ \de^{1/2} \, \|f\|_{\Lcal^{1}_{2,\ga + 1}(I)}.
	\ee
	Analogous estimates hold for the case of higher~$j$. Furthermore, $B_\de^*$ and~$A_\de$ enjoy the same properties.

	Finally, if~$g$ is a function bounded as $g' = O(m^{\lambda})$ for some $\lambda \in\Rbb$, then  
	$$
	[g, \, B_\de] f \in  \Lcal^{1}_{2, \ga - \lambda}(m), \qquad [g, \, B_\de^*] f \in  \Lcal^{1}_{2, \ga - \lambda}(m).
	$$	
\end{lemma}

\begin{proof} 
	We consider the case $j=0$, that is, $f\in\Lcal_{2,\ga}(I)$. 	The estimates for higher $j$ follow from very similar arguments.
	
	To  show that $B_\de (f)\in\Lcal_{2, \ga}(m)$, we write 
	$$
	\aligned
	|B_\de f(x)|^2 &\, \lesssim \, \int_{B_{\de\eta(x)}(x)} \phi_{\de\eta(x)}(x-y) m(y)^{-2\ga} dy \int_{B_{\de\eta(x)}(x)} \phi_{\de\eta(x)}(x-y) |f(y)|^2 m(y)^{2\ga} dy\\
	&\,\lesssim \, m(x)^{-2\ga}\int_{B_\de(x)} \phi_{\de\eta(x)}(x-y) |f(y)|^2 m(y)^{2\ga} dy, 
	\endaligned
	$$
	where we have used that $m(x) \sim m(y)$ uniformly for $|y-x|\leq \de\eta(x)$. Multiplying both sides by $m(x)^{2\ga}$ and integrating, and then using that $ B^*_\de (1) \lesssim  1$, we conclude that $B_\de$ is bounded in~$\Lcal_{2, \ga}(m)$. Since $\Ccont_c(I)$ is dense in $\Lcal_{2, \ga}(I)$ (see e.g.~\cite{Kufner}), we have $B_\de(f) \rightarrow f$ in $\Lcal_{2, \ga}(I)$ as $\de\to 0$. Although we shall not need this fact, for $f\in\Ccont(\overline{I})$, we actually have $\sup |B_\de f - f| \rightarrow 0$ as $\de \rightarrow 0$. It is not difficult to see that $B_\de f$ are smooth in $I$ and that taking $j$ derivatives on $\phi$ leads to factors bounded as $O((\de m)^{-j})$. All these facts hold for $B_\de^*$ and $A_\de$ as well.
	
	Let us now consider the commutator with a derivative and $B_\de$ (or $B_\de^*$). With $f\in \Lcal^1_{2,\ga + 1}(I)$, we write
	$$
	\aligned
	(B_\de f)'(x) &= \pa_x \int  \phi\left(v\right) f\left(x+\de v\eta(x)\right)dv \\
	&=\int  \phi\left(v\right) f'\left(x+\de v\eta(x)\right)(1+\de v \eta'(x))dv \\  
	&= \int \phi\left(v\right) f'\left(x+\de v\eta(x)\right)dv - \frac{\eta'(x)}{\eta(x)}\int \pa_v (v\phi\left(v\right))  f\left(x+\de v\eta(x)\right) dv  \\ 
	&=: B_\de(f') +  K_\de f,
	\endaligned
	$$ 
	where $K_\de$ is bounded in $\Lcal_{2,\ga }(I)\to \Lcal_{2,\ga + 1}(I)$ by the same argument as above. Similarly, 
	$$
	\aligned
	B_\de^* (g)'(y) 
	&=  \int  \phi\left(v\right) g\left(x(y, v)\right)\pa_y x \, dv \\
	&= \int  \phi\left(v\right) \left(g'(x(y, v))(\pa_y x)^2 + g(x(y, v))\pa_y^2 x\right) dv \\
	&= \int  \phi\left(v\right) g'(x(y, v))\pa_y x dv + \int \left(\pa_v\left(\phi\left(v\right)\frac{\pa_y x}{\pa_v x}(\pa_y x - 1)\right) + \phi(v)\pa_y^2 x\right) g(x(y, v)) dv \\
	&=: B_\de^* (g') + \t K_\de(g),
	\endaligned
	$$
	where note that $\frac{\pa_y x}{\pa_v x}(\pa_y x - 1) = \frac{\eta'(x)}{\eta(x)}v\pa_y x$. The estimate for the commutator of $A_\de$ with the derivative follows directly from these results and the formula
	\be\label{E.commAdepa}
	[A_\de, \pa] = [B_\de^*, \pa]B_\de  + B_\de^*[B_\de, \pa].
	\ee
	
	Finally, to prove the claim for the commutator with $g$, we use the formula 
	$$
	\pa_x  [ g, \, B_\de](f)(x) = \int_{B_{\de\eta(x)}(x)} \pa_x\phi_{\de\eta(x)}(x- y)(g(x) - g(y)) f(y)dy + g'(x)B_\de(f)(x).
	$$
	As we are integrating over a set where $m(x) \sim m(y)$, the claim trivially follows from the fact that
	$$
	(x - y)\pa_x\phi_{\de\eta(x)}(x- y) = O(1), \qquad \frac{|g(x) - g(y)|}{|x - y|} = O(m(x)^{\lambda}).
	$$
	The lemma is then proven.
	
\end{proof}

\end{appendices}



\end{document}